\documentclass[12pt]{elsarticle}
\usepackage{graphicx}
\usepackage{subcaption} 
\usepackage{epstopdf} 
\usepackage{subcaption} 
\usepackage{geometry} 
\usepackage{hhline}

\usepackage{booktabs}
\usepackage{multirow}
\usepackage{pdflscape}
\usepackage[table,xcdraw]{xcolor}
\usepackage{ucs}
\usepackage{latexsym}
\usepackage{setspace}
\usepackage{color}
\usepackage{amsmath}
\usepackage{amssymb}
\usepackage{latexsym}
\usepackage{amsmath}
\usepackage{mathrsfs}
\usepackage{amssymb}
\usepackage{natbib}
\setcitestyle{authoryear,open={(},close={)}}
\setcitestyle{numbers}
\setcitestyle{square}
\usepackage{hyperref}
\hypersetup{colorlinks,citecolor=blue}
\usepackage{graphics}
\usepackage{graphicx}
\makeatletter
\def\ps@pprintTitle{%
	\let\@oddhead\@empty
	\let\@evenhead\@empty
	\let\@oddfoot\@empty
	\let\@evenfoot\@oddfoot
}
\makeatother
\makeatletter \oddsidemargin  -.5cm \evensidemargin -.5cm \textwidth
17.5cm \topmargin -.65in \textheight 24cm

\newcommand{\bd}{\begin{definition}}
	\newcommand{\ed}{\end{definition}}
\newcommand{\br}{\begin{remark}}
	\newcommand{\er}{\end{remark}}
\newcommand{\bea}{\begin{eqnarray}}
	\newcommand{\eea}{\end{eqnarray}}
\newcommand{\beann}{\begin{eqnarray*}}
	\newcommand{\eeann}{\end{eqnarray*}}
\newtheorem{theorem}{Theorem}[section]

\newtheorem{lemma}[theorem]{Lemma}

\newtheorem{corollary}[theorem]{Corollary}
\newtheorem{remark}{Remark}[section]
\newtheorem{example}{Example}[section]

\numberwithin{equation}{section}
\numberwithin{equation}{section}
\title{Improved estimation of the positive powers ordered restricted standard deviation of two normal populations}
\author{Somnath Mondal  and  Lakshmi Kanta Patra\footnote
	{\baselineskip=10pt
		~lkpatra@iitbhilai.ac.in; patralakshmi@gmail.com}  \\
	\it  Department of Mathematics\\
	\it Indian Institute of Technology Bhilai, Durg, India-491001}
\spacing{1.1}
\begin{document}
	\begin{frontmatter}
		\date{}
		\begin{abstract}
		The present manuscript is concerned with component-wise estimation of the positive power of the ordered restricted standard deviation of two normal populations with certain restrictions on the means. We have obtained sufficient conditions to prove the dominance of equivariant estimators with respect to a general scale-invariant bowl-shaped loss function. Consequently, we propose various estimators that dominate the best affine equivariant estimator (BAEE). Also, we obtained a class of improved estimators and proved that the boundary estimator of this class is generalized Bayes. The improved estimators are derived for four special loss functions: quadratic loss, entropy loss, symmetric loss, and Linex loss function. We have conducted extensive Monte Carlo simulations to study and compare the risk performance of the proposed estimators. Finally, we have given a data analysis for implementation purposes. 
		
		\noindent\textbf{Keywords}: Decision theory; Improved estimators; Scale invariant loss function; Generalized Bayes; Relative risk improvement.
	\end{abstract}
\end{frontmatter}

\section{Introduction} 
The problem of estimating parameters under order restriction has received significant attention due to its practical applications across various fields, including bio-assays, economics, reliability, and life-testing studies. For instance, ranking employee pay based on their job description is reasonable. It is anticipated in agricultural research that the average yield of a particular crop will be higher when fertilizer is used than when it is not. Suppose we measure voltage using two voltmeters. One is an old version, and the other one is updated. In this case, it is reasonable to assume that the variability in the measurements taken by the old version is higher than that of the updated one. Also, voltages are usually positive, we can take the mean of these measurements as positive.
Thus, imposing an order restriction on some model parameters, such as average values and variance, makes sense. Estimators are more efficient when this prior knowledge of the order restriction on parameters is considered. 
Some early works in this direction are \cite{MR0326887}, \cite{MR961262} and \cite{MR2265239}. The problem of finding improved estimators of ordered parameters in various probability distributions has been extensively studied in the literature. For some important contributions in these directions are  \cite{kumar1988simultaneous},  \cite{misra1994estimation},  \cite{vijayasree1995componentwise}, \cite{misra1997estimation}, \cite{chang2002comparison},  \cite{misra2002smooth}, \cite{oono2005estimation}, \cite{oono2006class}, \cite{jana2015estimation}. The authors have used the approach of \cite{stein1964}, \cite{brewster1974improving}, and \cite{kubokawa1994unified} to derive the estimators that dominate usual estimators such as the maximum likelihood estimator (MLE), best affine equivariant estimator  (BAEE), etc. 

\cite{tripathy2013estimating} has considered estimating the common variance of two normal distributions with ordered location parameters under the quadratic loss function. They have shown that the usual estimators are inadmissible by proposing improved estimators. Estimation of ordered restricted normal means under a Linex loss function has been studied by  \cite{ma2013estimation}. The authors prove that plug-in estimators improve upon the unrestricted MLE. \cite{petropoulos2017estimation} studied component-wise estimation of ordered scale parameter two Lomax distributions with respect to the quadratic loss function. He proposed various estimators that dominate the BAEE. \cite{chang2017estimation} studied the estimation of two ordered normal means with a known covariance matrix using the Pitman nearness criterion. \cite{patra2021componentwise} discussed the component-wise estimation of the ordered scale parameter of two exponential distributions with respect to a general scale-invariant loss function. They have proved the inadmissibility of usual estimators by proposing several improved estimators.  
\cite{jana2022estimation} has considered the component-wise estimation of the ordered variance of two normal populations with common mean under a quadratic loss function. They have proposed various estimators that dominate some usual estimators. \cite{garg2023componentwise} has investigated improved estimation of ordered restricted location and scale parameters of a bivariate model with respect to a general bowl-shaped invariant loss function. They have used the techniques of \cite{brewster1974improving} to derive estimators that improve upon the usual estimators. 

In this paper, we consider the problem of estimating the positive powers of ordered scale parameters for two normal distributions. For this estimation problem, we consider a class of scale-invariant bowl-shaped loss functions $L\left(\frac{\delta}{\theta}\right)$, where $\delta$ is an estimator of $\theta$. We assume that the loss function $L(t)$ satisfies the following criteria:
\begin{enumerate}
	\item[(i)] $L(t)$ is strictly bowl shaped that is $L(t)$ decreasing for $t\leq 1$ and increasing for $t\geq 1$ and  reaching its minimum value $0$ at $t=1$.  
	\item[(ii)] The integrals involving $L(t)$ are finite and can be differentiated under integral sign.
	\item[(iii)] $L'(t)$ is increasing, almost everywhere.
\end{enumerate}
Let the random variables $X_1$, $X_2$, $S_1$ and $S_2$ are independent and distributed as,
\begin{equation}
	X_1\sim N\left(\mu_1,\frac{\sigma_1^2}{p_1}\right),~S_1 \sim \sigma_1^2\chi ^2_{p_1-1}~~\mbox{ and }~~
	X_2\sim N\left(\mu_2,\frac{\sigma_2^2}{p_2}\right),~S_2 \sim \sigma_2^2\chi ^2_{p_{2}-1}
\end{equation}
with unknown $\mu_i$, $\sigma_i$ for $i=1,2$ and $\sigma_1\leq \sigma_2$. Now onwards we denote $\underline{X}=(X_1,X_2)$, $\underline{S}=(S_1,S_2)$, $V_i=S_i/\sigma_i^2$ and $\underline{\theta}=(\sigma_1,\sigma_2,\mu_1,\mu_2)$. Here we will propose various estimators that improve upon the BAEE of $\sigma_1^k$ and $\sigma_2^k$, $k>0$. The main contributions of this article are as follows. 
	\begin{enumerate}
		\item [(i)] We have obtained BAEE of $\sigma^k_i$ with respect to a $L(t)$. We opposed Stein-type improved estimators, and as an application, we derive an estimator that dominates BAEE of $\sigma^k_i$ when there is no restriction on $\mu_1$ and $\mu_2$. Further, we have derived a class of improved estimators, and it is shown that the boundary estimator of this class is a generalized Bayes estimator for estimating $\sigma^k_i$. 
		\item [(ii)] Next we consider the improved estimation estimation of $\sigma_i^k$ (for $i = 1, 2$)  when both $\mu_1$ and $\mu_2$ are positive.  In this case, we have also proposed several estimators that dominate BAEE of $\sigma^k_i$ under $L(t)$. 
		\item [(iii)] Finally an improved estimator of $\sigma^k_i$ has been derived when $\mu_1 \le \mu_2$. As an application, we have obtained improved estimators with respect to four special loss functions: quadratic loss, entropy loss, symmetric loss, and Linex loss function. 
		\item [(iv)] A simulation study has been carried out to measure the risk performance of the proposed estimators of $\sigma^2_i$. We have plotted the relative risk improvement with respect to the BAEE of the proposed estimators to compare the risk performance. 
\end{enumerate}
We first apply the invariance principle to obtain BAEE. For this purpose, we consider the group of transformations as 
$\mathcal{G}= \left\{g_{a_1,a_2,b_1,b_2}:\ a_1>0,\ a_2>0,\ b_1 \in \mathbb{R},\ b_2 \in \mathbb{R}\right\}.$
The group $\mathcal{G}$ act on the $\mathscr{X}=\mathbb{R}\times \mathbb{R}\times\mathbb{R}_+\times \mathbb{R}_+$ in the following manner	
$$(X_1,X_2,S_1,S_2)\rightarrow(a_1X_1+b_1,a_2X_2+b_2,a_1^2S_1,a_2^2S_2).$$
Under this group of transformations, the problem of estimating $\sigma^k_i$ is invariant. After some simplification, the form of  an affine equivariant estimator is obtained as 
\begin{equation}
	\delta_{ic}(\underline{X},\underline{S})=cS_i^{\frac{k}{2}},\ i=1,2,
\end{equation}
where $c>0$ is a constant. The following lemma provides the BAEE of $\sigma_i^k$. 
\begin{lemma}
	Under a bowl-shaped loss function  $L(t)$ the best affine equivariant estimator of $\sigma_i^k$ is $\delta_{0i}\left(\underline{X},\underline{S}\right)=c_{0i}S_i^{\frac{k}{2}}$, where $c_{0i}$ is the unique solution of the equation 
	\begin{equation}\label{eq2.2}
		E\left[L'\left(c_{0i}V_i^{\frac{k}{2}}\right)V_i^{\frac{k}{2}}\right]=0.
	\end{equation}
\end{lemma}
\begin{example}\rm
	For $i=1,2$
	\begin{enumerate}
		\item [(i)] Under the quadratic loss function $L_1(t)=(t-1)^2$, the BAEE of $\sigma_i^k$ is obtained as  $\delta_{0i}^1=\frac{\Gamma\left(\frac{p_i+k-1}{2}\right)}{2^{\frac{k}{2}}\Gamma\left(\frac{p_i+2k-1}{2}\right)}S_i^{\frac{k}{2}}$.
		\item [(ii)] For the entropy loss function $L_2(t)=t-\ln t-1$, we get the BAEE of $\sigma_i^k$ is $\delta_{0i}^2=\frac{\Gamma\left(\frac{p_i-1}{2}\right)}{2^{\frac{k}{2}}\Gamma\left(\frac{p_i+k-1}{2}\right)}S_i^{\frac{k}{2}}$.
		\item [(iii)] The BAEE of  $\sigma_i^k$ is $\delta_{0i}^3=\sqrt{\frac{\Gamma\left(\frac{p_i-k-1}{2}\right)}{2^k\Gamma\left(\frac{p_i+k-1}{2}\right)}}S_i^{\frac{k}{2}}$ with respect ot a symmetric loss function $L_3(t)=t+\frac{1}{t}-2$.
		\item [(iv)] For linex loss function $L_4(t)=e^{a(t-1)}-a(t-1)-1$; $a \in \mathbb{R}-\{0\}$ the BAEE of  $\sigma_i^k$ is $\delta_{0i}^4=c_{0i}S_i^{\frac{k}{2}}$, where $c_{0i}$ is the solution of the equation 
		$$	\int_{0}^{\infty}v_i^{\frac{p_i+k-1}{2}-1}e^{-\frac{v_i}{2}+ac_{0i}v_i^{\frac{k}{2}}}dv_i=e^a\int_{0}^{\infty}v_i^{\frac{p_i+k-1}{2}-1}e^{-\frac{v_i}{2}}dv_i.$$
			In particular, for $k=2$ we have $c_{0i}=\frac{1}{2a}\left(1-e^{-\frac{2a}{1+p_i}}\right)$ and thus the BAEE of $\sigma_i^2$ is obtained as $\frac{1}{2a}\left(1-e^{-\frac{2a}{1+p_i}}\right) S_i^{\frac{k}{2}}$.
	\end{enumerate}
\end{example}
\begin{remark}
	The UMVUE of $\sigma_i^k$ is $\delta_{iMV}=\frac{\Gamma\left(\frac{p_i-1}{2}\right)}{2^{\frac{k}{2}} \Gamma\left(\frac{p_i+k-1}{2}\right)}S_i^{\frac{k}{2}}$. We observe that this is the BAEE with respect to entropy loss function $L_2(t)$. Also, the BAEE improves upon the  UMVUE under the loss function $L_1(t)$, $L_3(t)$ and $L_4(t)$. 
\end{remark}

The rest of the paper is organized as follows. In Section \ref{sec2},  we consider the estimation of $\sigma_1^k$ when $\sigma_1\le \sigma_2$. We have proposed estimators that dominate the BAEE. A class of improved estimators is obtained, and it is shown that the boundary estimator of this class is a generalized Bayes estimator. In Subsection \ref{subsec2.2}, we have considered improved estimation $\sigma_1^k$ when $\mu_1$ and $\mu_2$ are non-negative. Next, we have studied the estimation of $\sigma_1^k$ when $\mu_1 \le \mu_2$. Further, as an application, we have derived improved estimators for four special loss functions. In Section \ref{sec3}, we have obtained results similar to Section \ref{sec2} for estimating $\sigma_2^k$. A simulation has been carried out to compare the risk performance of the improved estimators in Section \ref{sec4}. Finally in Section \ref{sec5} we have presented a real life data analysis. 
\section{Improved estimation of $\sigma_1^k$ when $\sigma_1\le \sigma_2$}\label{sec2} 
In this section, we consider the problem of finding an improved estimation of $\sigma_1^k$ with the constraint $\sigma_1 \leq \sigma_2$. Similar to \cite{petropoulos2017estimation} we consider a class of estimators of the form 
\begin{equation}\label{bae2.1}
	\mathcal{C}_1=\left\{\delta_{\phi_1}=\phi_{1} \left(U\right)S_1^{\frac{k}{2}}:  U=S_2S_1^{-1}  \mbox{ and }\phi_{1}(.) \mbox{ is positive measurable function}\right\}.
\end{equation}
Now we analyse the risk function  $R\left(\underline{\theta},\delta_{\phi_1}\right)= E\left[E\left\{L\left(V_1^{k/2}\phi_1(U)\right)\big\rvert U\right\} \right]$ for $k>0$. The conditional risk function can be written as $R_1(\underline{\theta}, c)=E_{\eta}\left\{L\left(V_1^{\frac{k}{2}}c\right)\big\rvert U=u\right\}$, where $V_1|U=u \sim Gamma\left(\frac{p_1+p_2-2}{2}, \frac{2}{(1+\eta^2u)}\right)$ distribution, with $\eta=\frac{\sigma_1}{\sigma_2} \le 1$. The function $R_1(\underline{\theta},c)$ minimized at $c_{\eta}(w)$, where $c_{\eta}(u)$ be the unique solution of  $E_{\eta}\left(L'\left(V_1^{k/2}c_{\eta}(u)\right)V_1^{k/2}|U=u\right)=0$. Using Lemma $3.4.2.$ of \cite{lehmann2005testing}, we have 
\begin{align*}
	E_{\eta}\left(L'\left(V_1^{k/2}c_1(u)\right)V_1^{k/2}|U=u\right)&\geq E_1\left(L'\left(V_1^{k/2}c_1(u)\right)V_1^{k/2}|U=u\right)\\
	&=0
=E_{\eta}\left(L'\left(V_1^{k/2}c_{\eta}(u)\right)V_1^{k/2}|U=u\right).
\end{align*}
Consequently we get $c_{\eta}(u) \le c_{1}(u)$, where  $c_{1}(u)$ is the unique solution of $$E_1\left(L'\left(V_1^{k/2}c_1(u)\right)V_1^{k/2}|U=u\right)=0.$$ 
Making the transformation $z_1=v_1(1+u)$, we get  $E\left(L^{\prime}(Z_1^{k/2}c_1(u)(1+u)^{-k/2})\right)=0$ with  $Z_1 \sim \chi^2_{p_1+p_2+k-2}$. Comparing with (\ref{s1st1}), we obtain $c_1(u)=\alpha_1(1+u)^{\frac{k}{2}}$. Consider $\phi_{01}(u)=\min\{\phi_1(u), c_1(u)\}$, then for $P(c_{1}(U)<\phi_1(U)) \ne 0$ we get $c_{\eta}(u) \le c_{1}(u)=\phi_{01}(u) <\phi_1(u)$ on a set of positive probability. Hence we get  $R_{1}(\underline{\theta},\phi_{01}) < R_{1}(\underline{\theta},\phi_1)$. So we get the result as follows. 
\begin{theorem}	\label{th2.1}
Let $\alpha_1$ be a solution of the equation 
	\begin{equation}\label{s1st1}
		EL^{\prime} \left(Z_1^{k/2}\alpha_1\right)=0
	\end{equation}
	where $Z_1 \sim \chi^2_{p_1+p_2+k-2}$.
Then the risk of the estimator	$\delta_{\phi_{01}}=\phi_{01}(U)S_1^{\frac{k}{2}}$ is nowhere larger than the estimator $\delta_{\phi_1}$ provided $P(\phi_1(U)> c_{1}(U)) \ne 0$ holds true.
\end{theorem}
\begin{corollary}
	The risk of the estimator
	$\delta_{11}=\min \left\{c_{01}, \alpha_1(1+U)^{k/2}\right\} S_1^{k/2}$
	is nowhere larger than the BAEE $\delta_{01}$ provided $ \alpha_1 < c_{01}$.
\end{corollary}


\begin{example}
	\begin{enumerate}
		\item[(i)] For the quadratic loss function $L_1(t)$, we have $\alpha_{1}=\frac{\Gamma\left(\frac{p_1+p_2+k-2}{2}\right)}{2^{\frac{k}{2}}\Gamma\left(\frac{p_1+p_2+2k-2}{2}\right)}$. The improved estimator of $\sigma_1^k$ is obtained as 
		$$\delta^1_{11} = \min \left\{\frac{\Gamma\left(\frac{p_1+k-1}{2}\right)}{2^{\frac{k}{2}}\Gamma\left(\frac{p_1+2k-1}{2}\right)},\alpha_{1}(1+U)^{\frac{k}{2}}\right\} S_1^{\frac{k}{2}}.$$
		
		\item [(ii)] Under the entropy loss function $L_2(t)$, we get $\alpha_{1}=\frac{\Gamma\left(\frac{p_1+p_2-2}{2}\right)}{2^{\frac{k}{2}}\Gamma\left(\frac{p_1+p_2+k-2}{2}\right)}$. So the improved estimator is
		$$ \delta^2_{11}= \min \left\{\frac{\Gamma\left(\frac{p_1-1}{2}\right)}{2^{\frac{k}{2}}\Gamma\left(\frac{p_1+k-1}{2}\right)},\alpha_{1}(1+U)^{\frac{k}{2}}\right\} S_1^{\frac{k}{2}}.$$
	
		\item [(iii)] For the symmetric loss function $L_3(t)$ we obtain $\alpha_{1}=\sqrt{\frac{\Gamma\left(\frac{p_1+p_2-k-2}{2}\right)}{2^k\Gamma\left(\frac{p_1+p_2+k-2}{2}\right)}}$. The improved estimator of $\sigma_1^k$ is obtained as 
		$$\delta^3_{11}=\min \left\{\sqrt{\frac{\Gamma\left(\frac{p_1-k-1}{2}\right)}{2^k\Gamma\left(\frac{p_1+k-1}{2}\right)}},\alpha_{1}(1+U)^{\frac{k}{2}}\right\} S_1^{\frac{k}{2}}.$$
		
		\item [(iv)] With respect to linex loss function $L_4(t)$, the quantity $\alpha_{1}$ is defined as the solution to equation 
		$$\int_{0}^{\infty}z_1^{\frac{p_1+p_2+k-2}{2}-1}e^{a\alpha_{1}z_1^{\frac{k}{2}}-\frac{z_1}{2}}dz_1=e^a 2^{\frac{p_1+p_2+k-2}{2}} \Gamma\left(\frac{p_1+p_2+k-2}{2}\right).$$
			The improved estimator of $\sigma_1^k$ is obtained as 
			$\delta^4_{11}=\min \left\{c_{01},\alpha_{1}(1+U)^{\frac{k}{2}}\right\} S_1^{\frac{k}{2}}.$
			In particular for $k=2$, we have $\alpha_{1}=\frac{1}{2a}\left(1-e^{-\frac{2a}{p_1+p_2}}\right)$.
	\end{enumerate}
\end{example}

In the next theorem we have obtained a class of improved estimators using IERD approach \cite{kubokawa1994unified}.  
The joint density of $V_1$ and $U$ is
\begin{equation}
	f_{\eta}(v_1,u)\propto e^{-\frac{v_1}{2}(1+u\eta^2)}v_1^{\frac{p_1+p_2-2}{2}-1}u^{\frac{p_2-1}{2}-1}\eta^{p_2-1},~~~\forall ~~v_1>0,~ u>0,~ 0<\eta\leq 1.
\end{equation}
Define \begin{equation*}
	F_{\eta}(y,v_1)=\int_{0}^{y}f_{\eta}(s,v_1)ds ~~\text{and}~~	F_1(y,v_1)=\int_{0}^{y}f_1(s,v_1)ds.
\end{equation*}
\begin{theorem}\label{thkub1}
Suppose that the function $\phi_1$ satisfies the following conditions.
	\begin{enumerate}
		\item[(i)] $\phi_1(u)$ is increasing function in $u$ and $\lim\limits_{u\rightarrow\infty}\phi_1(u)=c_{01}$
		\item[(ii)] $\int_{0}^{\infty}L'(\phi_1(u)v_1^{\frac{k}{2}})v_1^{\frac{k}{2}}F_{\eta}(v_1,m)dv_1\geq0$
	\end{enumerate}
Then the risk of $\delta_{\phi_1}$ in (\ref{bae2.1}) is smaller than the $\delta_{01}$ under the loss function $L(t)$. 
\end{theorem}
\textbf{Proof:} Proof of this theorem is similar to the Theorem 4.1 of \cite{kubokawa1994double}. \\

\noindent Now, we obtain class of improved estimators for $\sigma_1^k$ under four special loss functions by applying Theorem \ref{thkub1} in the subsequent corollaries.
		\begin{corollary}
			Let us assume that the function $\phi_1(u)$ satisfies the following conditions 
			\begin{enumerate}
				\item[(i)] $\phi_1(u)$ is increasing function in $u$ and $\lim\limits_{u\rightarrow\infty}\phi_1(u)=\frac{\Gamma\left(\frac{p_1+k-1}{2}\right)}{2^{\frac{k}{2}}\Gamma\left(\frac{p_1+2k-1}{2}\right)}$
				\item[(ii)] $\phi_1(h)\geq \phi^1_{*}(u)$,  where
			\end{enumerate}
			\begin{align*}
				\phi_{*}^1(u)
				= \frac{ \Gamma\left(\frac{p_1+p_2+k-2}{2}\right)\displaystyle \int_{0}^{u}\frac{q^{\frac{p_2-3}{2}}}{(1+q)^\frac{p_1+p_2+k-2}{2}}dq}{2^{\frac{k}{2}} \Gamma\left(\frac{p_1+p_2+2k-2}{2}\right)\displaystyle\int_{0}^{u}\frac{q^{\frac{p_2-3}{2}}}{(1+q)^\frac{p_1+p_2+2k-2}{2}}dq}.
			\end{align*} 
			Then the risk of the estimator  $\delta_{\phi_1}$ given in (\ref{bae2.1}) is nowhere greater than that of $\delta^1_{01}$ under the quadratic loss function $L_1(t)$.
		\end{corollary}
		\begin{corollary}
			Under the loss function $L_2(t)$, the risk of the estimator  $\delta_{\phi_1}$ given in (\ref{bae2.1}) is nowhere greater than that of $\delta^2_{01}$ provided the function $\phi_1(u)$ satisfies  
			\begin{enumerate}
				\item [(i)] $\phi_1(u)$ is increasing function in $u$ and $\lim\limits_{u\rightarrow\infty}\phi_1(u)=\frac{\Gamma\left(\frac{p_1-1}{2}\right)}{2^{\frac{k}{2}}\Gamma\left(\frac{p_1+k-1}{2}\right)}$
				\item[(ii)] $\phi_1(u)\geq \phi^2_{*}(u)$, 	where
			\end{enumerate}
			\begin{align*}
				\phi_{*}^2(u)
				= \frac{ \Gamma\left(\frac{p_1+p_2-2}{2}\right)\displaystyle\int_{0}^{u}\frac{q^{\frac{p_2-3}{2}}}{(1+q)^\frac{p_1+p_2-2}{2}}dq}{2^{\frac{k}{2}} \Gamma\left(\frac{p_1+p_2+k-2}{2}\right)\displaystyle\int_{0}^{u}\frac{q^{\frac{p_2-3}{2}}}{(1+q)^\frac{p_1+p_2+k-2}{2}}dq}.
			\end{align*}
			
		\end{corollary}
		\begin{corollary}
			Suppose the following conditions hold true.  
			\begin{enumerate}
				\item[(i)] $\phi_1(u)$ is increasing function in $u$ and $\lim\limits_{u\rightarrow\infty}\phi_1(u)=\sqrt{\frac{\Gamma\left(\frac{p_1-k-1}{2}\right)}{2^k\Gamma\left(\frac{p_1+k-1}{2}\right)}}$
				\item[(ii)] $\phi_1(u)\geq \phi_{*}^3(u)$,  where
			\end{enumerate}
			\begin{align*}
				\phi_{*}^3(u)
				= \sqrt{\frac{ \Gamma\left(\frac{p_1+p_2-k-2}{2}\right)\displaystyle \int_{0}^{u}\frac{q^{\frac{p_2-3}{2}}}{(1+q)^\frac{p_1+p_2-k-2}{2}}dq}{2^k \Gamma\left(\frac{p_1+p_2+k-2}{2}\right)\displaystyle\int_{0}^{u}\frac{q^{\frac{p_2-3}{2}}}{(1+q)^\frac{p_1+p_2+k-2}{2}}dq}}.
			\end{align*}
			Then the risk of the estimator  $\delta_{\phi_1}$ given in (\ref{bae2.1}) is nowhere greater than that of $\delta^3_{01}$ under a symmetric loss function $L_3(t)$. 
\end{corollary}
\begin{corollary}
					Under the Linex loss function $L_4(t)$, the risk of the estimator  $\delta_{\phi_1}$ given in (\ref{bae2.1}) is nowhere greater than that of $\delta^4_{01}$ provided the function $\phi_1(u)$ satisfies  
				\begin{enumerate}
					\item [(i)] $\phi_1(u)$ is increasing function in $u$ and $\lim\limits_{u\rightarrow\infty}\phi_1(u)=c_{01}$
					\item[(ii)] $\phi_1(u)\geq \phi_{*}^4(u)$
				\end{enumerate}
				where the quantity $\phi_{*}^4(u)$ is defined as the solution of the  inequality
			\begin{equation*}\label{a1}
				\int_{0}^{\infty}\int_{0}^{u}v_1^{\frac{p_1+p_2+k-2}{2}-1}e^{a\phi_{1}(u)v_1^{\frac{k}{2}}-\frac{v_1}{2}(1+q)}q^{\frac{p_2-1}{2}}dqdv_1\geq e^a 	\int_{0}^{\infty}\int_{0}^{u}v_1^{\frac{p_1+p_2+k-2}{2}-1}e^{-\frac{v_1}{2}(1+q)}q^{\frac{p_2-1}{2}}dqdv_1
			\end{equation*} 
		\end{corollary}
		\begin{remark}
			In the above corollaries, we have obtained a class of improved estimators for $L_1,~L_2$ and $L_3$.  
			The boundary estimators of this class are obtained as $\delta_{\phi_{*}^1}=\phi_{*}^1S_1^{\frac{k}{2}}$, $\delta_{\phi_{*}^2}=\phi_{*}^2S_1^{\frac{k}{2}}$, $\delta_{\phi_{*}^3}=\phi_{*}^3S_1^{\frac{k}{2}}$ and $\delta_{\phi_{*}^4}=\phi_{*}^4S_1^{\frac{k}{2}}$. These estimators are \cite{brewster1974improving} type estimators. 
		\end{remark}
		\subsection{Generalized Bayes estimator of $\sigma_1^k$}
		In this subsection, we will derive generalized Bayes estimator of $\sigma_1^k$. We will prove that \cite{brewster1974improving} type estimator is a generalized Bayes. 
		Consider an improper prior 
		\begin{equation*}
			\pi(\underline{\theta})=\frac{1}{\sigma_1^4\sigma_2^4}, ~~ 0<\sigma_1\leq\sigma_2,~ \mu_1, \mu_2 \in \mathbb{R}.
		\end{equation*}
		For the quadratic loss function $L_1(t)$  the generalized Bayes estimator of $\sigma_1^k$ is obtained as
		\begin{equation*}
			\delta^1_{B1}=\frac{\int_{0}^{\infty} \int_{\sigma_1^2}^{\infty} \int_{0}^{\infty} \int_{0}^{\infty} 
					\frac{1}{\sigma_1^k} \, \pi(\underline{\theta} \mid x_1, x_2, s_1, s_2) \, d\mu_1 \, d\mu_2 \, d\sigma_2^2 \, d\sigma_1^2}{\int_{0}^{\infty} \int_{\sigma_1^2}^{\infty} \int_{0}^{\infty} \int_{0}^{\infty} 
					\frac{1}{\sigma_1^{2k}} \, \pi(\underline{\theta} \mid x_1, x_2, s_1, s_2) \,
					d\mu_1 \, d\mu_2 \, d\sigma_2^2 \, d\sigma_1^2}.
		\end{equation*}
		After simplification, we obtain the generalized Bayes estimator of $\sigma_1^k$ is
		\begin{equation*}
			\delta^1_{B1}=\frac{\int_{0}^{\infty}\int_{\sigma_1^2}^{\infty} \frac{1}{\sigma_1^{k+4}\sigma_2^4}\frac{1}e^{-\frac{s_1}{2\sigma_1^2}-\frac{s_2}{2\sigma_2^2}}\left(\frac{s_1}{\sigma_1^2}\right)^{\frac{p_1-3}{2}}\left(\frac{s_2}{\sigma_2^2}\right)^{\frac{p_2-3}{2}}d\sigma_2^2d\sigma_1^2}{\int_{0}^{\infty}\int_{\sigma_1^2}^{\infty} \frac{1}{\sigma_1^{2k+4}\sigma_2^4}\frac{1}e^{-\frac{s_1}{2\sigma_1^2}-\frac{s_2}{2\sigma_2^2}}\left(\frac{s_1}{\sigma_1^2}\right)^{\frac{p_1-3}{2}}\left(\frac{s_2}{\sigma_2^2}\right)^{\frac{p_2-3}{2}}d\sigma_2^2d\sigma_1^2}.
		\end{equation*}
		Using the transformation $v_1=\frac{s_1}{\sigma_1^2}$, $t_1=\frac{s_2}{s_1}\frac{\sigma_1^2}{\sigma_2^2}$, we get 
		\begin{equation*}
			\delta^1_{B1}=s_1^{\frac{k}{2}}\frac{\int_{0}^{\infty}\int_{0}^{u}e^{-\frac{v_1}{2}(1+t_1)}v_1^{\frac{p_1+p_2+k-4}{2}}t_1^{\frac{p_2-3}{2}}dt_1dv_1}{\int_{0}^{\infty}\int_{0}^{u}e^{-\frac{v_1}{2}(1+t_1)}v_1^{\frac{p_1+p_2+2k-4}{2}}q^{\frac{p_2-3}{2}}dt_1dv_1}
		\end{equation*} 
		which is $\delta_{\phi_{*}^1}(u)$, with $u=\frac{s_2}{s_1}$. Similarly we get the generalize Bayes estimator for $L_2(t)$ loss
		\begin{equation*}
			\delta^2_{B1}=s_1^{\frac{k}{2}}\frac{\int_{0}^{\infty}\int_{0}^{u}e^{-\frac{v_1}{2}(1+t_1)}v_1^{\frac{p_1+p_2-4}{2}}t_1^{\frac{p_2-3}{2}}dt_1dv_1}{\int_{0}^{\infty}\int_{0}^{u}e^{-\frac{v_1}{2}(1+t_1)}v_1^{\frac{p_1+p_2+k-4}{2}}q^{\frac{p_2-3}{2}}dt_1dv_1}
		\end{equation*} 
		which is $\delta_{\phi_{*}^2}(u)$. For the symmetric  $L_3(t)$, we obtain the generalized Bayes estimator as
		\begin{equation*}
			\delta^3_{B1}=s_1^{\frac{k}{2}}\sqrt{\frac{\int_{0}^{\infty}\int_{0}^{u}e^{-\frac{v_1}{2}(1+t_1)}v_1^{\frac{p_1+p_2-k-4}{2}}t_1^{\frac{p_2-3}{2}}dt_1dv_1}{\int_{0}^{\infty}\int_{0}^{u}e^{-\frac{v_1}{2}(1+t_1)}v_1^{\frac{p_1+p_2+k-4}{2}}q^{\frac{p_2-3}{2}}dt_1dv_1}}
		\end{equation*} 
		which is $\delta_{\phi_{*}^3}(u)$. 
\subsection{Improved estimation of $\sigma_1^k$ when $\mu_1\geq0$ and  $\mu_2\geq0$}\label{subsec2.2}
		In the above, we have obtained improved estimators of $\sigma_1^k$ when there is no restriction on the means. In this subsection, we consider the improved estimation of $\sigma_1^k$ when both means are non-negative, i.e., $\mu_1\geq0$ and $\mu_2\geq0$. In this context, we propose some more estimators that dominate BAEE. For this purpose, we consider a wider class of estimators similar to \cite{petropoulos2017estimation} as
	$$ \mathcal{C}_2=\left\{\delta_{\phi_2}=\phi_2(U,U_1)S_1^{\frac{k}{2}}:\ U_1=\frac{X_1}{\sqrt{S_1}},\ ~\phi_2(.) \mbox{ is a positive measurable function}\right\}.$$ 
		\begin{theorem}\label{the2.8}
		Let $Z_2 \sim \chi^2_{p_1+p_2+k-1}$  and $\alpha_2$ be a solution of the equation  $EL^{\prime} \left(Z_2^{k/2}\alpha_2\right)=0$. The risk of the estimator
		$$\delta_{\phi_{02}}= \begin{cases}\min \left\{\phi_2(U,U_1), c_{1,0}(U,U_1)\right\} S_1^{\frac{k}{2}}, & U_1>0 \\ \phi_2(U,U_1) S_1^{\frac{k}{2}}, & \text { otherwise }\end{cases}$$
		is nowhere larger than the estimator $\delta_{\phi_2}$ provided $P(c_{1,0}(U,U_1)<\phi_2(U,U_1)) > 0$ under a general scale invariant loss function $L(t)$, where $c_{1,0}(U,U_1)=\alpha_2\left(1+U+ p_1U_1^2\right)^{k/2}$. 	
	\end{theorem}
	\noindent \textbf{Proof:} Proof  is similar to Theorem \ref{th2.13}. We will prove Theorem \ref{th2.13}. 
	\begin{corollary}
		The estimator 
		$$\delta_{12}= \begin{cases}\min \left\{c_{01}, \alpha_2\left(1+U+p_1U_1^2\right)^{\frac{k}{2}}\right\} S_1^{\frac{k}{2}}, & U_1>0 \\ 
			c_{01} S_1^{\frac{k}{2}}, & \text { otherwise }\end{cases}$$
		dominates the BAEE under a general scale invariant loss function $L(t)$ provided $\alpha_2< c_{01}$.
	\end{corollary}
	
	Now using the information contained in both the sample, we consider a larger class of estimators as
	$$\mathcal{C}_3= \left\{\delta_{\phi_3}=\phi_3(U,U_1,U_2)S_1^{\frac{k}{2}}: U_1=\frac{X_1}{\sqrt{S_1}}, U_2=\frac{X_2}{\sqrt{S_1}},~\phi_3(.) \mbox{ is a positive measurable function}\right\}.$$ 
	In the following theorem, we give sufficient conditions under which we will get an improved estimator. 
	\begin{theorem}\label{th2.13}
		Let $Z_3 \sim \chi^2_{p_1+p_2+k}$ and $\alpha_3$ be a solution of the equation 
		\begin{equation}\label{st2.6}
			EL^{\prime} \left(Z_3^{k/2}\alpha_3\right)=0.
		\end{equation}
		Then the risk of the estimator
		$$\delta_{\phi_{03}}= \begin{cases}\min \left\{\phi_3(U,U_1,U_2), c_{1,0,0}(U,U_1,U_2)\right\} S_1^{\frac{k}{2}}, & U_1>0,U_2>0 \\ \phi_3(U,U_1,U_2) S_1^{\frac{k}{2}}, & \text { otherwise }\end{cases}$$
		is nowhere larger than the estimator $\delta_{\phi_3}$ under a general scale invariant loss function $L(t)$ provided $P(c_{1,0,0}(U,U_1,U_2)<\phi_3(U,U_1,U_2)) >0$, where $c_{1,0,0}(U,U_1,U_2)=\alpha_3\left(1+U+ p_1U_1^2+p_2U_2^2\right)^{k/2}$. 	
	\end{theorem}
	\noindent \textbf{Proof:} The risk function of the estimator $\delta_{\phi_3}$ is
	\begin{eqnarray*}
		R\left(\underline{\theta},\delta_{\phi_3}\right)
		=E\left[E\left\{L\left(V_1^{\frac{k}{2}}\phi_3(U,U_1,U_2)\right)\big\rvert U, U_1,U_2\right\} \right].
	\end{eqnarray*}
	The conditional risk can be written as  $R_1(\underline{\theta},c)=E\left\{L\left(V_1^{\frac{k}{2}}c\right)\big\rvert U=u,U_1=u_1,U_2=u_2\right\}$. We have conditional density of $V_1$ given $U=u$, $U_1=u_1$, $U_2=u_2$ is $$g_{\eta,\eta_1,\eta_2}(v_1) \propto \displaystyle{e^{-\frac{v_1}{2}(1+u\eta^2)-\frac{p_1}{2}(u_1\sqrt{v_1}-\eta_1)^2 -\frac{p_2}{2}(u_2\sqrt{v_1}\eta-\eta_2)^2}v_1^{\frac{p_1+p_2-2}{2}}}$$ $v_1>0$, $u>0$, $u_1\in \mathbb{R}$, $u_2\in \mathbb{R}$, where $\eta=\frac{\sigma_1}{\sigma_2}<1$, $\eta_1=\frac{\mu_1}{\sigma_1}\geq0$ and $\eta_2=\frac{\mu_2}{\sigma_2}\geq0$.
Applying Lemma $3.4.2.$ from \cite{lehmann2005testing} repeatedly, we get  for all $c>0$
			 \begin{align*}
				E_{\eta,\eta_1,\eta_2}\left[L'\left(V_1^{k/2}c\right)V_1^{k/2}\right] \geq E_{\eta,\eta_1,0}\left[L'\left(V_1^{k/2}c\right)V_1^{k/2}\right] &\geq E_{\eta,0,0}\left[L'\left(V_1^{k/2}c\right)V_1^{k/2}\right] \\
				&\geq E_{1,0,0}\left[L'\left(V_1^{k/2}c\right)V_1^{k/2}\right]
			\end{align*}
 Let $c_{\eta,\eta_1,\eta_2}(u,u_1,u_2)$ is the unique minimizer of $R_1(\underline{\theta},c)$. Now take $c=c_{1,0,0}(u,u_1,u_2)$ we have
\begin{align*}
E_{\eta,\eta_1,\eta_2}\left[L'\left(V_1^{\frac{k}{2}}c_{1,0,0}(u,u_1,u_2)\right)V_1^{k/2}\right]  &\geq	E_{1,0,0}\left[L'\left(V_1^{\frac{k}{2}}c_{1,0,0}(u,u_1,u_2)\right)V_1^{k/2}\right]\\
&=0\\
&=E_{\eta,\eta_1,\eta_2}\left[L'\left(V_1^{\frac{k}{2}}c_{\eta,\eta_1,\eta_2}(u,u_1,u_2)\right)V_1^{k/2}\right].
\end{align*}
Since $L^{\prime}(t)$ is increasing then from the above the inequality, we have $c_{\eta,\eta_1,\eta_2}(u,u_1,u_2)\leq c_{1,0,0}(u,u_1,u_2)$, where  $ c_{1,0,0}(u,u_1,u_2)$ is the unique solution of $$E_{1,0,0}\left[L'\left(V_1^{\frac{k}{2}}c_{1,0,0}(u,u_1,u_2)\right)V_1^{k/2}\right]=0.$$  
Using the transformation $z_3=v_1\left(1+u+p_1u_1^2+p_2u_2^2\right)$ we obtain 
\begin{equation}
			EL'\left(Z_3^{k/2}c_{1,0,0}(u,u_1,u_2)\left(1+u+p_1u_1^2+p_2u_2^2\right)^{-k/2}\right)=0,
		\end{equation} 
		where $Z_3\sim \chi^2_{p_1+p_2+k}$. Comparing with equation (\ref{st2.6}) we get  
		$$c_{1,0,0}(u,u_1,u_2)=\alpha_3\left(1+u+p_1u_1^2+p_2u_2^2\right)^{k/2}.$$ Define a function  $\phi_{03}(u,u_1,u_2)=\min\left\{\phi_3(u,u_1,u_2),c_{1,0,0}(u,u_1,u_2)\right\}$. Now we have $$c_{\eta,\eta_1,\eta_2}(u,u_1,u_2)\le  c_{1,0,0}(u,u_1,u_2)=\phi_{03} < \phi_3(u,u_1,u_2)$$  provided $P(c_{1,0,0}(U,U_1,U_2)$ $<\phi_3(U,U_1,U_2))>0$. Hence we get $R_1(\underline{\theta},\phi_3) > R_{1}(\underline{\theta}, \phi_{03})$. This complete the proof of the result.
		\begin{corollary}
			The estimator 
			$$\delta_{13}= \begin{cases}\min \left\{c_{01}, \alpha_3\left(1+U+p_1U_1^2+p_2U_2^2\right)^{\frac{k}{2}}\right\} S_1^{\frac{k}{2}}, & U_1>0,U_2>0 \\ 
				c_{01} S_1^{\frac{k}{2}}, & \text { otherwise }\end{cases}$$
			dominates the BAEE under a general scale invariant loss function $L(t)$ provided $\alpha_3<c_{01}$.
		\end{corollary}
		\begin{example}
			\begin{enumerate}
				\item[(i)] For the quadratic loss function $L_1(t)$ we have $\alpha_2=\frac{\Gamma\left(\frac{p_1+p_2+k-1}{2}\right)}{2^{\frac{k}{2}}\Gamma\left(\frac{p_1+p_2+2k-1}{2}\right)}$ and $\alpha_{3}=\frac{\Gamma\left(\frac{p_1+p_2+k}{2}\right)}{2^{\frac{k}{2}}\Gamma\left(\frac{p_1+p_2+2k}{2}\right)}$. The improved estimators of $\sigma_1^k$ can be obtained as follows 
				$$\delta^1_{12}= \begin{cases}\min \left\{\frac{\Gamma\left(\frac{p_1+k-1}{2}\right)}{2^{\frac{k}{2}}\Gamma\left(\frac{p_1+2k-1}{2}\right)},\alpha_{2}\left(1+U+p_1U_1^2\right)^{\frac{k}{2}}\right\} S_1^{\frac{k}{2}}, & U_1>0 \\ \frac{\Gamma\left(\frac{p_1+k-1}{2}\right)}{2^{\frac{k}{2}}\Gamma\left(\frac{p_1+2k-1}{2}\right)} S_1^{\frac{k}{2}}, & \text { otherwise }\end{cases}$$
				$$\delta^1_{13}= \begin{cases}\min \left\{\frac{\Gamma\left(\frac{p_1+k-1}{2}\right)}{2^{\frac{k}{2}}\Gamma\left(\frac{p_1+2k-1}{2}\right)},\alpha_{3}\left(1+U+p_1U_1^2+p_2U_2^2\right)^{\frac{k}{2}}\right\} S_1^{\frac{k}{2}}, & U_1>0,U_2>0 \\ \frac{\Gamma\left(\frac{p_1+k-1}{2}\right)}{2^{\frac{k}{2}}\Gamma\left(\frac{p_1+2k-1}{2}\right)} S_1^{\frac{k}{2}}, & \text { otherwise }\end{cases}$$
				\item [(ii)] Under the entropy loss function $L_2(t)$ we get $\alpha_{2}=\frac{\Gamma\left(\frac{p_1+p_2-1}{2}\right)}{2^{\frac{k}{2}}\Gamma\left(\frac{p_1+p_2+k-1}{2}\right)}$, $\alpha_{3}=\frac{\Gamma\left(\frac{p_1+p_2}{2}\right)}{2^{\frac{k}{2}}\Gamma\left(\frac{p_1+p_2+k}{2}\right)}$. So we get the improved estimators as
				$$\delta^2_{12}= \begin{cases}\min \left\{\frac{\Gamma\left(\frac{p_1-1}{2}\right)}{2^{\frac{k}{2}}\Gamma\left(\frac{p_1+k-1}{2}\right)},\alpha_{2}\left(1+U+p_1U_1^2\right)^{\frac{k}{2}}\right\} S_1^{\frac{k}{2}}, &U_1>0 \\ \frac{\Gamma\left(\frac{p_1-1}{2}\right)}{2^{\frac{k}{2}}\Gamma\left(\frac{p_1+k-1}{2}\right)} S_1^{\frac{k}{2}}, & \text { otherwise }\end{cases}$$
				
				$$\delta^2_{13}= \begin{cases}\min \left\{\frac{\Gamma\left(\frac{p_1-1}{2}\right)}{2^{\frac{k}{2}}\Gamma\left(\frac{p_1+k-1}{2}\right)},\alpha_{3}\left(1+U+p_1U_1^2+p_2U_2^2\right)^{\frac{k}{2}}\right\} S_1^{\frac{k}{2}}, &U_1>0, U_2>0 \\ \frac{\Gamma\left(\frac{p_1-1}{2}\right)}{2^{\frac{k}{2}}\Gamma\left(\frac{p_1+k-1}{2}\right)} S_1^{\frac{k}{2}}, & \text { otherwise }\end{cases}$$
				\item [(iii)] For the symmetric loss function $L_3(t)$ we obtain $\alpha_{2}=\sqrt{\frac{\Gamma\left(\frac{p_1+p_2-k-1}{2}\right)}{2^k\Gamma\left(\frac{p_1+p_2+k-1}{2}\right)}}$ and $\alpha_{3}=\sqrt{\frac{\Gamma\left(\frac{p_1+p_2-k}{2}\right)}{2^k\Gamma\left(\frac{p_1+p_2+k}{2}\right)}}$. The improved estimators of $\sigma_1^k$ are obtained as 
				$$\delta^3_{12}= \begin{cases}\min \left\{\sqrt{\frac{\Gamma\left(\frac{p_1-k-1}{2}\right)}{2^k\Gamma\left(\frac{p_1+k-1}{2}\right)}},\alpha_{2}\left(1+U+p_1U_1^2\right)^{\frac{k}{2}}\right\} S_1^{\frac{k}{2}}, &U_1>0 \\ \sqrt{\frac{\Gamma\left(\frac{p_1-k-1}{2}\right)}{2^k\Gamma\left(\frac{p_1+k-1}{2}\right)}} S_1^{\frac{k}{2}}, & \text { otherwise }\end{cases}$$
				
				$$\delta^3_{13}= \begin{cases}\min \left\{\sqrt{\frac{\Gamma\left(\frac{p_1-k-1}{2}\right)}{2^k\Gamma\left(\frac{p_1+k-1}{2}\right)}},\alpha_{3}\left(1+U+p_1U_1^2+p_2U_2^2\right)^{\frac{k}{2}}\right\} S_1^{\frac{k}{2}}, &U_1>0, U_2>0 \\ \sqrt{\frac{\Gamma\left(\frac{p_1-k-1}{2}\right)}{2^k\Gamma\left(\frac{p_1+k-1}{2}\right)}} S_1^{\frac{k}{2}}, & \text { otherwise }\end{cases}$$
				
				\item [(iv)] Under the Linex loss function $L_4(t)$, the quantities $\alpha_{2}$ and $\alpha_{3}$ are defined as the solutions to equations \
				$$\int_{0}^{\infty}z_2^{\frac{p_1+p_2+k-1}{2}-1}e^{a\alpha_{2}z_2^{\frac{k}{2}}-\frac{z_2}{2}}dz_2=e^a 2^{\frac{p_1+p_2+k-1}{2}} \Gamma\left(\frac{p_1+p_2+k-1}{2}\right)$$ 
					and 
				$$\int_{0}^{\infty}z_3^{\frac{p_1+p_2+k}{2}-1}e^{a\alpha_{3}z_3^{\frac{k}{2}}-\frac{z_3}{2}}dz_3=e^a 2^{\frac{p_1+p_2+k}{2}} \Gamma\left(\frac{p_1+p_2+k}{2}\right)$$
				respectively. 
	Then the improved estimators of $\sigma_1^k$ are obtained as follows 
					$$\delta^4_{12}= \begin{cases}\min \left\{c_{01},\alpha_{2}\left(1+U+p_1U_1^2\right)^{\frac{k}{2}}\right\} S_1^{\frac{k}{2}}, &U_1>0 \\ c_{01} S_1^{\frac{k}{2}}, & \text { otherwise }\end{cases}$$
					$$\delta^4_{13}= \begin{cases}\min \left\{c_{01},\alpha_{3}\left(1+U+p_1U_1^2+p_2U_2^2\right)^{\frac{k}{2}}\right\} S_1^{\frac{k}{2}}, &U_1>0, U_2>0 \\ c_{01} S_1^{\frac{k}{2}}, & \text { otherwise }\end{cases}$$
					In particular for $k=2$, we obtained $\alpha_{2}=\frac{1}{2a}\left(1-e^{-\frac{2a}{p_1+p_2+1}}\right)$ and $\alpha_{3}=\frac{1}{2a}\left(1-e^{-\frac{2a}{p_1+p_2+2}}\right)$.
			\end{enumerate} 	
		\end{example}
\subsection{Improved estimation of $\sigma_1^k$ when $\mu_1\le \mu_2$}\label{subsec2.3}
In this subsection, we address the problem of estimating the parameter $\sigma_1^k$ under the order restriction $\mu_1 \le \mu_2$ and $\sigma_1 \le \sigma_2$. By incorporating this restriction on the parameter, we aim to construct estimators that dominate the BAEE. Now we consider a subgroup of the affine group $\mathcal{G}$ as 
$$\mathcal{G}_1=\{g_{a,b}: a>0, b\in \mathbb{R}\}$$ and this group act as follows
$$(X_1,X_2,S_1,S_2)\rightarrow(aX_1+b,aX_2+b,a^2S_1,a^2S_2).$$
Under this group a  class of $\mathcal{G}_1$ equivariant estimators is obtained as 

$$\mathcal{C}_4=\left\{\delta_{\phi_4}=\phi_4(U,U_3)S_1^{k/2}:\ U_3=(X_2-X_1) S_1^{-1/2}\  \mbox{and}\ \phi_4(.)\ \mbox{is a positive measurable function}\right\}.$$
\begin{theorem} \label{thst2sigma1}
	Let $\alpha_4$ be a solution of the equation 
	\begin{equation} \label{st2eq2}
		EL^{\prime} \left(Z_4^{k/2}\alpha_4\right)=0.
	\end{equation}
	where $Z_4 \sim \chi^2_{p_1+p_2+k-1}$.
	Then, the risk function of the estimator
	$$\delta_{\phi_{04}}= \begin{cases}\min \left\{\phi_4(U,U_3), c_{1,0}(U,U_3)\right\} S_1^{\frac{k}{2}}, &U_3>0 \\ \phi_4(U,U_3) S_1^{\frac{k}{2}}, & \text { otherwise },\end{cases}$$
	is nowhere larger than the estimator $\delta_{\phi_4}$ under a general scale invariant loss function $L(t)$ provided $P(c_{1,0}(U,U_3)<\phi_4(U,U_3)) >0$, where $c_{1,0}(U,U_3)=\alpha_4\left(1+U+ U_3^2 \left(1/p_1+1/p_2\right)^{-1}\right)^{k/2}$. 	
\end{theorem}
\textbf{Proof:}  The risk function of the estimator $\delta_{\phi_4}\left(\underline{X},\underline{S}\right)$ can be written as 
\begin{eqnarray*}
	R\left(\underline{\theta},\delta_{\phi_4}\right)
	=E\left[E\left\{L\left(V_1^{\frac{k}{2}}\phi_4(U,U_3)\right)\big\rvert U,U_3\right\} \right].
\end{eqnarray*}
We denote the conditional risk as  $R_1(\underline{\theta},c)=E\left\{L\left(V_1^{\frac{k}{2}}c\right)\big\rvert U=u,U_3=u_3\right\}$.
We have conditional distribution of $V_1$ given $U=u$, $U_3=u_3$ is 
$$g_{\eta,\rho_1}(v_1) \propto \displaystyle{e^{-\frac{v_1}{2}(1+u\eta^2)-\frac{1}{2\left(\frac{1}{p_1}+\frac{1}{p_2\eta^2}\right)}(u_3\sqrt{v_1}-\rho_1)^2}v_1^{\frac{p_1+p_2-1}{2}-1}},~v_1>0, u_3\in \mathbb{R},~ u>0,$$ 
where $\eta=\frac{\sigma_1}{\sigma_2}<1$ and $\rho_1=\frac{\mu_2-\mu_1}{\sigma_1}\geq0$. 
Now, for all $u_3>0$ we have  $\frac{g_{\eta,\rho_1}(v_1)}{g_{\eta,0}(v_1)}$ and  $\frac{g_{\eta,0}(v_1)}{g_{1,0}(v_1)}$ is increasing in $v_1$. Hence applying the Lemma $3.4.2$ from \cite{lehmann2005testing}, it follows that for all $c>0$ $$E_{\eta,\rho_1}\left[L'\left(V_1^{\frac{k}{2}}c\right)V_1^{k/2}\right]\geq	E_{\eta,0}\left[L'\left(V_1^{\frac{k}{2}}c\right)V_1^{k/2} \right]\geq	E_{1,0}\left[L'\left(V_1^{\frac{k}{2}}c\right)V_1^{k/2}\right]=0.$$
Let $c_{\eta,\rho_1}(u,u_3)$ is the unique minimizer of $R_1(\underline{\theta},c)$. For $c=c_{1,0}(u,u_3)$ we get 
\begin{align*}
	E_{\eta,\rho_1}\left[L'\left(V_1^{\frac{k}{2}}c_{1,0}(u,u_3)\right)V_1^{k/2}\right]\geq	E_{\eta,0}\left[L'\left(V_1^{\frac{k}{2}}c_{1,0}(u,u_3)\right)V_1^{k/2} \right]&\geq	E_{1,0}\left[L'\left(V_1^{\frac{k}{2}}c_{1,0}(u,u_3)\right)V_1^{k/2}\right]\\
	&=0\\
	&=E_{\eta,\rho_1}\left[L'\left(V_1^{\frac{k}{2}}c_{\eta,\rho_1}(u,u_3)\right)V_1^{k/2}\right].
\end{align*}
Since $L^{\prime}(t)$ is increasing then from the above the inequality, we have $c_{\eta,\rho_1}(u,u_3)\leq c_{1,0}(u,u_3)$, where  $ c_{1,0}(u,u_3)$ is the unique solution of $E_{1,0}\left[L'\left(V_1^{\frac{k}{2}}c_{1,0}(u,u_3)\right)V_1^{k/2}\right]=0$. Using the transformation $z_4=v_1\left(1+u+u_3^2\left(\frac{1}{p_1}+\frac{1}{p_2}\right)^{-1}\right)$ we obtain 
\begin{equation}
	EL'\left(Z_4^{k/2}c_{1,0}(u,u_3)\left(1+u+u_3^2\left(1/p_1+1/p_2\right)^{-1}\right)^{-k/2}\right)=0,
\end{equation} 
where $Z_4\sim \chi^2_{p_1+p_2+k-1}$. Comparing with equation (\ref{st2eq2}) we get  
$$c_{1,0}(u,u_3)=\alpha_4\left(1+u+u_3^2\left(1/p_1+1/p_2\right)^{-1}\right)^{k/2}.$$ Consider a function 
$\phi_{04}(u,u_3)=\min\left\{\phi_4(u,u_3),c_{1,0}(u,u_3)\right\}$. Now we have $c_{\eta,\rho_1}(u,u_3) \le c_{1,0}(u,u_3)=\phi_{04} <\phi_4(u,u_3)$  provided $P(c_{1,0}(U,U_3)<\phi_4(U,U_3))>0$. Hence we get $R_1(\underline{\theta},\phi_4) > R(\underline{\theta}, \phi_{04})$. This completes the proof of the result.

\begin{corollary}
	The estimator 
	$$\delta_{14}= \begin{cases}\min \left\{c_{01}, \alpha_4\left(1+U+U_3^2\left(1/p_1+1/p_2\right)^{-1}\right)^{\frac{k}{2}}\right\} S_1^{\frac{k}{2}}, &U_3>0 \\ 
		c_{01} S_1^{\frac{k}{2}}, & \text { otherwise }\end{cases}$$
	dominates $\delta_{01}$ under a general scale invariant loss function $L(t)$ provided $\alpha_4<c_{01}$.
\end{corollary}

\begin{example}\rm
	\begin{enumerate}
		\item[(i)] For the quadratic loss function $L_1(t)$ we have $\alpha_{4}=\frac{\Gamma\left(\frac{p_1+p_2+k-1}{2}\right)}{2^{\frac{k}{2}}\Gamma\left(\frac{p_1+p_2+2k-1}{2}\right)}$. The improved estimator of $\sigma_1^k$ is obtained as 
		$$	\delta^1_{14}= \begin{cases}\min \left\{\frac{\Gamma\left(\frac{p_1+k-1}{2}\right)}{2^{\frac{k}{2}}\Gamma\left(\frac{p_1+2k-1}{2}\right)},\alpha_{4}\left(1+U+U_3^2\left(1/p_1+1/p_2\right)^{-1}\right)^{\frac{k}{2}}\right\} S_1^{\frac{k}{2}}, & U_3>0 \\ \frac{\Gamma\left(\frac{p_1+k-1}{2}\right)}{2^{\frac{k}{2}}\Gamma\left(\frac{p_1+2k-1}{2}\right)} S_1^{\frac{k}{2}}, & \text { otherwise }\end{cases}$$
		
		\item [(ii)] Under the entropy loss function $L_2(t)$ we get $\alpha_{4}=\frac{\Gamma\left(\frac{p_1+p_2-1}{2}\right)}{2^{\frac{k}{2}}\Gamma\left(\frac{p_1+p_2+k-1}{2}\right)}$. So we get the improved estimator as 
		$$\delta^2_{14}= \begin{cases}\min \left\{\frac{\Gamma\left(\frac{p_1-1}{2}\right)}{2^{\frac{k}{2}}\Gamma\left(\frac{p_1+k-1}{2}\right)},\alpha_{4}\left(1+U+U_3^2\left(1/p_1+1/p_2\right)^{-1}\right)^{\frac{k}{2}}\right\} S_1^{\frac{k}{2}}, &U_3>0 \\ \frac{\Gamma\left(\frac{p_1-1}{2}\right)}{2^{\frac{k}{2}}\Gamma\left(\frac{p_1+k-1}{2}\right)} S_1^{\frac{k}{2}}, & \text { otherwise }\end{cases}$$
		\item [(iii)] For the symmetric loss function $L_3(t)$ we obtain $\alpha_{4}=\sqrt{\frac{\Gamma\left(\frac{p_1+p_2-k-1}{2}\right)}{2^k\Gamma\left(\frac{p_1+p_2+k-1}{2}\right)}}$. The improved estimator of $\sigma_1^k$ is obtained as 
		$$\delta^3_{14}= \begin{cases}\min \left\{\sqrt{\frac{\Gamma\left(\frac{p_1-k-1}{2}\right)}{2^k\Gamma\left(\frac{p_1+k-1}{2}\right)}},\alpha_{4}\left(1+U+U_3^2\left(1/p_1+1/p_2\right)^{-1}\right)^{\frac{k}{2}}\right\} S_1^{\frac{k}{2}}, &U_3>0 \\ \sqrt{\frac{\Gamma\left(\frac{p_1-k-1}{2}\right)}{2^k\Gamma\left(\frac{p_1+k-1}{2}\right)}} S_1^{\frac{k}{2}}, & \text { otherwise }\end{cases}$$
		
		\item [(iv)]  Under the Linex loss function $L_4(t)$,  the quantity $\alpha_{4}$ is defined as the solution to equation	 
		$$\int_{0}^{\infty}z_4^{\frac{p_1+p_2+k-1}{2}-1}e^{a\alpha_{4}z_4^{\frac{k}{2}}-\frac{z_4}{2}}dz_4=e^a 2^{\frac{p_1+p_2+k-1}{2}} \Gamma\left(\frac{p_1+p_2+k-1}{2}\right).$$ 
		Then the improved estimator of $\sigma_1^k$ is obtained as 
		$$\delta^4_{14}= \begin{cases}\min \left\{c_{01},\alpha_{4}\left(1+U+U_3^2\left(1/p_1+1/p_2\right)^{-1}\right)^{\frac{k}{2}}\right\} S_1^{\frac{k}{2}}, &U_3>0 \\ c_{01} S_1^{\frac{k}{2}}, & \text { otherwise }\end{cases}$$
		In particular for $k=2$, we have $\alpha_{4}=\frac{1}{2a}\left(1-e^{-\frac{2a}{p_1+p_2+1}}\right)$.
	\end{enumerate}
\end{example}

\section{Improved estimation for $\sigma_2^k$ when $\sigma_1 \le \sigma_2$} \label{sec3}
	In this section, we address the problem of estimating $\sigma_2^k$ under the restriction $\sigma_1\le \sigma_2$. Using the information from the first sample, we can consider estimators of the form
\begin{equation}\label{bae3.1}
	\mathcal{D}_1=\left\{\delta_{\psi_1}=\psi_{1} \left(W\right)S_2^{\frac{k}{2}}:  W=S_1S_2^{-1}  \mbox{ and }\psi_{1}(.) \mbox{ is positive measurable function}\right\}
\end{equation} 
We propose a \cite{stein1964} type improved estimator in the following theorem. 

	\begin{theorem}	\label{th4.3}
		Suppose $k>0$. Let $\alpha_1$ be a solution of the equation 
		\begin{equation}\label{s2st1}
			EL^{\prime} \left(Z_1^{k/2}\alpha_1\right)=0
		\end{equation}
		where $Z_1 \sim \chi^2_{p_1+p_2+k-2}$.
		Consider $\psi_{01}(W)=\max\{\psi_1(W), d_1(W)\}$, then the risk function of the estimator
		$\delta_{\psi_{01}}=\psi_{01}(W)S_2^{\frac{k}{2}}$ is nowhere larger than the estimator $\delta_{\psi_1}$ provided $P(\psi_1(W) <d_{1}(W)) > 0$ holds true.
	\end{theorem}
\noindent \textbf{Proof:} Proof of this theorem is similar to the Theorem \ref{th2.1}.\\

	\noindent In the following corollary we propose an estimator which improves upon the BAEE.  
	\begin{corollary}\label{coro3.2}
		The risk function of the estimator
		$\delta_{21}=\max \left\{c_{02}, \alpha_1(1+W)^{k/2}\right\} S_2^{k/2}$
		is nowhere larger than the estimator $\delta_{02}$ provided $ \alpha_1 < c_{02}$.
	\end{corollary}
	
		\begin{example}
		\begin{enumerate}
			\item[(i)] Under the quadratic loss function $L_1(t)$, we obtain $\alpha_1=\frac{\Gamma\left(\frac{p_1+p_2+k-2}{2}\right)}{2^{\frac{k}{2}}\Gamma\left(\frac{p_1+p_2+2k-2}{2}\right)}$ and  the improved estimator is obtained as 
			$$\delta_{21}^1=\max \left\{\frac{\Gamma\left(\frac{p_2+k-1}{2}\right)}{2^{\frac{k}{2}}\Gamma\left(\frac{p_2+2k-1}{2}\right)}, \alpha_1(1+W)^{\frac{k}{2}}\right\} S_2^{\frac{k}{2}}.$$

			\item[(ii)] For the entropy loss function $L_2(t)$, we found that  $\alpha_1=\frac{\Gamma\left(\frac{p_1+p_2-2}{2}\right)}{2^{\frac{k}{2}}\Gamma\left(\frac{p_1+p_2+k-2}{2}\right)}$ then the improved estimator is as follows 
			$$\delta_{21}^2=\max \left\{\frac{\Gamma\left(\frac{p_2-1}{2}\right)}{2^{\frac{k}{2}}\Gamma\left(\frac{p_2+k-1}{2}\right)},\alpha_1(1+W)^{\frac{k}{2}}\right\} S_2^{\frac{k}{2}}.$$
			
			\item[(iii)] Under the symmetric loss function $L_3(t)$. We have  $\alpha_1=\sqrt{\frac{\Gamma\left(\frac{p_1+p_2-k-2}{2}\right)}{2^k\Gamma\left(\frac{p_1+p_2+k-2}{2}\right)}}$, then the improve estimator is as follows
			$$\delta_{21}^3= \max \left\{\sqrt{\frac{\Gamma\left(\frac{p_2-k-1}{2}\right)}{2^k\Gamma\left(\frac{p_2+k-1}{2}\right)}}, \alpha_1(1+W)^{\frac{k}{2}}\right\} S_2^{\frac{k}{2}}.$$
		
	\item[(iv)]Under the  Linex loss function $L_4(t)$, the improve estimator is
		$$\delta_{21}^4= \max \left\{c_{02}, \alpha_1(1+W)^{\frac{k}{2}}\right\} S_2^{\frac{k}{2}}$$ where $\alpha_1$ is the solution to equation 
				$$\int_{0}^{\infty}z_1^{\frac{p_1+p_2+k-2}{2}-1}e^{a\alpha_1z_1^{\frac{k}{2}}-\frac{z_1}{2}}dz_1=e^a 2^{\frac{p_1+p_2+k-2}{2}} \Gamma\left(\frac{p_1+p_2+k-2}{2}\right).$$  
		In particular for $k=2$, then we have, $\alpha_1=\frac{1}{2a}\left(1-e^{-\frac{2a}{p_1+p_2}}\right)$.
		\end{enumerate}
	\end{example}
	
	In the following theorem, we derive a class of improved estimators using the IERD approach \cite{kubokawa1994unified}.		
\begin{theorem}\label{th3.3}
		Let the function $\psi_1$ satisfies the following conditions.
			
		\begin{enumerate}
			\item[(i)] $\psi_1(w)$ is increasing function in $w$ and $\lim\limits_{w\rightarrow 0}\psi_1(w)=c_{02}$.
			\item[(ii)] $\int_{0}^{\infty}\int_{v_2w}^{\infty}L'(\psi_1(w)v_2^{k/2})v_2^{\frac{k}{2}}\nu_1(y)\nu_2(v_2)dydv_2\leq0$.
		\end{enumerate}
	
		where $\nu_i$ is pdf of $\chi^2_{p_i-1}$ for $i=1,2$. Then the risk of $\delta_{\psi_1}$ in (\ref{bae3.1}) is uniformly smaller than the estimator $\delta_{02}$ under $L(t)$.
	\end{theorem}
	\noindent \textbf{Proof:} Proof of this theorem is similar to the Theorem 4.3 of \cite{kubokawa1994double}\\

In the following, we have obtained  improved estimators for $\sigma_2^k$ under three special loss functions by applying Theorem \ref{th3.3}. 
					\begin{corollary}
						Let us assume that the function $\psi_1(w)$  satisfies the subsequent criterion:
\begin{enumerate}
	\item[(i)]$\psi_1(w)$ is increasing function in $w$ and $\lim\limits_{w\rightarrow0}\psi_1(w)=\frac{\Gamma\left(\frac{p_2+k-1}{2}\right)}{2^{\frac{k}{2}}\Gamma\left(\frac{p_2+2k-1}{2}\right)} $.
	\item[(ii)]$\psi_1(w)\leq\psi_{*}^1(w)$
\end{enumerate}
	where
\begin{align*}
\psi_{*}^1(w)= \frac{\Gamma\left(\frac{p_1+p_2+k-2}{2}\right)\int_{w}^{\infty}\frac{q^{\frac{p_1-3}{2}}}{(1+q)^\frac{p_1+p_2+k-2}{2}}dq}{2^{\frac{k}{2}}\Gamma\left(\frac{p_1+p_2+2k-2}{2}\right)\int_{w}^{\infty}\frac{q^{\frac{p_1-3}{2}}}{(1+q)^\frac{p_1+p_2+2k-2}{2}}dq}
							.
						\end{align*}
Then under the loss function $L_1(t)$, the risk of the estimator  $\delta_{\psi_1}$ is nowhere larger  than that of $\delta^1_{02}$.
					\end{corollary}
					
\begin{corollary}
Let us assume that the function $\psi_1(w)$ satisfies the following conditions 
 \begin{enumerate}
	\item[(i)] $\psi_1(w)$ is increasing function in $w$ and $\lim\limits_{w\rightarrow0}\psi_1(w)=\frac{\Gamma\left(\frac{p_2-1}{2}\right)}{2^{\frac{k}{2}}\Gamma\left(\frac{p_2+k-1}{2}\right)}$.
	\item[(ii)] $\psi_1(w)\leq\psi_{*}^2(w)$
\end{enumerate}
where
	$$\psi_{*}^2(w)=\frac{\Gamma\left(\frac{p_1+p_2-2}{2}\right)\int_{w}^{\infty}\frac{q^{\frac{p_1-3}{2}}}{\left(1+q\right)^{\frac{p_1+p_2-2}{2}}}dq}{\Gamma\left(\frac{p_1+p_2+k-2}{2}\right)\int_{w}^{\infty}\frac{q^{\frac{p_1-3}{2}}}{\left(1+q\right)^{\frac{p_1+p_2+k-2}{2}}}dq}$$	
						
	The risk of the estimator  $\delta_{\psi_1}$ is uniformly smaller than that of $\delta^2_{02}$ with respect to $L_2(t)$, .
\end{corollary}
\begin{corollary}
	Let us assume that the following conditions holds true
\begin{enumerate}
	\item[(i)] $\psi_1(w)$ is increasing function in $w$ and $\lim\limits_{w\rightarrow0}\psi_1(w)=\sqrt{\frac{\Gamma\left(\frac{p_2-k-1}{2}\right)}{2^k\Gamma\left(\frac{p_2+k-1}{2}\right)}}$.
     \item[(ii)] $\psi_1(w)\leq\psi_{*}^3(w)$
\end{enumerate}
where
\begin{align*}
	\psi_{*}^3(w)
		=\sqrt{\frac{\Gamma\left(\frac{p_1+p_2-k-2}{2}\right)\int_{w}^{\infty}\frac{q^{\frac{p_1-3}{2}}}{(1+q)^\frac{p_1+p_2-k-2}{2}}dq}{2^k \Gamma\left(\frac{p_1+p_2+k-2}{2}\right)\int_{w}^{\infty}\frac{q^{\frac{p_1-3}{2}}}{(1+q)^\frac{p_1+p_2+k-2}{2}}dq}}
\end{align*}	
 The risk of the estimator  $\delta_{\psi_1}$ is nowhere larger  than that of $\delta^3_{02}$ with respect to the loss function $L_3(t)$.
\end{corollary}
	\begin{corollary}
	For the loss function $L_4(t)$, the risk of the estimator  $\delta_{\psi_1}$ given in (\ref{bae3.1}) is nowhere greater than that of $\delta^4_{02}$ provided the function $\psi_1(w)$ satisfies  
		\begin{enumerate}
			\item [(i)] $\psi_1(w)$ is increasing function in $w$ and $\lim\limits_{w\rightarrow0}\psi_1(w)=c_{02}$
			\item[(ii)] $\psi_1(w)\leq \psi_{*}^4(w)$
		\end{enumerate}
		where the quantity $\psi_{*}^4(w)$ is defined as the solution to inequality 
		$$	\int_{0}^{\infty}\int_{v_2w}^{\infty}v_2^{\frac{p_2+k-3}{2}}y^{\frac{p_1-3}{2}}e^{a\psi_{1}(w)v_2^{\frac{k}{2}}-\frac{v_2}{2}-\frac{y}{2}}dydv_2\leq e^a 	\int_{0}^{\infty}\int_{v_2w}^{\infty}v_2^{\frac{p_2+k-3}{2}}y^{\frac{p_1-3}{2}}e^{-\frac{v_2}{2}-\frac{y}{2}}dydv_2.$$
\end{corollary}
					
\begin{remark}
	In the above corollaries, we obtained a class of improved estimators for $L_1,~L_2$ and $L_3$. The boundary estimators of this class are obtained as $\delta_{\psi_{*}^1}=\psi_{*}^1S_2^{\frac{k}{2}}$, $\delta_{\psi_{*}^2}=\psi_{*}^2S_2^{\frac{k}{2}}$, $\delta_{\psi_{*}^3}=\psi_{*}^3S_2^{\frac{k}{2}}$ and $\delta_{\psi_{*}^4}=\psi_{*}^4S_2^{\frac{k}{2}}$. These estimators are \cite{brewster1974improving} type estimators. 
\end{remark}
\subsection{Generalized Bayes estimator of $\sigma_2^k$}
	Here we find the generalized Bayes estimator for $\sigma_2^k$, and we have proved that the \cite{brewster1974improving} type estimator is a generalized Bayes estimator. Consider an improper prior
					\begin{equation*}
						\pi(\underline{\theta})=\frac{1}{\sigma_1^4\sigma_2^4}, ~~ 0<\sigma_1\leq\sigma_2,~ \mu_1, \mu_2 \in\mathbb{R}.
					\end{equation*} 
					For the quadratic loss function $L_1(t)$ the generalized Bayes estimator  of $\sigma_2^k$ is obtain as 
					\begin{equation*}
					\delta^1_{B2}=\frac{\int_{0}^{\infty} \int_{\sigma_1^2}^{\infty} \int_{0}^{\infty} \int_{0}^{\infty} \frac{1}{\sigma_2^k} \, \pi(\underline{\theta}\mid x_1, x_2, s_1, s_2) \, d\mu_1 \, d\mu_2 \, d\sigma_2^2 \, d\sigma_1^2}{\int_{0}^{\infty} \int_{\sigma_1^2}^{\infty} \int_{0}^{\infty} \int_{0}^{\infty} \frac{1}{\sigma_2^{2k}} \, \pi(\underline{\theta} \mid x_1, x_2, s_1, s_2) \,d\mu_1 \, d\mu_2 \, d\sigma_2^2 \, d\sigma_1^2}.
					\end{equation*}
					After performing some calculations by taking the transformation $v_2=\frac{s_2}{\sigma_2^2}$, $t_2=\frac{s_1}{s_2}\frac{\sigma_2^2}{\sigma_1^2}$, we obtain the  generalized Bayes estimator of $\sigma_2^k$ is
					\begin{equation*}
						\delta^1_{B2}=s_2^{\frac{k}{2}}\frac{\int_{0}^{\infty}\int_{w}^{\infty}e^{-\frac{v_2}{2}(1+t_2)}v_2^{\frac{p_1+p_2+k-4}{2}}t_2^{\frac{p_1-3}{2}}dt_2dv_2}{\int_{0}^{\infty}\int_{w}^{\infty}e^{-\frac{v_2}{2}(1+t_2)}v_2^{\frac{p_1+p_2+2k-4}{2}}t_2^{\frac{p_1-3}{2}}dt_2dv_2}
					\end{equation*} 
					which is $\delta_{\psi_{*}^1}(w)$, where  $w=\frac{s_1}{s_2}$. By using the similar argument as for $L_2(t)$ we get the generalized Bayes for $L_2$  is 
					\begin{equation*}
						\delta^2_{B2}=s_2^{\frac{k}{2}}\frac{\int_{0}^{\infty}\int_{w}^{\infty}e^{-\frac{v_2}{2}(1+t_2)}v_2^{\frac{p_1+p_2-4}{2}}t_2^{\frac{p_1-3}{2}}dt_2dv_2}{\int_{0}^{\infty}\int_{w}^{\infty}e^{-\frac{v_1}{2}(1+t_2)}v_2^{\frac{p_1+p_2+k-4}{2}}t_2^{\frac{p_1-3}{2}}dt_2dv_2}
					\end{equation*} 
					which is $\delta_{\psi_{*}^2}(w)$. For the symmetric loss $L_3(t)$ we obtain the generalized Bayes estimator as 
					\begin{equation*}
						\delta^3_{B2}=s_2^{\frac{k}{2}}\sqrt{\frac{\int_{0}^{\infty}\int_{w}^{\infty}e^{-\frac{v_2}{2}(1+t_2)}v_2^{\frac{p_1+p_2-k-4}{2}}t_2^{\frac{p_1-3}{2}}dt_2dv_2}{\int_{0}^{\infty}\int_{w}^{\infty}e^{-\frac{v_2}{2}(1+t_2)}v_2^{\frac{p_1+p_2+k-4}{2}}t_2^{\frac{p_1-3}{2}}dt_2dv_2}}
					\end{equation*} 
					which is $\delta_{\psi_{*}^3}(w)$.

\subsection{Improved estimation of $\sigma_2^k$ when $\mu_1\geq0 , \mu_2\geq0$}
In this previous subsection we found improved estimators of $\sigma_2^k$ without any restriction  on the means. Now, we consider estimation of $\sigma_2^k$ when $\mu_1\geq 0$ and $\mu_2\geq 0$. In this setting, we propose some estimators that perform better than BAEE. Similar to \cite{petropoulos2017estimation}, we consider a class of estimators of the form
	$$\mathcal{D}_2=\left\{\delta_{\psi_2}=\psi_2(W,W_1)S_2^{\frac{k}{2}}: W=\frac{S_1}{S_2}, W_1=\frac{X_1}{\sqrt{S_2}}, \psi_2(.) \mbox{ is a positive measurable function.}\right\}$$ 
	\begin{theorem}		
Let $Z_2 \sim \chi^2_{p_1+p_2+k-1}$  and $\alpha_2$ be a solution of the equation 
		\begin{equation*}
		EL^{\prime} \left(Z_2^{k/2}\alpha_2\right)=0.
		\end{equation*}
Then the risk function of the estimator
$$\delta_{\psi_{02}}= \begin{cases}\max \left\{\psi_2(W,W_1), d_{1,0}(W,W_1)\right\} S_2^{\frac{k}{2}}, & W_1<0 \\ \psi_2(W,W_1) S_2^{\frac{k}{2}}, & \text { otherwise }\end{cases}$$
is nowhere larger than the estimator $\delta_{\psi_2}$ under a general scale invariant loss function $L(t)$ provided $P(d_{1,0}(W,W_1)>\psi_2(W,W_1)) > 0$, where $d_{1,0}(W,W_1)=\alpha_2\left(1+W+ p_1W_1^{2}\right)^{k/2}$. 	
    \end{theorem}
\noindent \textbf{Proof:} Proof  is similar to Theorem \ref{th2.13}.
					
	\begin{corollary}
	The estimator 
		$$\delta_{22}= \begin{cases}\max \left\{c_{02}, \alpha_2\left(1+W+p_1W_1^{2}\right)^{\frac{k}{2}}\right\} S_2^{\frac{k}{2}}, & W_1<0 \\ 
		c_{02} S_2^{\frac{k}{2}}, & \text { otherwise }\end{cases}$$
	dominates the BAEE under a general scale invariant loss function $L(t)$ provided $\alpha_2< c_{02}$.
	\end{corollary}
Next we consider another class of estimators of the form
  $$\mathcal{D}_3=\left\{\delta_{\psi_3}=\psi_3(W,W_1,W_2)S_2^{\frac{k}{2}}: W_1=\frac{X_1}{\sqrt{S_2}}, W_2=\frac{X_2}{\sqrt{S_2}},~\psi_3(.) \mbox{ is a positive measurable function.}\right\}$$ 
					\begin{theorem}\label{th3.11}
						Let $Z_3 \sim \chi^2_{p_1+p_2+k}$ and $\alpha_3$ be a solution of the equation 
						\begin{equation}\label{Y3.4}
							EL^{\prime} \left(Z_3^{k/2}\alpha_3\right)=0.
						\end{equation}
						Then the estimator
						$$\delta_{\psi_{03}}= \begin{cases}\max \left\{\psi_3(W,W_1,W_2), d_{1,0,0}(W,W_1,W_2)\right\} S_2^{\frac{k}{2}}, & W_1<0,W_2<0 \\ \psi_3(W,W_1,W_2) S_2^{\frac{k}{2}}, & \text { otherwise }\end{cases}$$
						dominates  $\delta_{\psi_3}$ provided $P(d_{1,0,0}(W,W_1,W_2)>\psi_3(W,W_1,W_2)) >0$ under a general scale invariant loss function $L(t)$, where $d_{1,0,0}(W,W_1,W_2)=\alpha_3\left(1+W+ p_1W_1^{2}+p_2W_2^{2}\right)^{k/2}$. 	
					\end{theorem}
					
\noindent\noindent\textbf{Proof:} Proof is similar to Theorem \ref{th2.13}.

					\begin{corollary}\label{c3.12}
						The estimator 
						$$\delta_{23}= \begin{cases}\max \left\{c_{02}, \alpha_3\left(1+W+p_1W_1^{2}+p_2W_2^{2}\right)^{\frac{k}{2}}\right\} S_2^{\frac{k}{2}}, & W_1<0, W_2<0 \\ 
							c_{02} S_2^{\frac{k}{2}}, & \text { otherwise }\end{cases}$$
						dominates the BAEE under a general scale invariant loss function provided $\alpha_3< c_{02}$.
					\end{corollary}
					\begin{example}
						\begin{enumerate}
							\item[(i)] For the quadratic loss function $L_1(t)$ we have $\alpha_2=\frac{\Gamma\left(\frac{p_1+p_2+k-1}{2}\right)}{2^{\frac{k}{2}}\Gamma\left(\frac{p_1+p_2+2k-1}{2}\right)}$ and  $\alpha_3=\frac{\Gamma\left(\frac{p_1+p_2+k}{2}\right)}{2^{\frac{k}{2}}\Gamma\left(\frac{p_1+p_2+2k}{2}\right)}$. The improved estimators of $\sigma_2^k$ are obtained as follows 
							
							$$	\delta^1_{22}= \begin{cases}\max \left\{c_{02},\alpha_2\left(1+W+p_1W_1^{2}\right)^{\frac{k}{2}}\right\} S_2^{\frac{k}{2}}, & W_1<0 \\
								c_{02} S_2^{\frac{k}{2}}, & \text { otherwise }\end{cases}$$
							$$	\delta^1_{23}= \begin{cases}\max \left\{c_{02},\alpha_3\left(1+W+p_1W_1^{2}+p_2W_2^{2}\right)^{\frac{k}{2}}\right\} S_2^{\frac{k}{2}}, & W_1<0,W_2<0 \\
								c_{02} S_2^{\frac{k}{2}}, & \text { otherwise }\end{cases}$$
							\item [(ii)] Under the entropy loss function $L_2(t)$ we get $\alpha_2=\frac{\Gamma\left(\frac{p_1+p_2-1}{2}\right)}{2^{\frac{k}{2}}\Gamma\left(\frac{p_1+p_2+k-1}{2}\right)}$ and $\alpha_3=\frac{\Gamma\left(\frac{p_1+p_2}{2}\right)}{2^{\frac{k}{2}}\Gamma\left(\frac{p_1+p_2+k}{2}\right)}$. So we get the improved estimators as
							$$\delta^2_{22}= \begin{cases}\max \left\{c_{02},\alpha_2\left(1+W+p_1W_1^{2}\right)^{\frac{k}{2}}\right\} S_1^{\frac{k}{2}}, &W_1<0 \\ 
								c_{02} S_2^{\frac{k}{2}}, & \text { otherwise }\end{cases}$$
							
							$$\delta^2_{23}= \begin{cases}\max \left\{c_{02},\alpha_3\left(1+W+p_1W_1^{2}+p_2W_2^{2}\right)^{\frac{k}{2}}\right\} S_2^{\frac{k}{2}}, &W_1<0, W_2<0 \\ 
								c_{02} S_2^{\frac{k}{2}}, & \text { otherwise }\end{cases}$$
					
						\item [(iii)] For the symmetric loss function $L_3(t)$ we obtain $\alpha_2=\sqrt{\frac{\Gamma\left(\frac{p_1+p_2-k-1}{2}\right)}{2^k\Gamma\left(\frac{p_1+p_2+k-1}{2}\right)}}$ and $\alpha_3=\sqrt{\frac{\Gamma\left(\frac{p_1+p_2-k}{2}\right)}{2^k\Gamma\left(\frac{p_1+p_2+k}{2}\right)}}$. The improved estimators of $\sigma_2^k$ are obtained as 
						$$\delta^3_{23}= \begin{cases}\max \left\{c_{02},\alpha_2\left(1+W+p_1W_1^{2}\right)^{\frac{k}{2}}\right\} S_2^{\frac{k}{2}}, &W_1<0 \\ 
							c_{02} S_2^{\frac{k}{2}}, & \text { otherwise }\end{cases}$$
						
						$$\delta^3_{23}= \begin{cases}\max \left\{c_{02},\alpha_3\left(1+W+p_1W_1^{2}+p_2W_2^{2}\right)^{\frac{k}{2}}\right\} S_2^{\frac{k}{2}}, &W_1<0, W_2<0 \\
							c_{02} S_2^{\frac{k}{2}}, & \text { otherwise }\end{cases}$$
						\item [(iv)]
						Under the Linex loss function $L_4(t)$.	The improved estimators of $\sigma_2^k$ are obtained as follows 
     					$$	\delta^4_{22}= \begin{cases}\max \left\{c_{02},\alpha_2\left(1+W+p_1W_1^{2}\right)^{\frac{k}{2}}\right\} S_2^{\frac{k}{2}}, & W_1<0 \\
							c_{02} S_2^{\frac{k}{2}}, & \text { otherwise }\end{cases}$$
						$$	\delta^4_{23}= \begin{cases}\max \left\{c_{02},\alpha_3\left(1+W+p_1W_1^{2}+p_2W_2^{2}\right)^{\frac{k}{2}}\right\} S_2^{\frac{k}{2}}, & W_1<0,W_2<0 \\
							c_{02} S_2^{\frac{k}{2}}, & \text { otherwise }\end{cases}$$ where the quantities $\alpha_2$ and $\alpha_3$ are defined as the solutions to equations
							$$\displaystyle\int_{0}^{\infty}z_2^{\frac{p_1+p_2+k-1}{2}-1}e^{a\alpha_2z_2^{\frac{k}{2}}-\frac{z_2}{2}}dz_2=e^a 2^{\frac{p_1+p_2+k-1}{2}} \Gamma\left(\frac{p_1+p_2+k-1}{2}\right)$$ 
						and  
						$$\displaystyle \int_{0}^{\infty}z_3^{\frac{p_1+p_2+k}{2}-1}e^{a\alpha_3z_3^{\frac{k}{2}}-\frac{z_3}{2}}dz_3=e^a 2^{\frac{p_1+p_2+k}{2}} \Gamma\left(\frac{p_1+p_2+k}{2}\right)$$
						respectively.
						In particular for $k=2$, we have $\alpha_2=\frac{1}{2a}\left(1-e^{-\frac{2a}{p_1+p_2+1}}\right)$, $\alpha_3=\frac{1}{2a}\left(1-e^{-\frac{2a}{p_1+p_2+2}}\right)$.
					\end{enumerate}
					\end{example}

\subsection{Improved estimation of $\sigma_2^k$ when $\mu_1\le \mu_2$} 
In this subsection, we develop improved estimators for  $\sigma_2^k$ under the ordered restriction $\mu_1 \le \mu_2$. Similar to Subsection \ref{subsec2.3}, we consider a class of estimators as 
\begin{equation*}
	\mathcal{D}_4=\left\{\delta_{\psi_4}=\psi_{4} \left(W,W_3\right)S_2^{\frac{k}{2}}:\ W_3=\left(X_1-X_2\right)S_2^{-\frac{1}{2}},\ \psi_4(.) \mbox{ is a positive measurable function.}\right\}
\end{equation*}
In the following, we propose a sufficient condition to derive an improved estimator. 
\begin{theorem}
	Let $k$ be a positive real number and $\alpha_4$ be a solution of the equation 
	\begin{equation}\label{eq4.3}
		EL^{\prime} \left(Z_4^{k/2}\alpha_4\right)=0
	\end{equation}
	where $Z_4 \sim \chi^2_{p_1+p_2+k-1}$. Then the risk function of the estimator
	\begin{equation*}
		\delta_{\psi_{04}}= \begin{cases}\max \left\{\psi_{4}\left(W,W_3\right), d_{1,0}(W,W_3)\right\} S_2^{\frac{k}{2}}, & W_3>0 \\ \psi_{4}\left(W,W_3\right) S_2^{\frac{k}{2}}, & \text { otherwise }\end{cases}
	\end{equation*}
	is nowhere larger than the estimator $\delta_{\psi_4}$ provided $P(d_{1,0}(W,W_3)>\psi_4(W,W_3)) >0$ under a general scale invariant loss function $L(t)$, where $d_{1,0}(W,W_3)=\alpha_4(1+W+W_3\left(1/p_1+1/p_2\right)^{-1})^{k/2}$. 	
\end{theorem}

	\noindent \textbf{Proof:}  Proof of this Theorem is similar to Theorem \ref{thst2sigma1}.

		\begin{corollary}
			The estimator 
			$$\delta_{24}= \begin{cases}\max \left\{c_{02}, \alpha_4\left(1+W+W_3^{2}\left(1/p_1+1/p_2\right)^{-1}\right)^{\frac{k}{2}}\right\} S_2^{\frac{k}{2}}, & W_3>0 \\ 
				c_{02} S_2^{\frac{k}{2}}, & \text { otherwise }\end{cases}$$
			dominates $\delta_{02}$ under a general scale invariant loss function provided $\alpha_4<c_{02}$.
		\end{corollary}
		\begin{example}
			\begin{enumerate}
				\item[(i)] Under the quadratic loss function $L_1(t)$, we get $\alpha_4=\frac{\Gamma\left(\frac{p_1+p_2+k-1}{2}\right)}{2^{\frac{k}{2}}\Gamma\left(\frac{p_1+p_2+2k-1}{2}\right)}$ and the  improved estimator is obtained as 
				
				$$\delta_{25}^1= \begin{cases}\max \left\{\frac{\Gamma\left(\frac{p_2+k-1}{2}\right)}{2^{\frac{k}{2}}\Gamma\left(\frac{p_2+2k-1}{2}\right)}, \alpha_4\left(1+W+W_3^{2}\left(1/p_1+1/p_2\right)^{-1}\right)^{\frac{k}{2}}\right\} S_2^{\frac{k}{2}}, & W_3>0 \\ \frac{\Gamma\left(\frac{p_2+k-1}{2}\right)}{2^{\frac{k}{2}}\Gamma\left(\frac{p_2+2k-1}{2}\right)} S_2^{\frac{k}{2}}, & \text { otherwise }\end{cases}$$
				
				\item[(ii)] For the entropy loss function $L_2(t)$, we found that  $\alpha_4=\frac{\Gamma\left(\frac{p_1+p_2-1}{2}\right)}{2^{\frac{k}{2}}\Gamma\left(\frac{p_1+p_2+k-1}{2}\right)}$ and the improved estimator is 
				
				$$\delta_{25}^2= \begin{cases}\max \left\{\frac{\Gamma\left(\frac{p_2-1}{2}\right)}{2^{\frac{k}{2}}\Gamma\left(\frac{p_2+k-1}{2}\right)}, \alpha_4\left(1+W+W_3^{2}\left(1/p_1+1/p_2\right)^{-1}\right)^{\frac{k}{2}}\right\} S_2^{\frac{k}{2}}, &W_3>0 \\ \frac{\Gamma\left(\frac{p_2-1}{2}\right)}{2^{\frac{k}{2}}\Gamma\left(\frac{p_2+k-1}{2}\right)} S_2^{\frac{k}{2}}, & \text { otherwise }\end{cases}$$
				\item[(iii)] For symmetric loss function $L_3(t)$, we have  $\alpha_4=\sqrt{\frac{\Gamma\left(\frac{p_1+p_2-k-1}{2}\right)}{2^k\Gamma\left(\frac{p_1+p_2+k-1}{2}\right)}}$ and the improve estimator is obtained as
				
				$$\delta_{25}^3= \begin{cases}\max \left\{\sqrt{\frac{\Gamma\left(\frac{p_2-k-1}{2}\right)}{2^k\Gamma\left(\frac{p_2+k-1}{2}\right)}},\alpha_4\left(1+W+W_3^{2}\left(1/p_1+1/p_2\right)^{-1}\right)^{\frac{k}{2}}\right\} S_2^{\frac{k}{2}}, & W_3>0 \\ \sqrt{\frac{\Gamma\left(\frac{p_2-k-1}{2}\right)}{2^k\Gamma\left(\frac{p_2+k-1}{2}\right)}}S_2^{\frac{k}{2}}, & \text { otherwise }\end{cases}$$
				\item[(iv)]Under the Linex loss function $L_4(t)$, the improved estimator is 
					$$\delta_{25}^4= \begin{cases}\max \left\{c_{02}, \textcolor{blue}{\alpha_4}\left(1+W+W_3^{2}\left(1/p_1+1/p_2\right)^{-1}\right)^{\frac{k}{2}}\right\} S_2^{\frac{k}{2}}, & W_3>0 \\ c_{02} S_2^{\frac{k}{2}}, & \text { otherwise },\end{cases}$$ where the quantity $\alpha_4$ is defined as the solution to equation	 
					\begin{equation*}\label{bet2}
						\int_{0}^{\infty}z_4^{\frac{p_1+p_2+k-1}{2}-1}e^{a\alpha_4z_4^{\frac{k}{2}}-\frac{y_4}{2}}dz_4=e^a 2^{\frac{p_1+p_2+k-1}{2}} \Gamma\left(\frac{p_1+p_2+k-1}{2}\right).
					\end{equation*}	 
					In particular for $k=2$, we have $\alpha_4=\frac{1}{2a}\left(1-e^{-\frac{2a}{p_1+p_2+1}}\right)$.
			\end{enumerate}
		\end{example}

\section{A simulation study} \label{sec4}
\begin{figure}[htbp]
	\centering
	\begin{subfigure}{\textwidth}
		\centering
		\begin{subfigure}{0.32\textwidth}
			\centering
			\includegraphics[height=3.9cm, width=\textwidth]{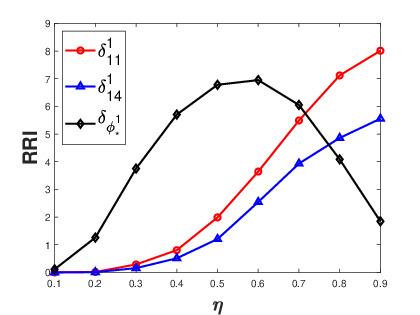}
			\caption{$(p_1,p_2)=(4,7)$, $(\mu_1,\mu_2)=(-0.5,-0.3)$}
		\end{subfigure}
		\hfill
		\begin{subfigure}{0.32\textwidth}
			\centering
			\includegraphics[height=3.9cm, width=\textwidth]{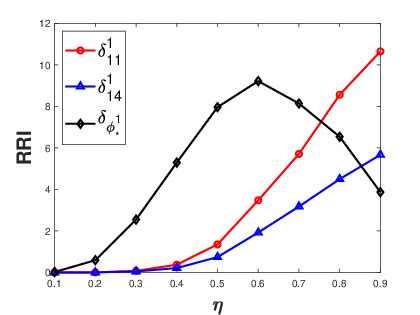}
			\caption{$(p_1,p_2)=(6,9)$, $(\mu_1,\mu_2)=(0,0)$}
		\end{subfigure}
		\hfill
		\begin{subfigure}{0.32\textwidth}
			\centering
			\includegraphics[height=3.9cm, width=\textwidth]{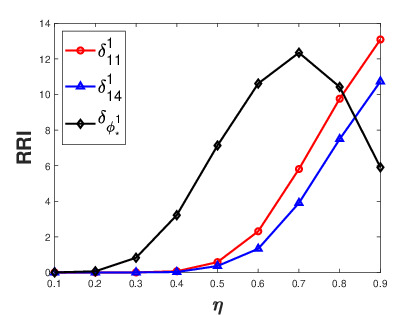}
			\caption{$(p_1,p_2)=(10,14)$, $(\mu_1,\mu_2)=(1.5,2)$}
		\end{subfigure}
	\end{subfigure}
	\hfill
	\begin{subfigure}{0.32\textwidth}
		\centering
		\includegraphics[height=3.9cm, width=\textwidth]{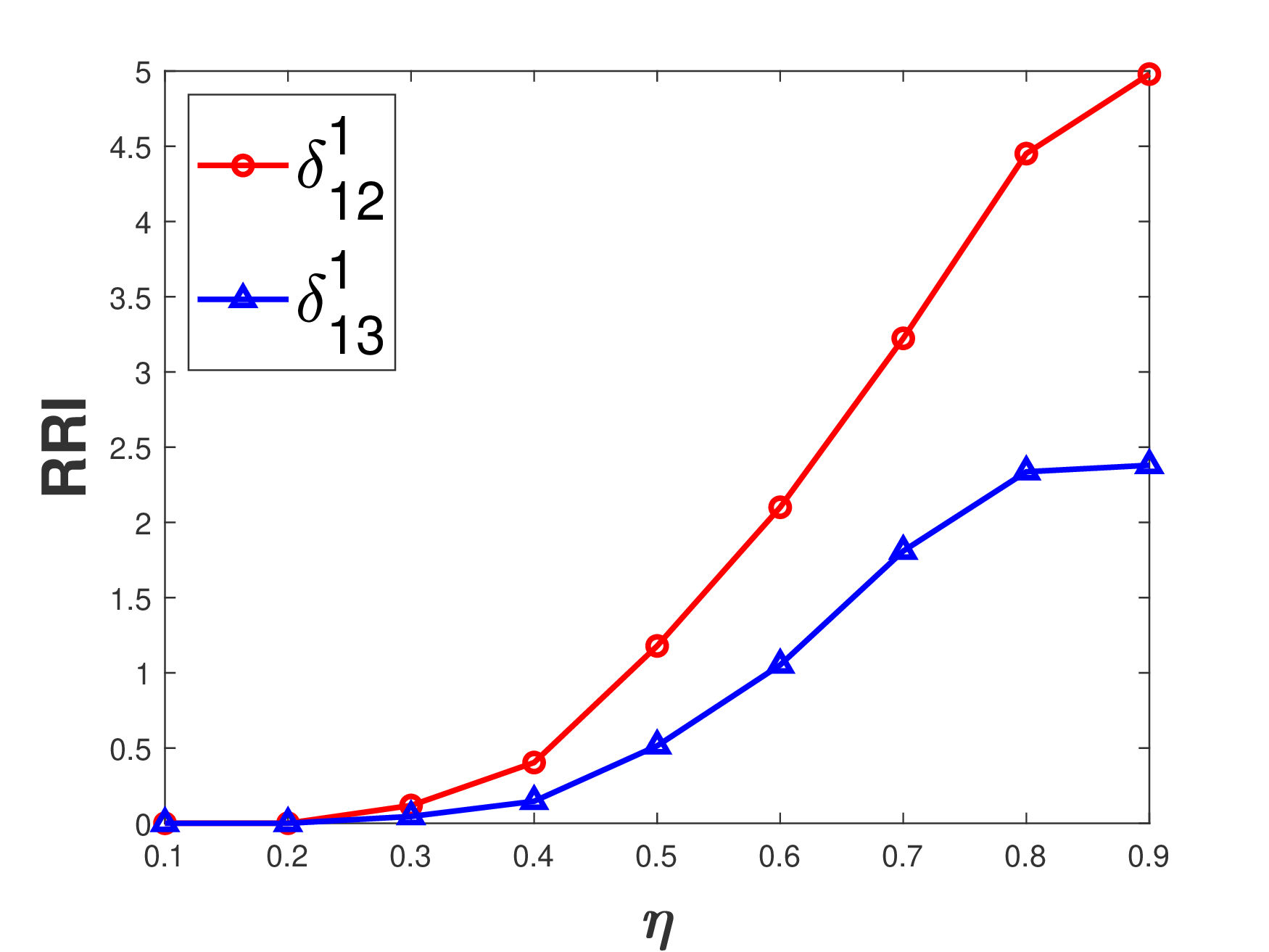}
		\caption{$(p_1,p_2)=(4,8)$, $(\mu_1,\mu_2)=(0,0)$}
	\end{subfigure}
	\hfill
	\begin{subfigure}{0.32\textwidth}
		\centering
		\includegraphics[height=3.9cm, width=\textwidth]{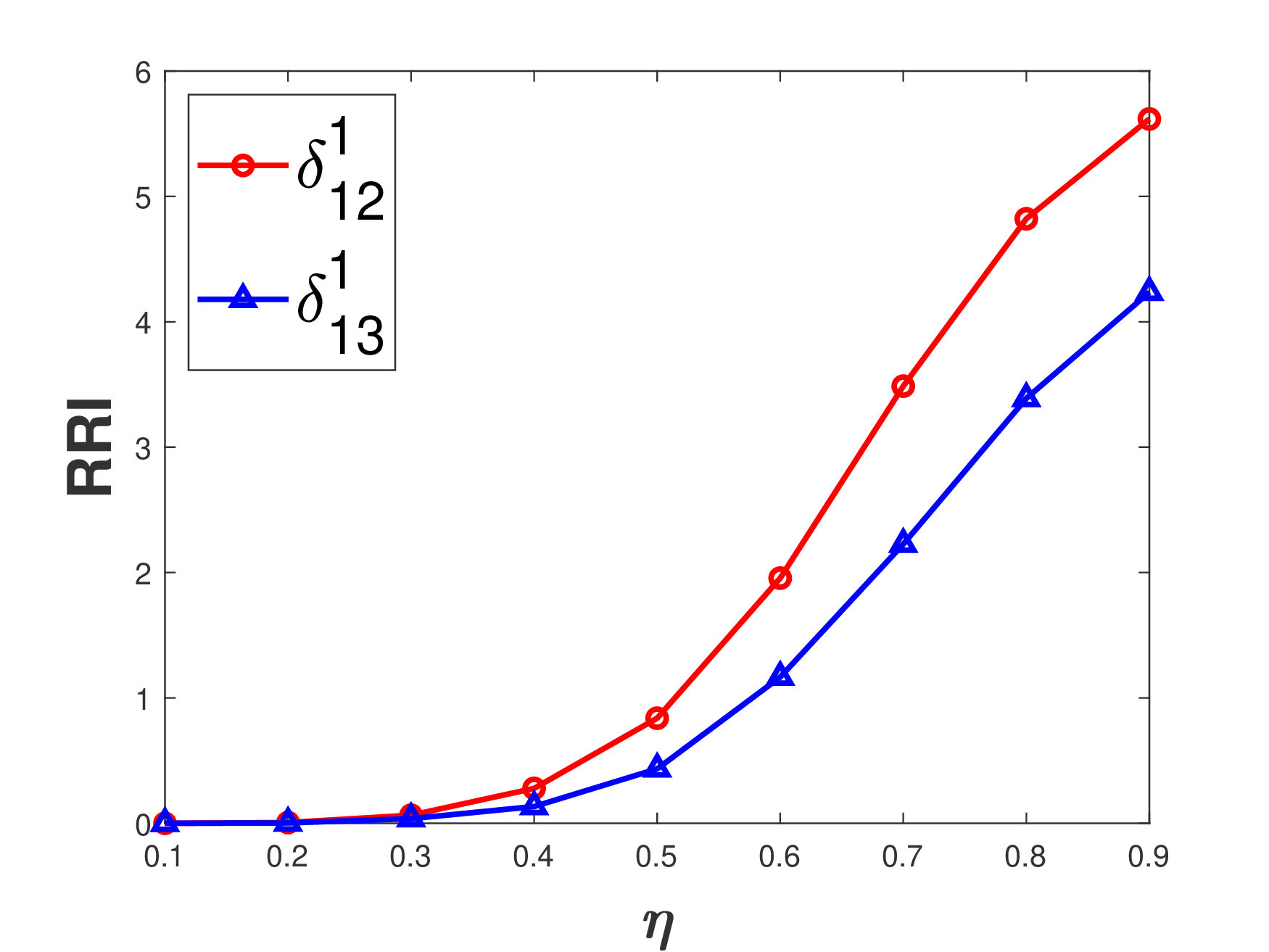}
		\caption{$(p_1,p_2)=(5,9)$, $(\mu_1,\mu_2)=(0,0.2)$}
	\end{subfigure}
	\hfill
	\begin{subfigure}{0.32\textwidth}
		\centering
		\includegraphics[height=3.9cm, width=\textwidth]{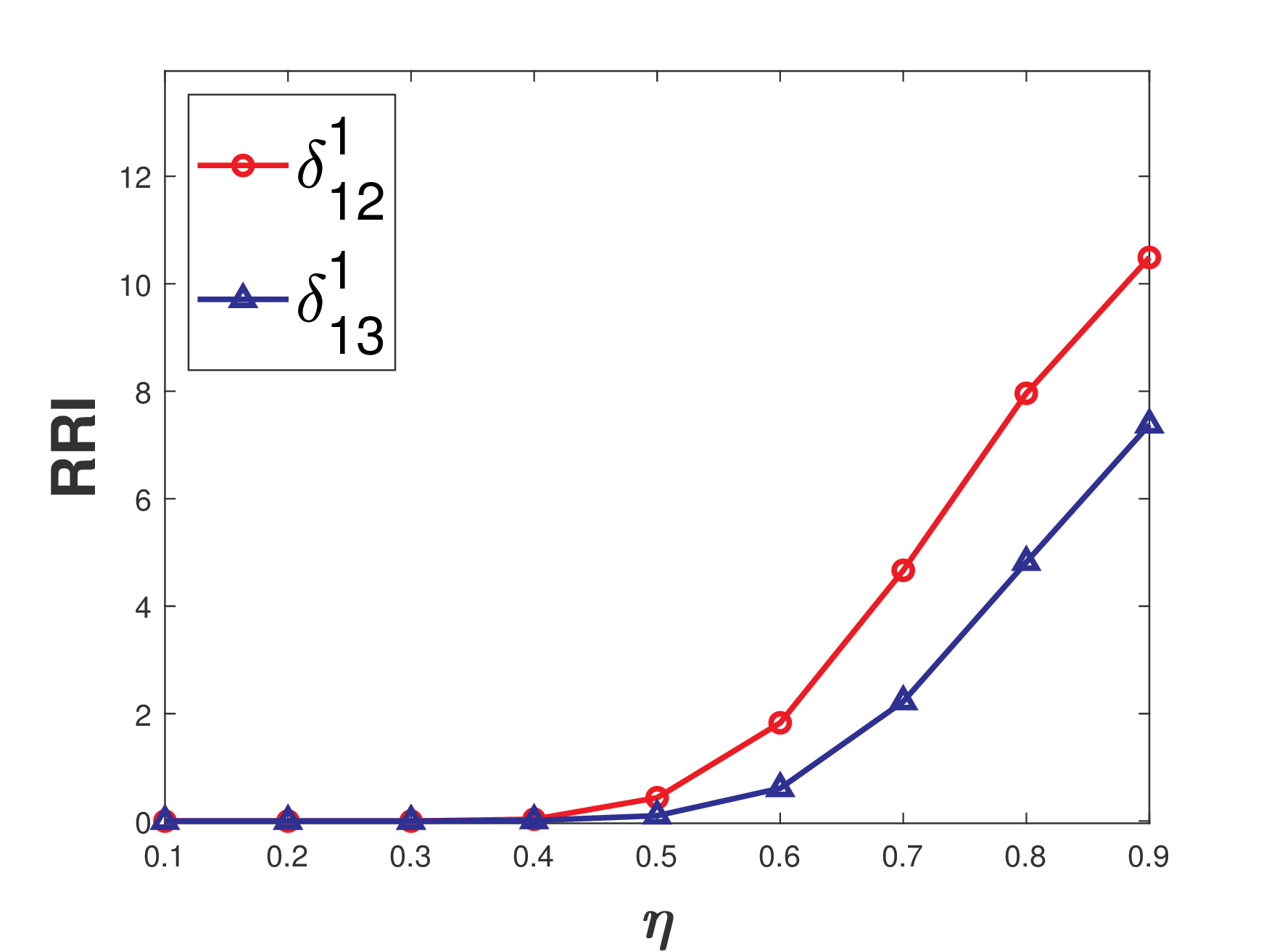}
		\caption{$(p_1,p_2)=(10,13)$, $(\mu_1,\mu_2)=(0.3,0.5)$}
	\end{subfigure}
	\hfill
\caption{ RRI of different estimators with respect to BAEE for $\sigma_1^2$ under $L_1(t)$.}\label{fig1}
\end{figure}

\begin{figure}[htbp]
	\begin{subfigure}{\textwidth}
		\centering
		\begin{subfigure}{0.32\textwidth}
			\centering
			\includegraphics[height=3.9cm, width=\textwidth]{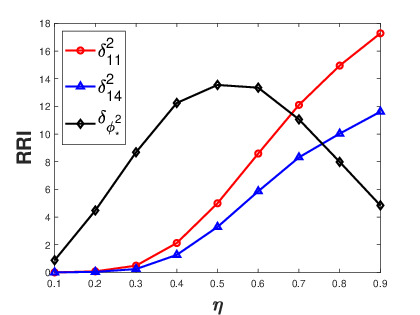}
			\caption{$(p_1,p_2)=(4,7)$, $(\mu_1,\mu_2)=(-0.5,-0.3)$}
		\end{subfigure}
		\hfill
		\begin{subfigure}{0.32\textwidth}
			\centering
			\includegraphics[height=3.9cm, width=\textwidth]{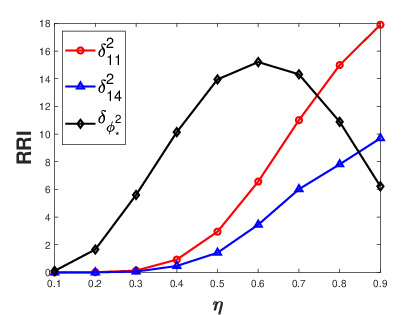}
			\caption{$(p_1,p_2)=(6,9)$, $(\mu_1,\mu_2)=(0,0)$}
		\end{subfigure}
		\hfill
		\begin{subfigure}{0.32\textwidth}
			\centering
			\includegraphics[height=3.9cm, width=\textwidth]{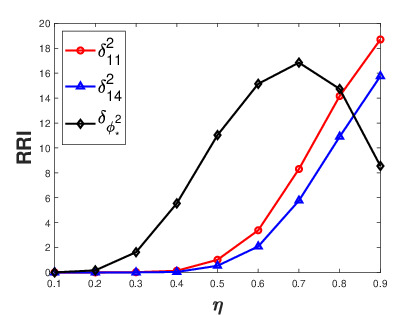}
			\caption{$(p_1,p_2)=(10,14)$, $(\mu_1,\mu_2)=(1.5,2)$}
		\end{subfigure}
        \hfill
       	\begin{subfigure}{0.32\textwidth}
       	\centering
       	\includegraphics[height=3.9cm, width=\textwidth]{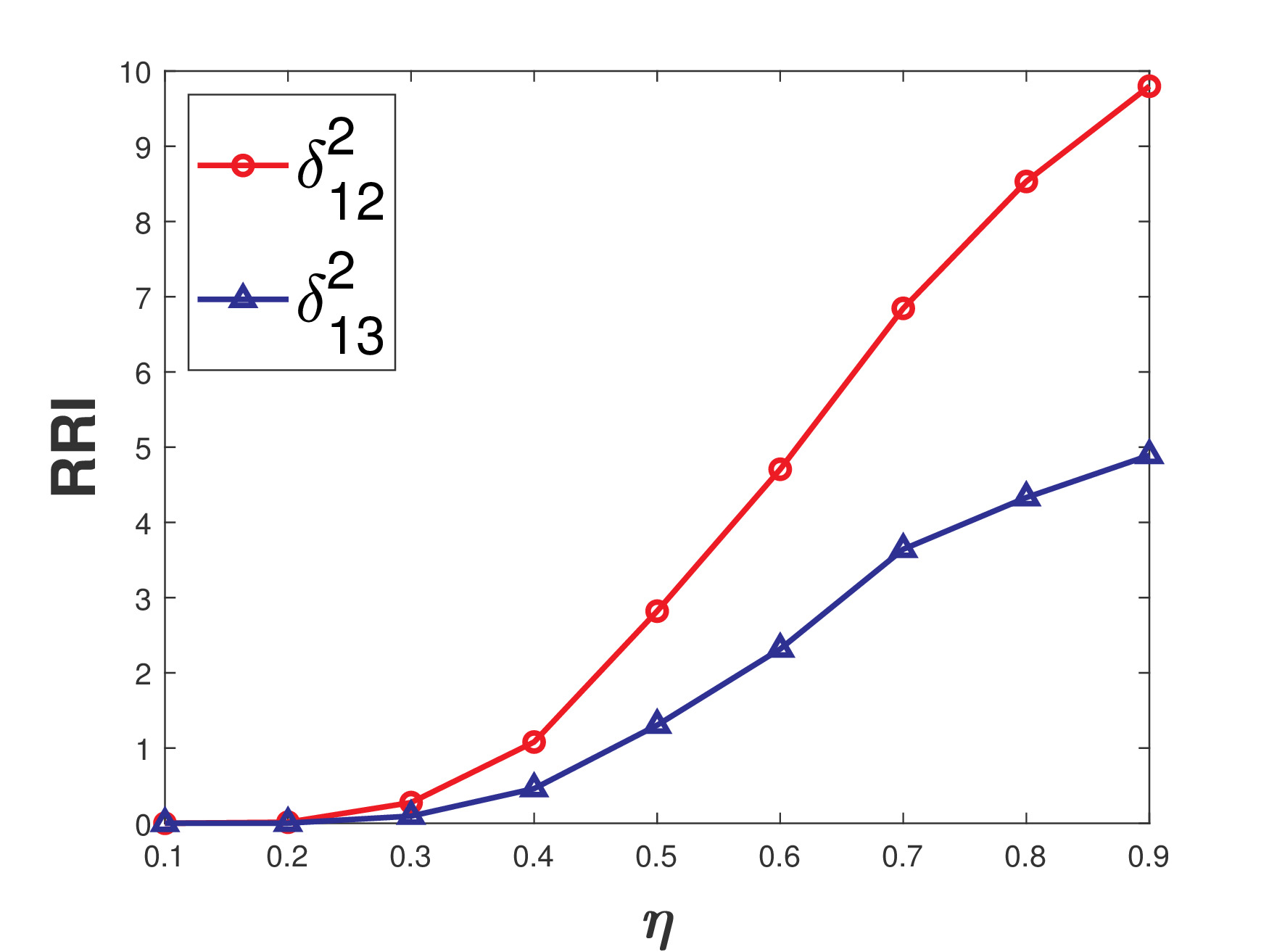}
       	\caption{$(p_1,p_2)=(4,8)$, $(\mu_1,\mu_2)=(0,0)$}
       \end{subfigure}
       \hfill
       \begin{subfigure}{0.32\textwidth}
       	\centering
       	\includegraphics[height=3.9cm, width=\textwidth]{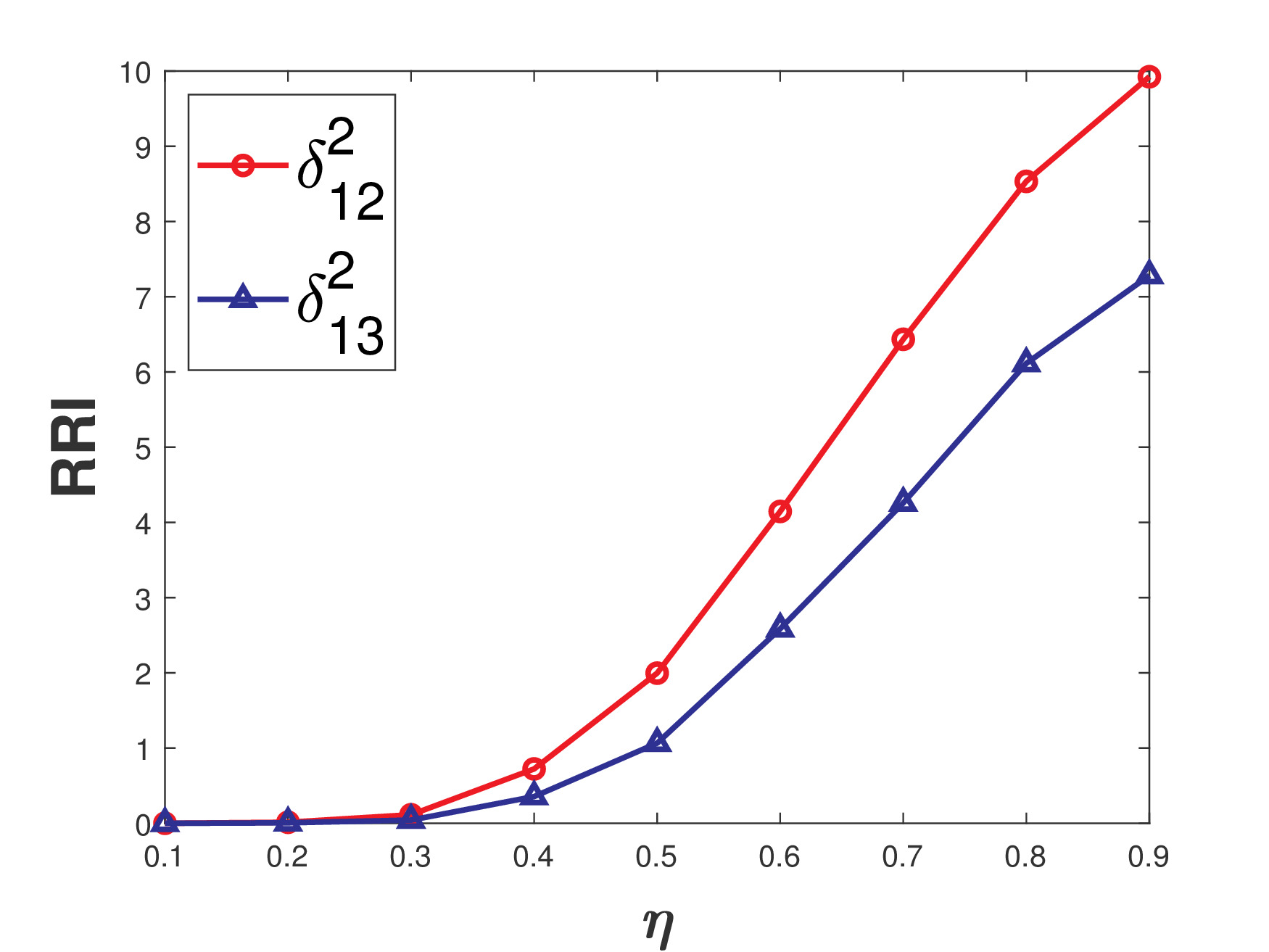}
       	\caption{$(p_1,p_2)=(5,9)$, $(\mu_1,\mu_2)=(0,0.2)$}
       \end{subfigure}
       \hfill
       \begin{subfigure}{0.32\textwidth}
       	\centering
       	\includegraphics[height=3.9cm, width=\textwidth]{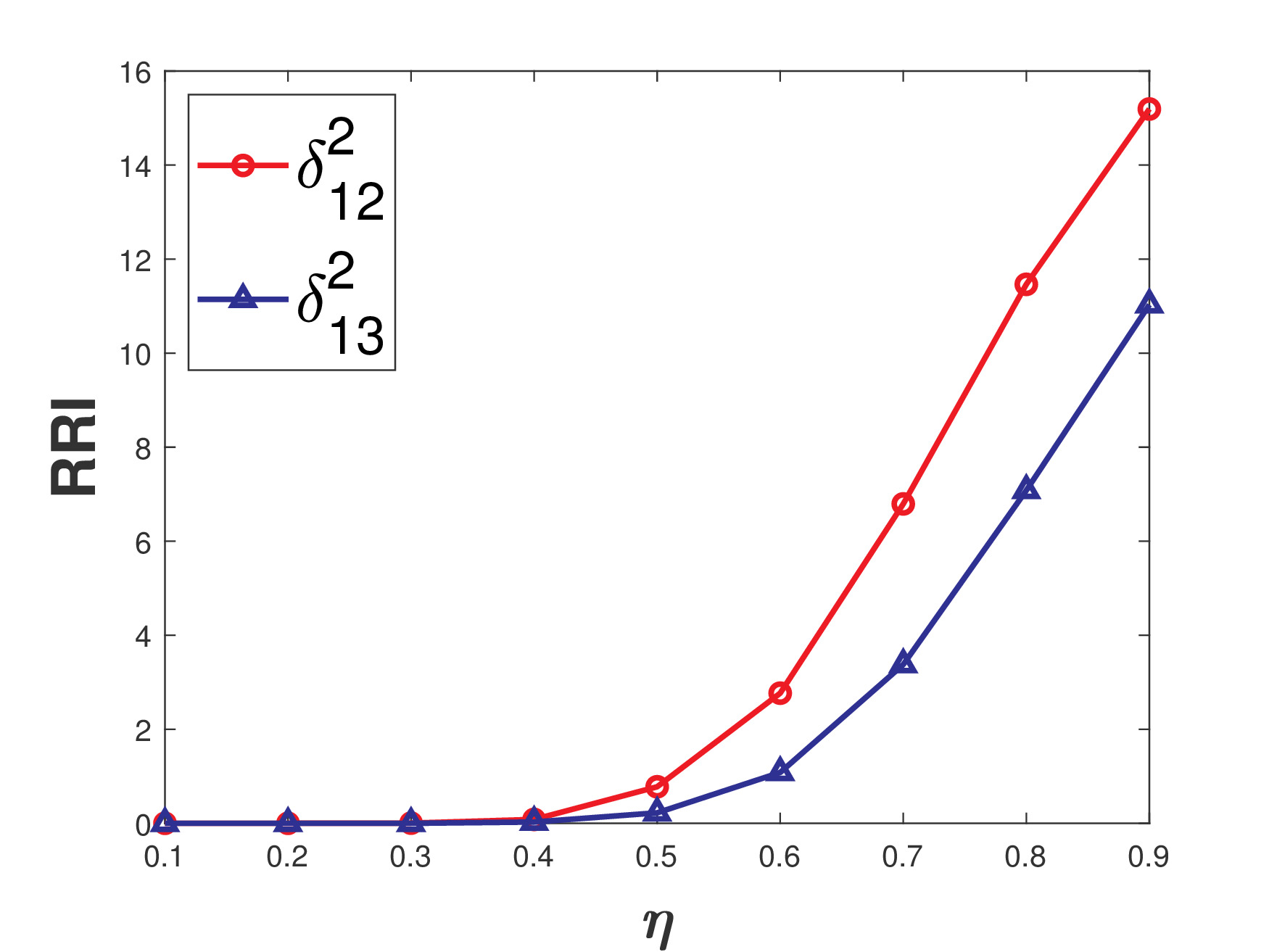}
       	\caption{$(p_1,p_2)=(10,13)$, $(\mu_1,\mu_2)=(0.3,0.5)$}
       \end{subfigure}
       
	\end{subfigure}
\caption{RRI of different estimators with respect to BAEE for $\sigma_1^2$ under $L_2(t)$.} \label{fig2}
\end{figure}

\begin{figure}[htbp]
	\begin{subfigure}{\textwidth}
		\centering
		\begin{subfigure}{0.32\textwidth}
			\centering
			\includegraphics[height=3.9cm, width=\textwidth]{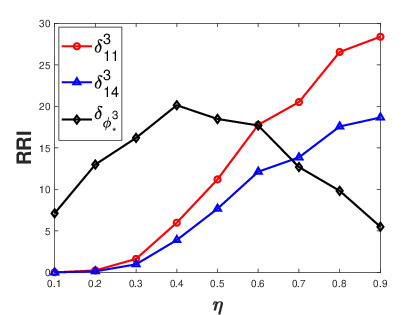}
			\caption{$(p_1,p_2)=(4,7)$, $(\mu_1,\mu_2)=(-0.5,-0.3)$}
		\end{subfigure}
		\hfill
		\begin{subfigure}{0.32\textwidth}
			\centering
			\includegraphics[height=3.9cm, width=\textwidth]{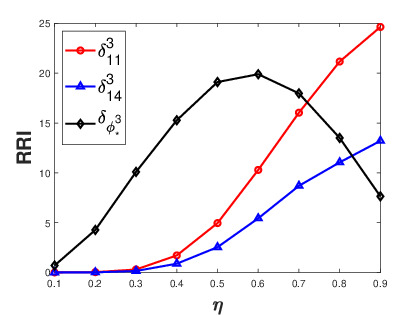}
			\caption{$(p_1,p_2)=(6,9)$, $(\mu_1,\mu_2)=(0,0)$}
		\end{subfigure}
		\hfill
		\begin{subfigure}{0.32\textwidth}
			\centering
			\includegraphics[height=3.9cm, width=\textwidth]{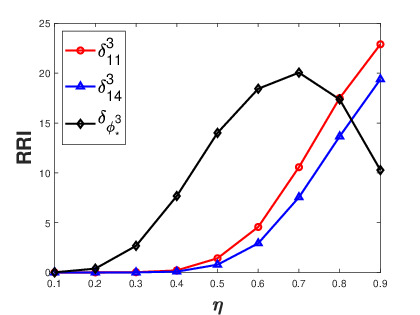}
			\caption{$(p_1,p_2)=(10,14)$, $(\mu_1,\mu_2)=(1.5,2)$}
		\end{subfigure}
       \hfill
       	\begin{subfigure}{0.32\textwidth}
       	\centering
       	\includegraphics[height=3.9cm, width=\textwidth]{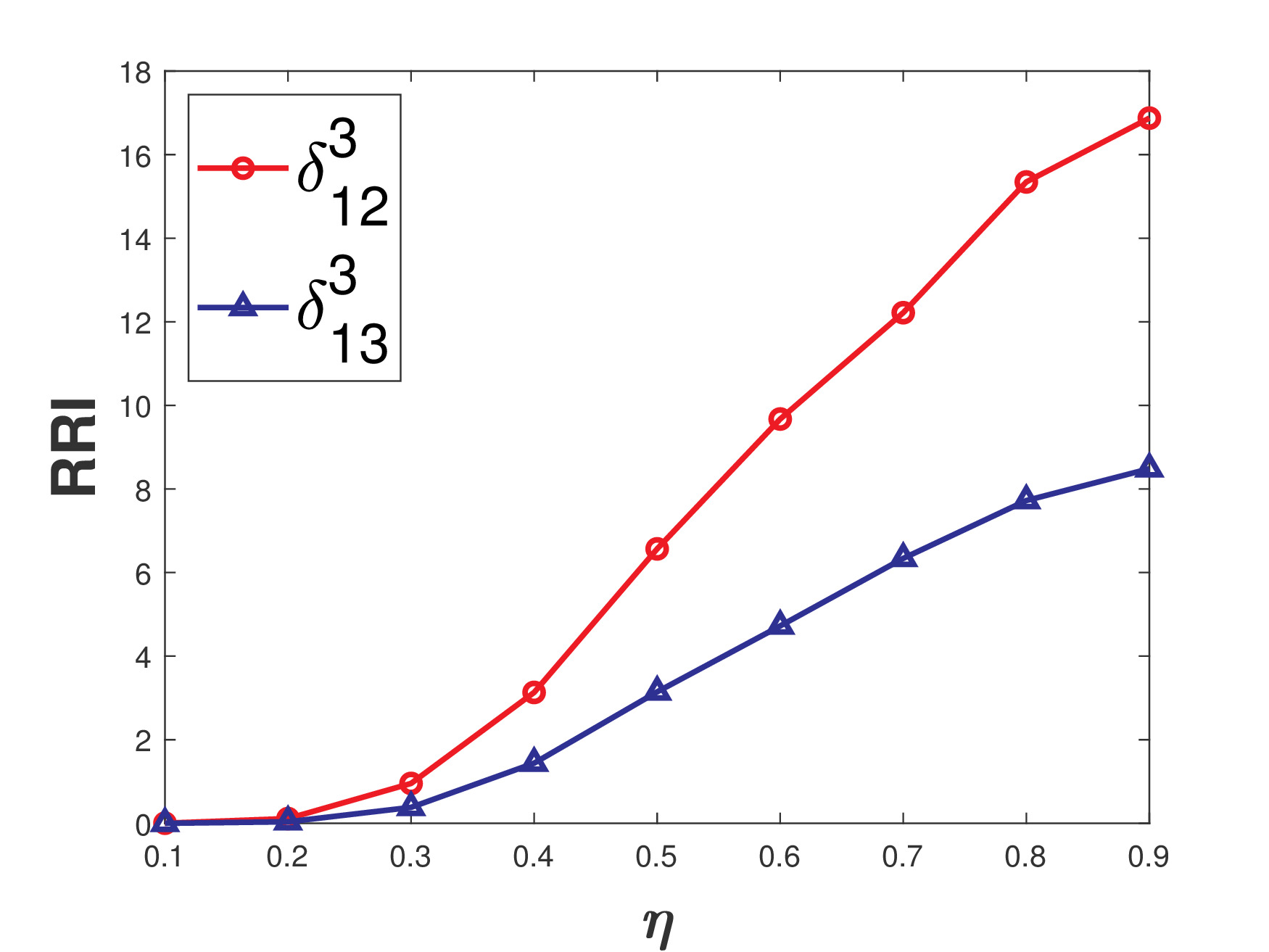}
       	\caption{$(p_1,p_2)=(4,8)$, $(\mu_1,\mu_2)=(0,0)$}
       \end{subfigure}
       \hfill
       \begin{subfigure}{0.32\textwidth}
       	\centering
       	\includegraphics[height=3.9cm, width=\textwidth]{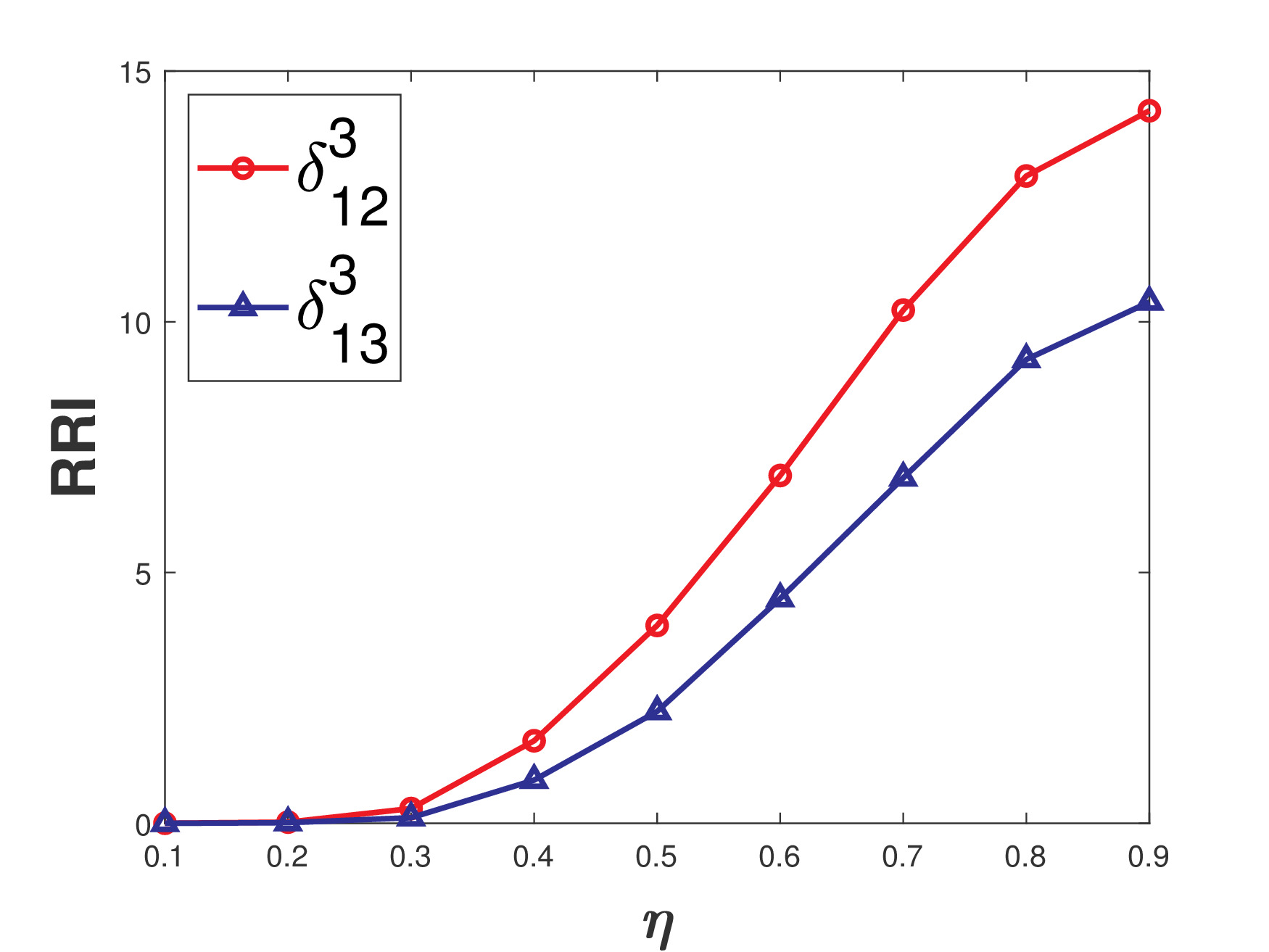}
       	\caption{$(p_1,p_2)=(5,9)$, $(\mu_1,\mu_2)=(0,0.2)$}
       \end{subfigure}
       \hfill
       \begin{subfigure}{0.32\textwidth}
       	\centering
       	\includegraphics[height=3.9cm, width=\textwidth]{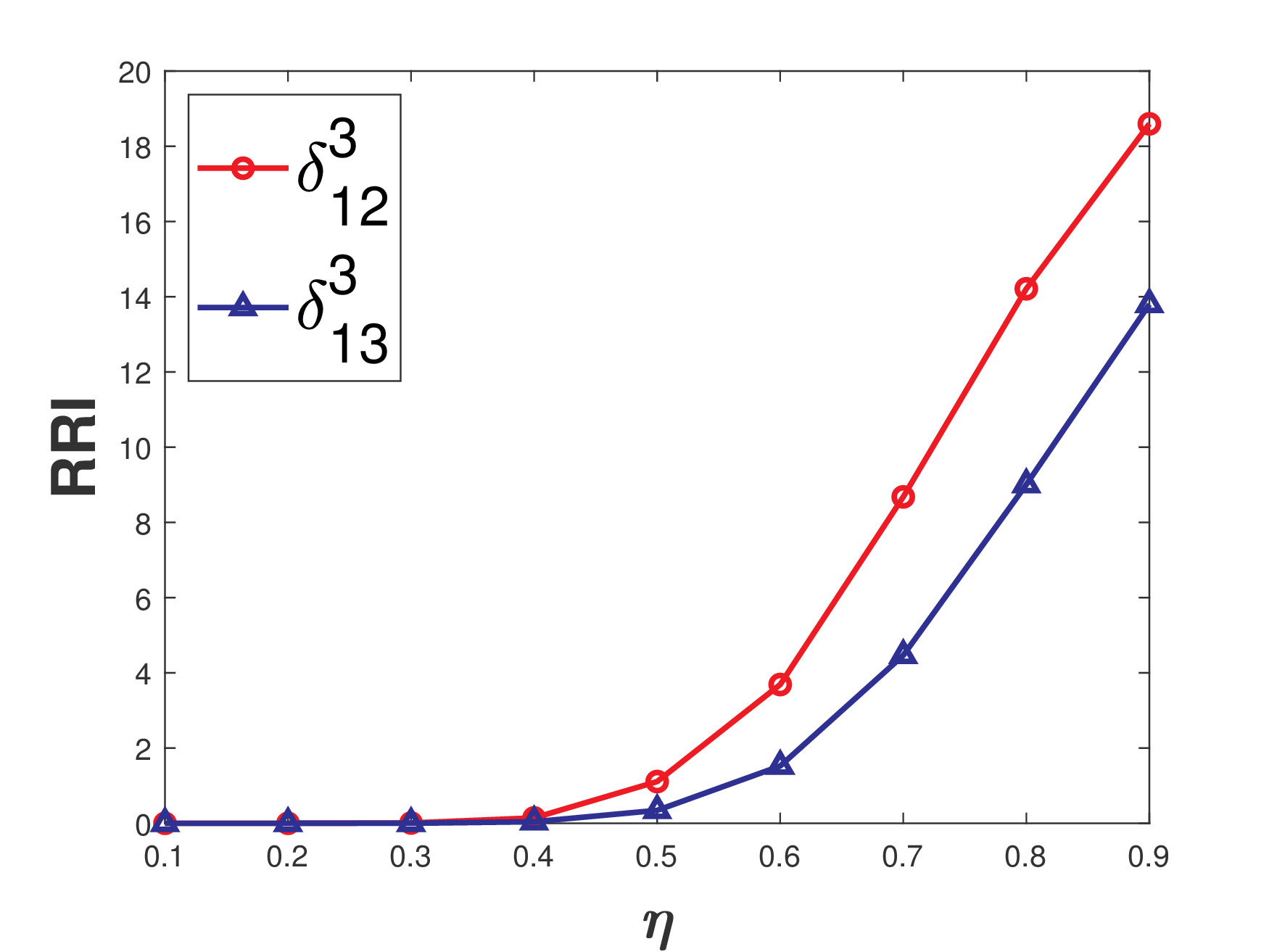}
       	\caption{$(p_1,p_2)=(10,13)$, $(\mu_1,\mu_2)=(0.3,0.5)$}
       \end{subfigure}
       
	\end{subfigure}
	\caption{ RRI of different estimators with respect to BAEE for $\sigma_1^2$ under $L_3(t)$.}\label{fig3}
	
\end{figure}

\begin{figure}[htbp]
\centering
	\begin{subfigure}{\textwidth}
		\centering
		\begin{subfigure}{0.32\textwidth}
			\centering
			\includegraphics[height=3.9cm, width=\textwidth]{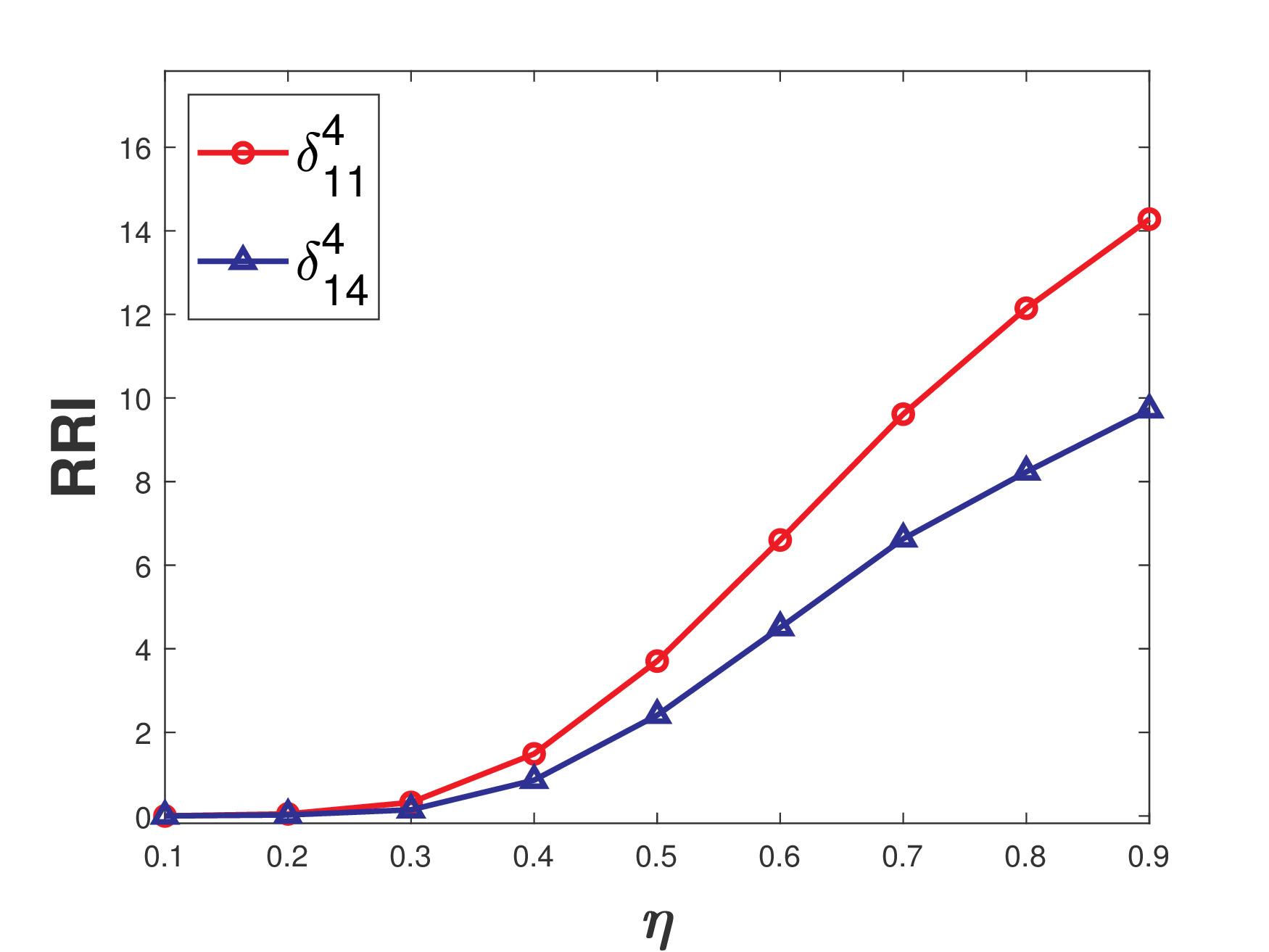}
			\caption{$a=-2$, $(p_1,p_2)=(4,7)$, $(\mu_1,\mu_2)=(-0.5,-0.3)$}
		\end{subfigure}
		\hfill
		\begin{subfigure}{0.32\textwidth}
			\centering
			\includegraphics[height=3.9cm, width=\textwidth]{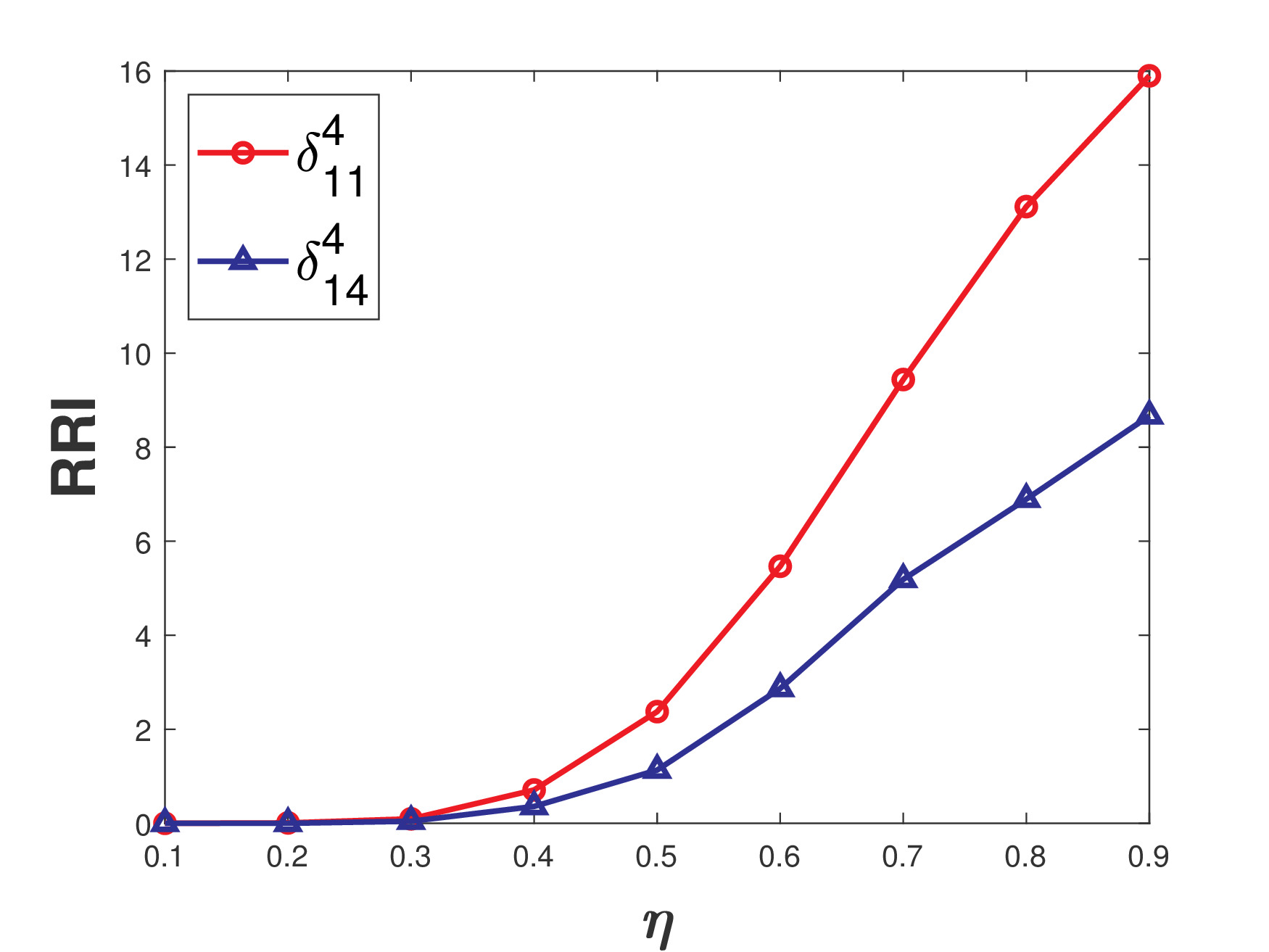}
			\caption{$a=-2$, $(p_1,p_2)=(6,9)$, $(\mu_1,\mu_2)=(0,0)$}
		\end{subfigure}
		\hfill
		\begin{subfigure}{0.32\textwidth}
			\centering
			\includegraphics[height=3.9cm, width=\textwidth]{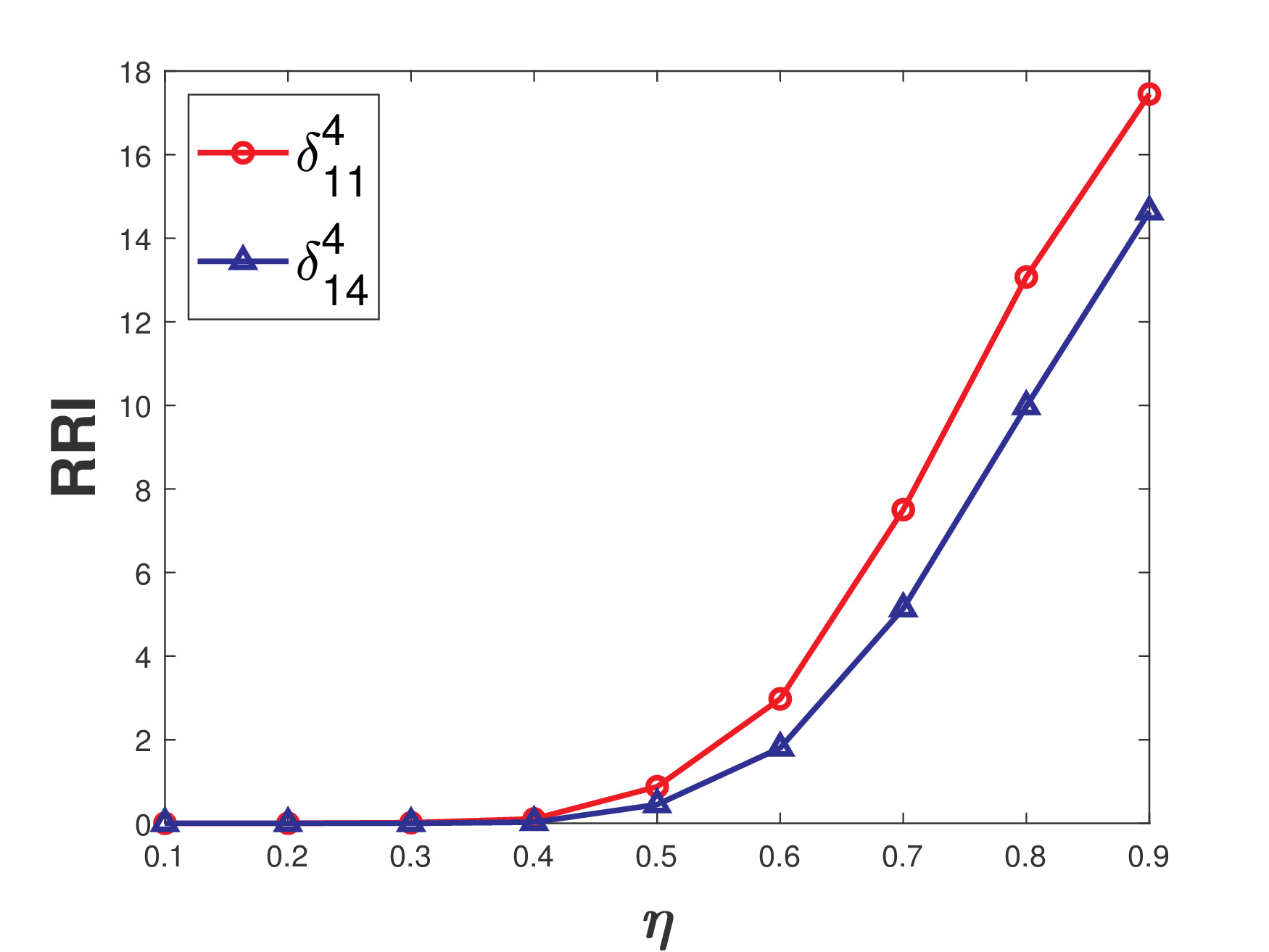}
			\caption{$a=-2$, $(p_1,p_2)=(10,14)$, $(\mu_1,\mu_2)=(1.5,2)$}
		\end{subfigure}
       \hfill
       	\begin{subfigure}{0.32\textwidth}
       	\centering
       	\includegraphics[height=3.9cm, width=\textwidth]{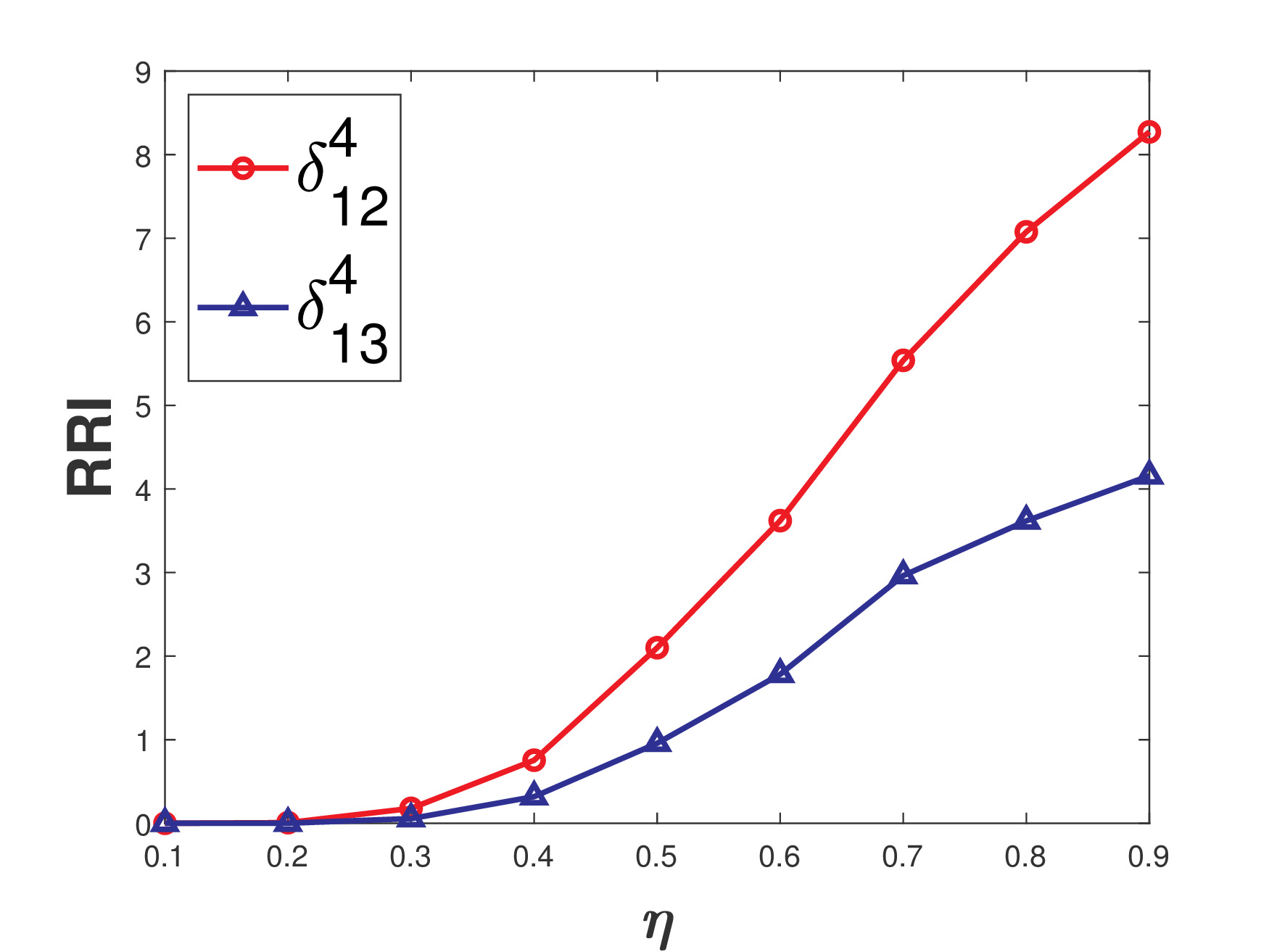}
       	\caption{$a=-2$, $(p_1,p_2)=(4,8)$, $(\mu_1,\mu_2)=(0,0)$}
       \end{subfigure}
       \hfill
       \begin{subfigure}{0.32\textwidth}
       	\centering
       	\includegraphics[height=3.9cm, width=\textwidth]{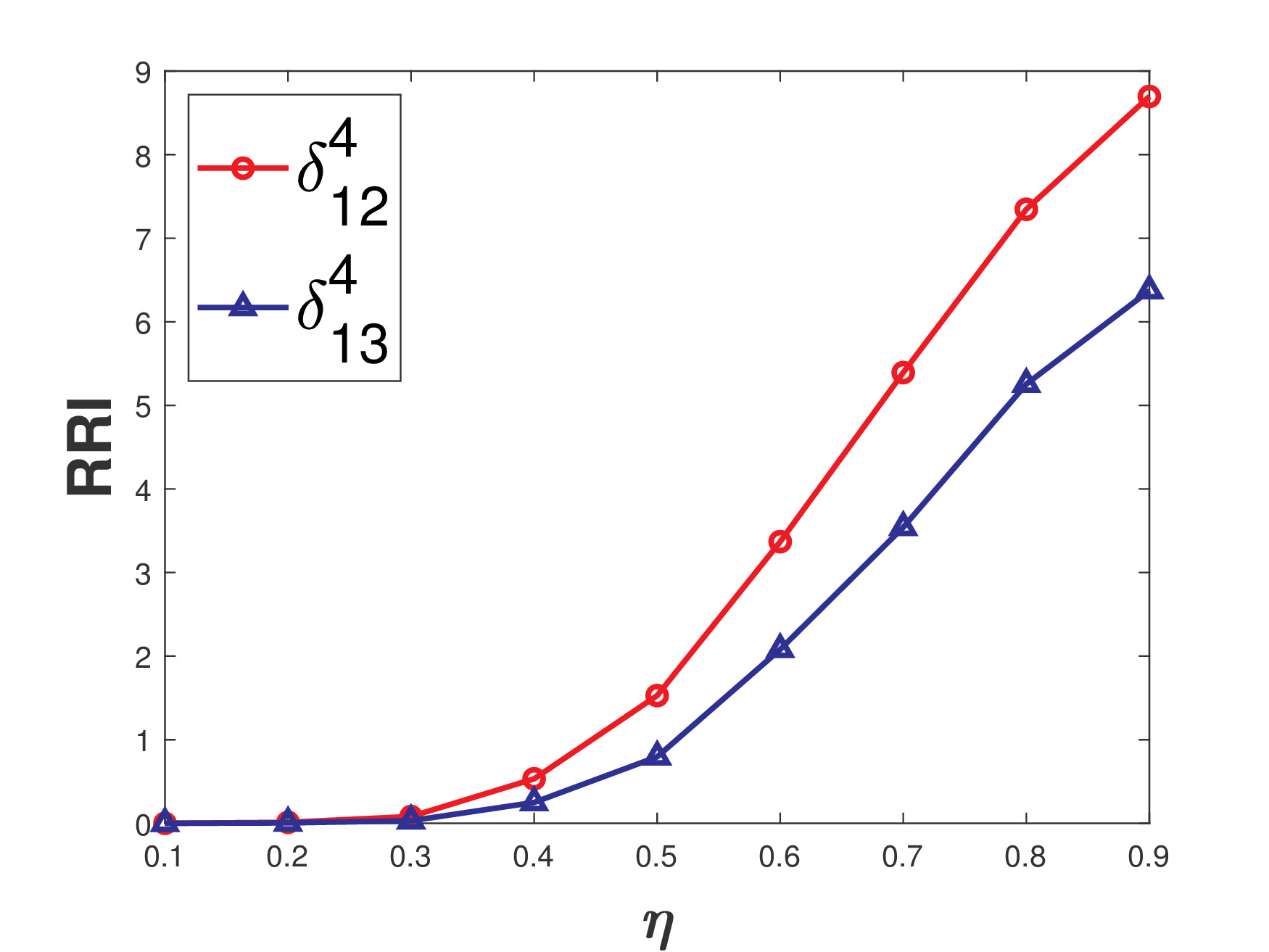}
       	\caption{$a=-2$, $(p_1,p_2)=(5,9)$, $(\mu_1,\mu_2)=(0,0.2)$}
       \end{subfigure}
       \hfill
       \begin{subfigure}{0.32\textwidth}
       	\centering
       	\includegraphics[height=3.9cm, width=\textwidth]{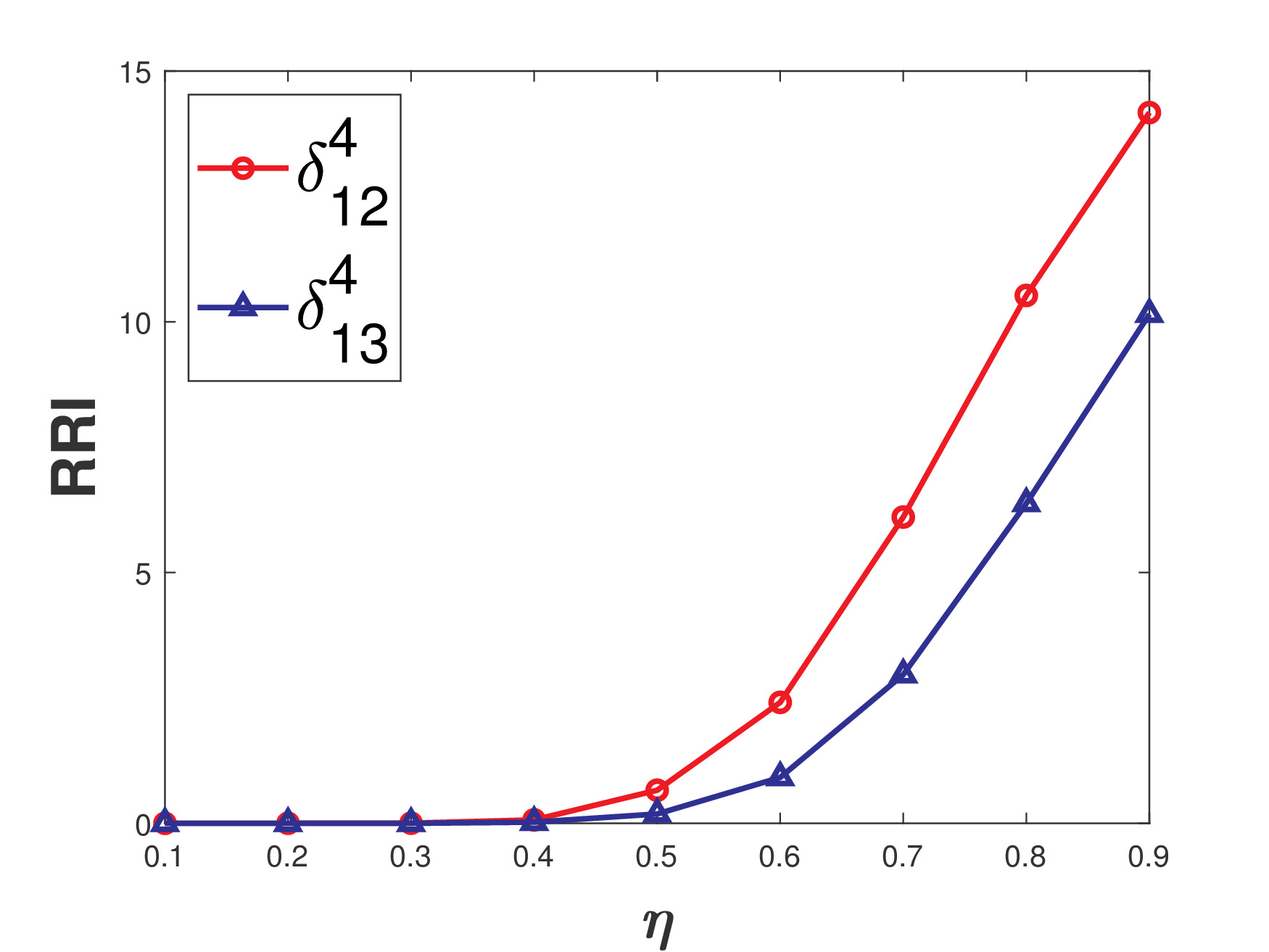}
       	\caption{$a=-2$, $(p_1,p_2)=(10,13)$, $(\mu_1,\mu_2)=(0.3,0.5)$}
       \end{subfigure}
	\end{subfigure}
	\caption{ RRI of different estimators with respect to BAEE for $\sigma_1^2$ under $L_4(t)$.}\label{fig4}
\end{figure}

\begin{figure}[htbp]
	\centering
	\begin{subfigure}{\textwidth}
		\centering
		\begin{subfigure}{0.32\textwidth}
			\centering
			\includegraphics[height=3.9cm, width=\textwidth]{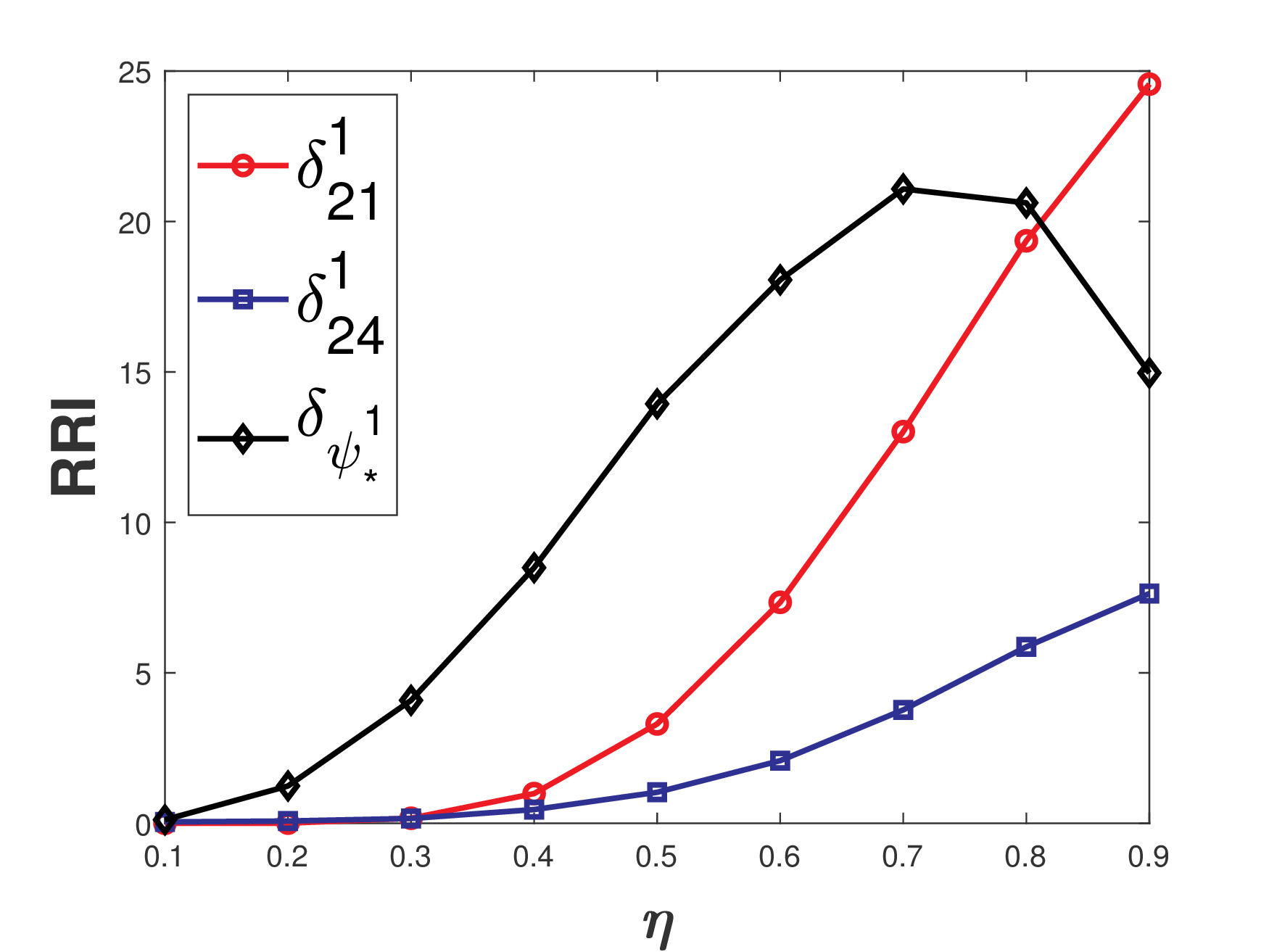}
			\caption{$(p_1,p_2)=(5,7)$, $(\mu_1,\mu_2)=(-0.5,-0.2)$}
		\end{subfigure}
		\hfill
		\begin{subfigure}{0.32\textwidth}
			\centering
			\includegraphics[height=3.9cm, width=\textwidth]{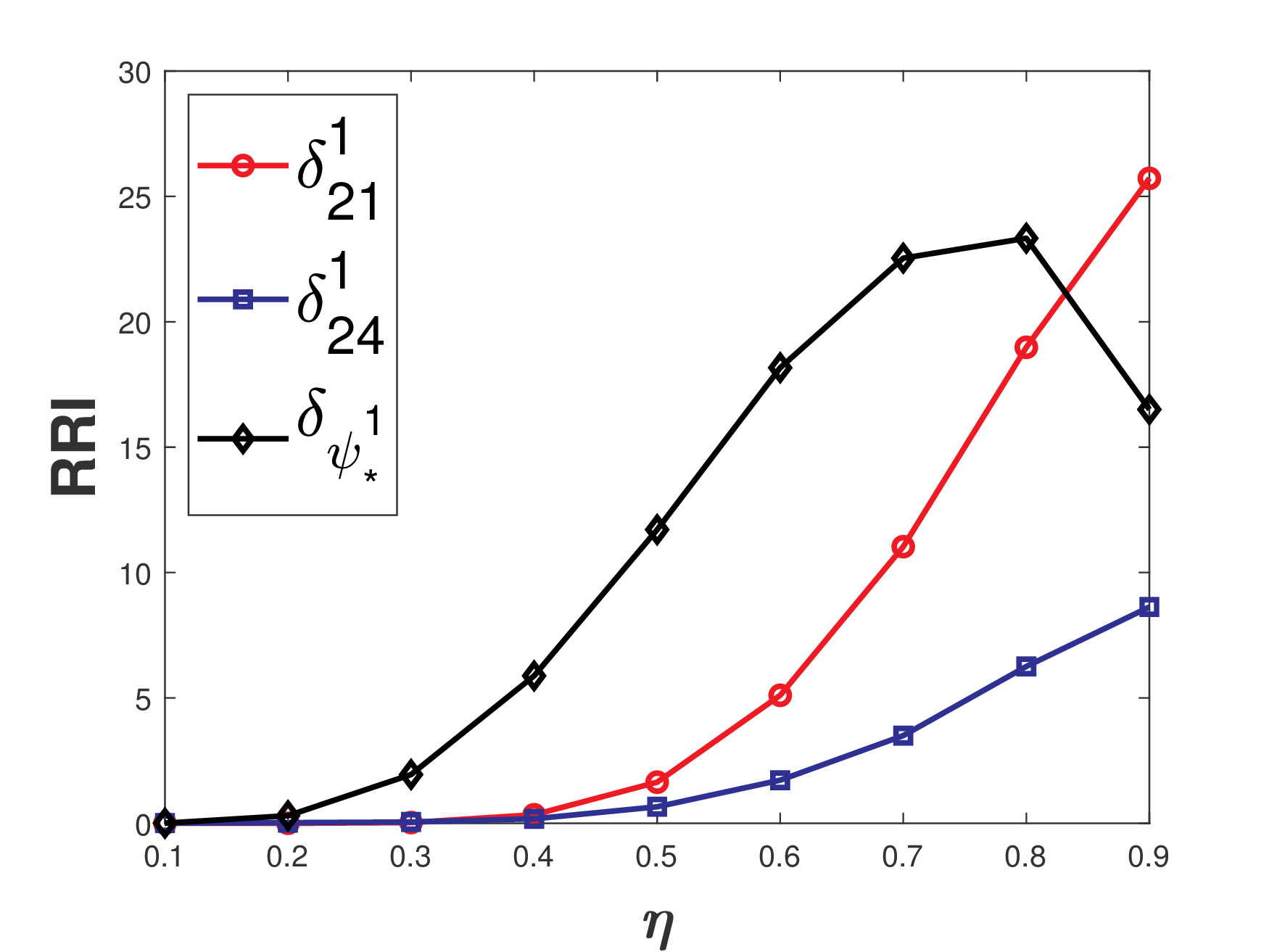}
			\caption{$(p_1,p_2)=(8,10)$, $(\mu_1,\mu_2)=(0,0.2)$}
		\end{subfigure}
		\hfill
		\begin{subfigure}{0.32\textwidth}
			\centering
			\includegraphics[height=3.9cm, width=\textwidth]{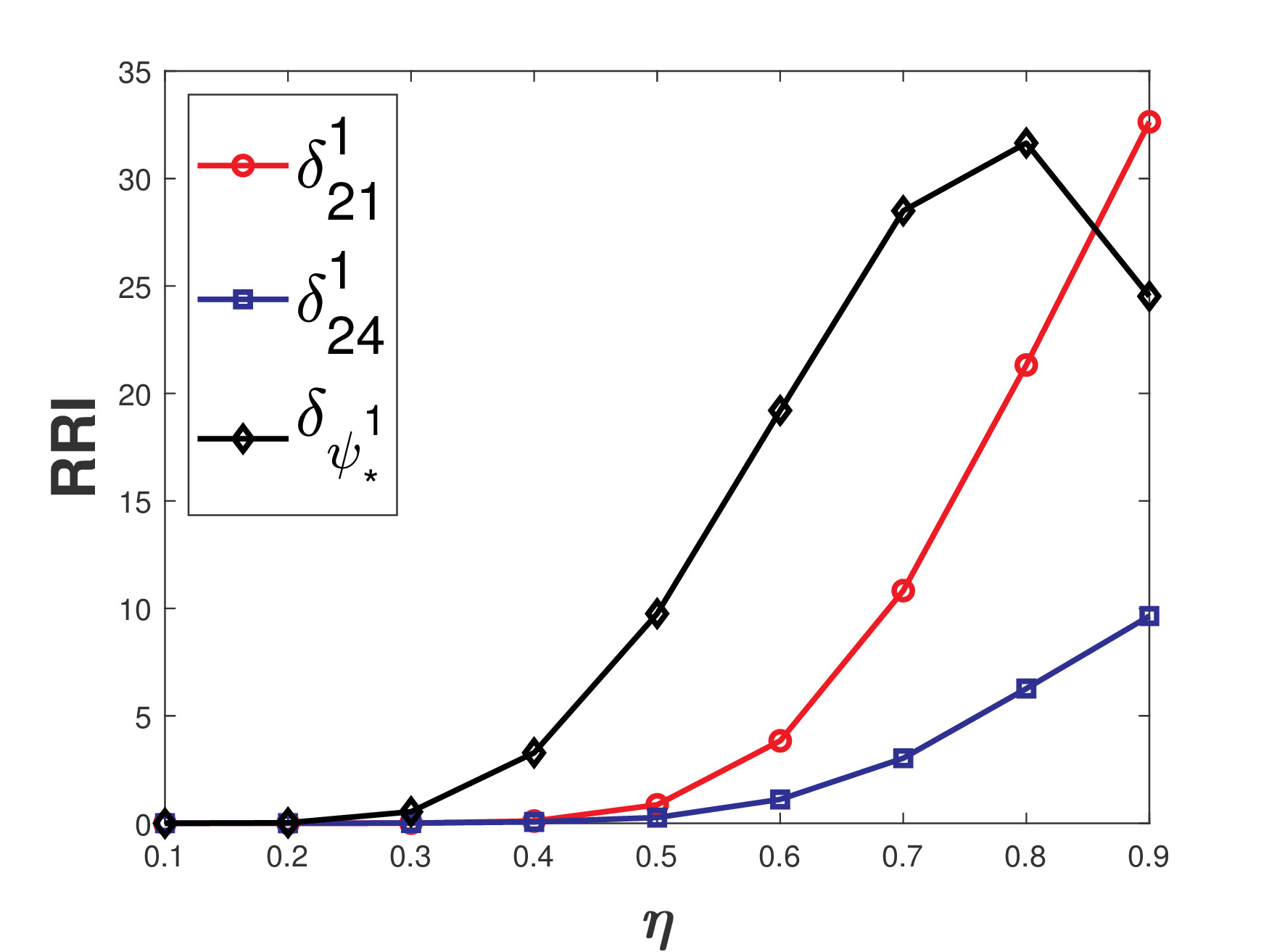}
			\caption{$(p_1,p_2)=(16,12)$, $(\mu_1,\mu_2)=(0.3,0.5)$}
		\end{subfigure}
        \hfill
        	\begin{subfigure}{0.32\textwidth}
        	\centering
        	\includegraphics[height=3.9cm, width=\textwidth]{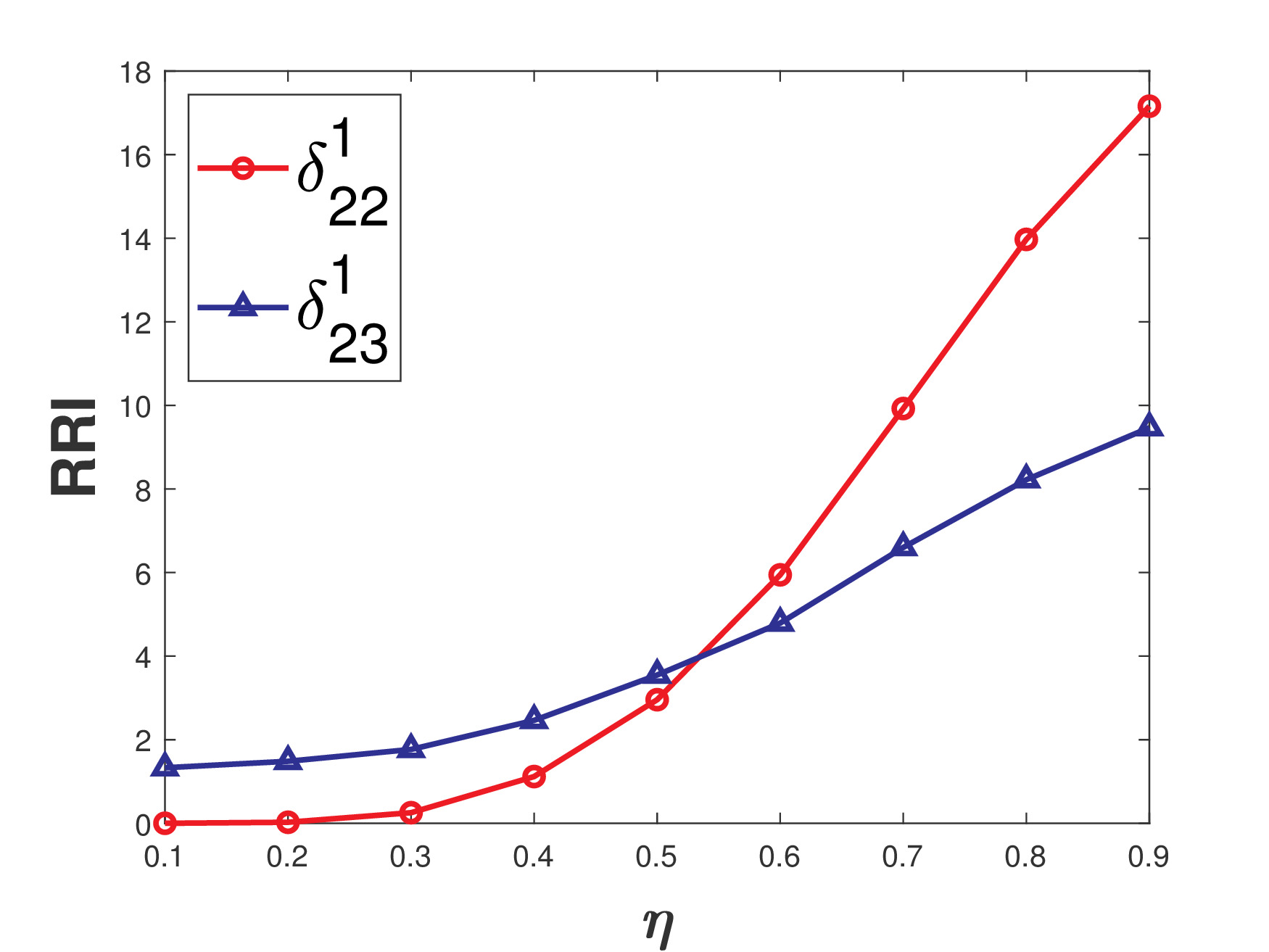}
        	\caption{$(p_1,p_2)=(5,5)$, $(\mu_1,\mu_2)=(0,0)$}
        \end{subfigure}
        \hfill
        \begin{subfigure}{0.32\textwidth}
        	\centering
        	\includegraphics[height=3.9cm, width=\textwidth]{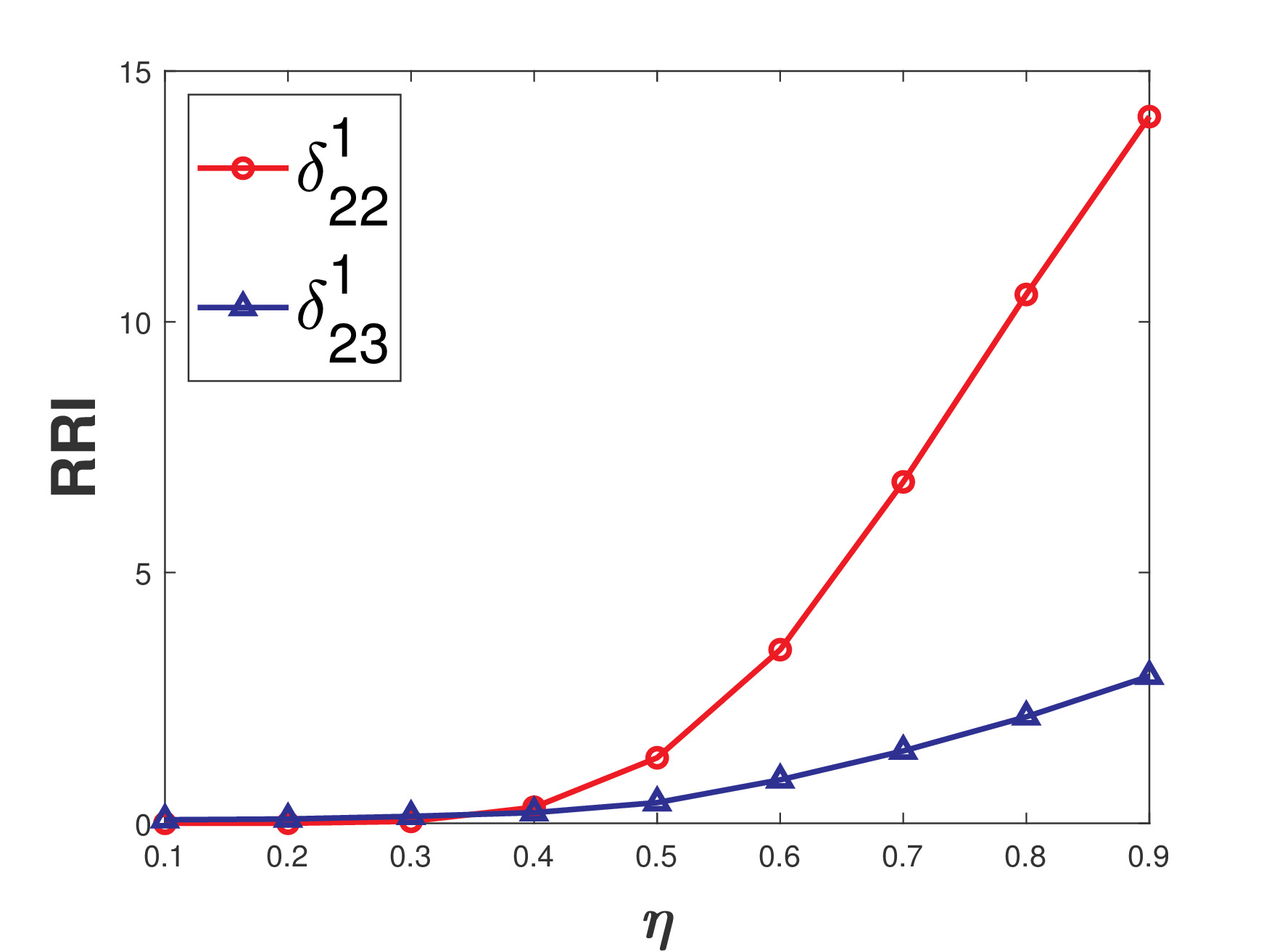}
        	\caption{$(p_1,p_2)=(6,8)$, $(\mu_1,\mu_2)=(0,0.3)$}
        \end{subfigure}
        \hfill
        \begin{subfigure}{0.32\textwidth}
        	\centering
        	\includegraphics[height=3.9cm, width=\textwidth]{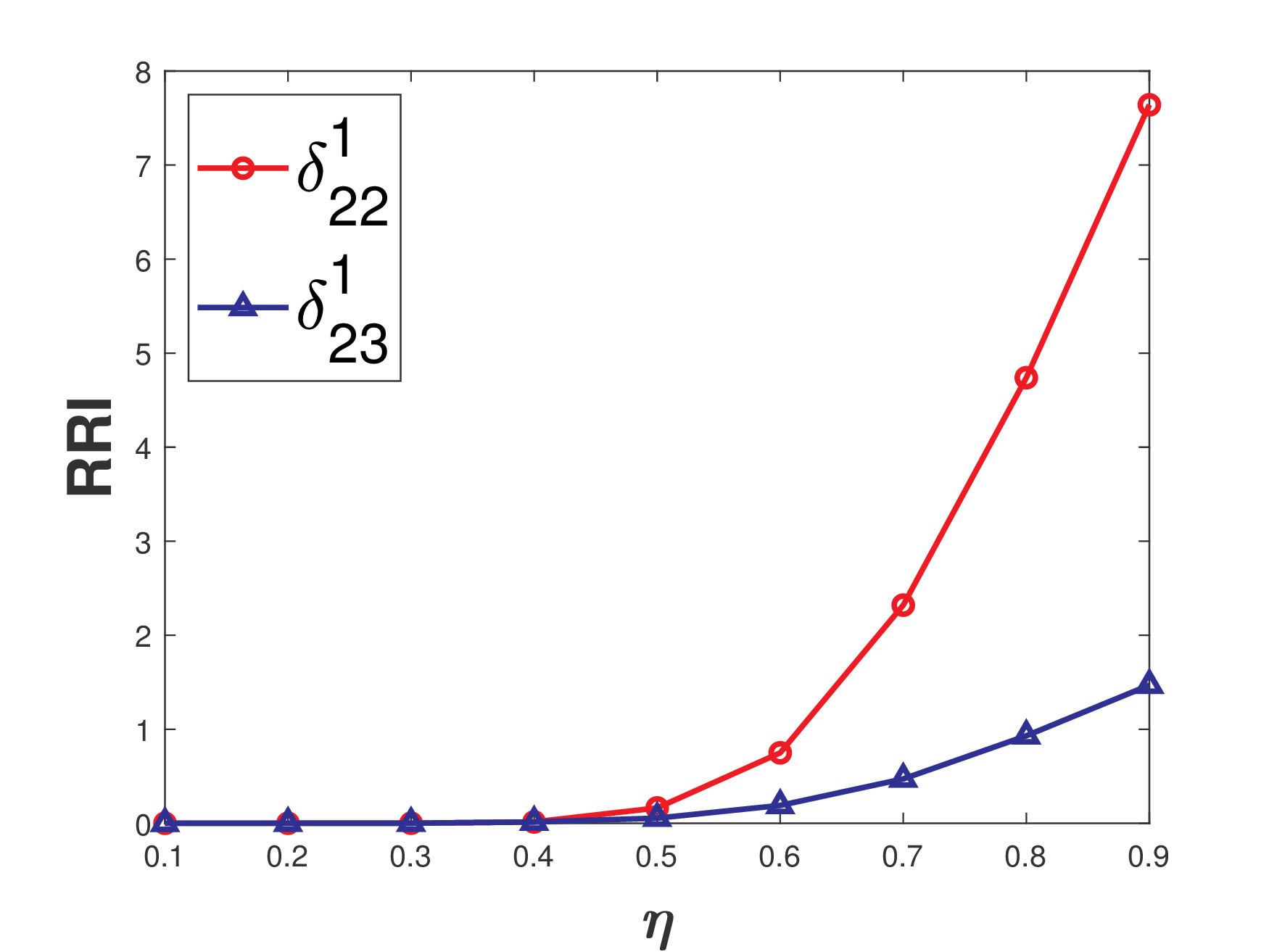}
        	\caption{$(p_1,p_2)=(9,12)$, $(\mu_1,\mu_2)=(0.15,0.25)$}
        \end{subfigure}
	\end{subfigure}
	\caption{RRI of different estimators with respect to BAEE for $\sigma_2^2$ under $L_1(t)$.}\label{fig5}
\end{figure}

\begin{figure}[htbp]
	\centering
	\begin{subfigure}{\textwidth}
	\centering
	\begin{subfigure}{0.32\textwidth}
		\centering
		\includegraphics[height=3.9cm, width=\textwidth]{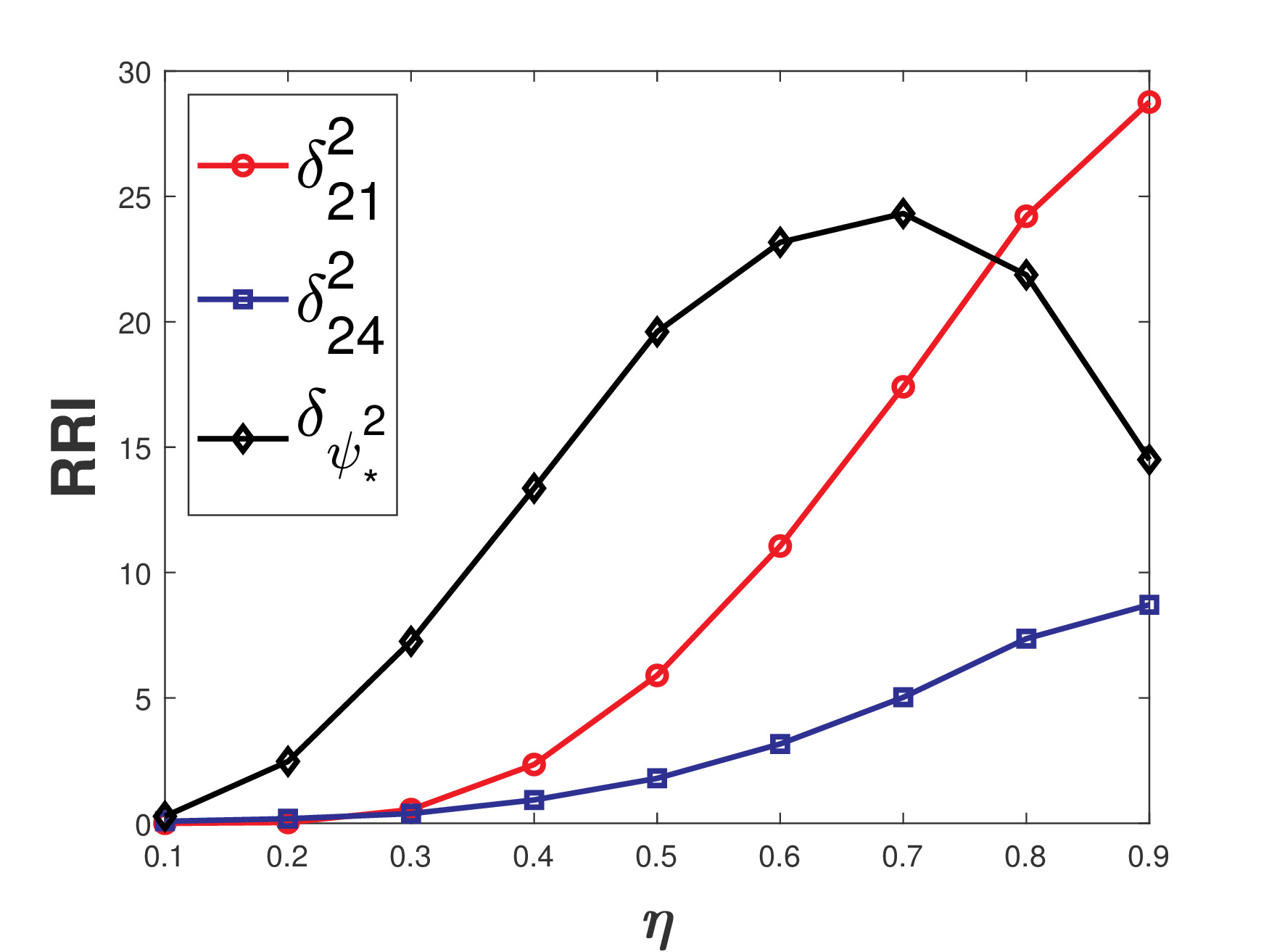}
		\caption{$(p_1,p_2)=(5,7)$, $(\mu_1,\mu_2)=(-0.5,-0.2)$}
	\end{subfigure}
	\hfill
	\begin{subfigure}{0.32\textwidth}
		\centering
		\includegraphics[height=3.9cm, width=\textwidth]{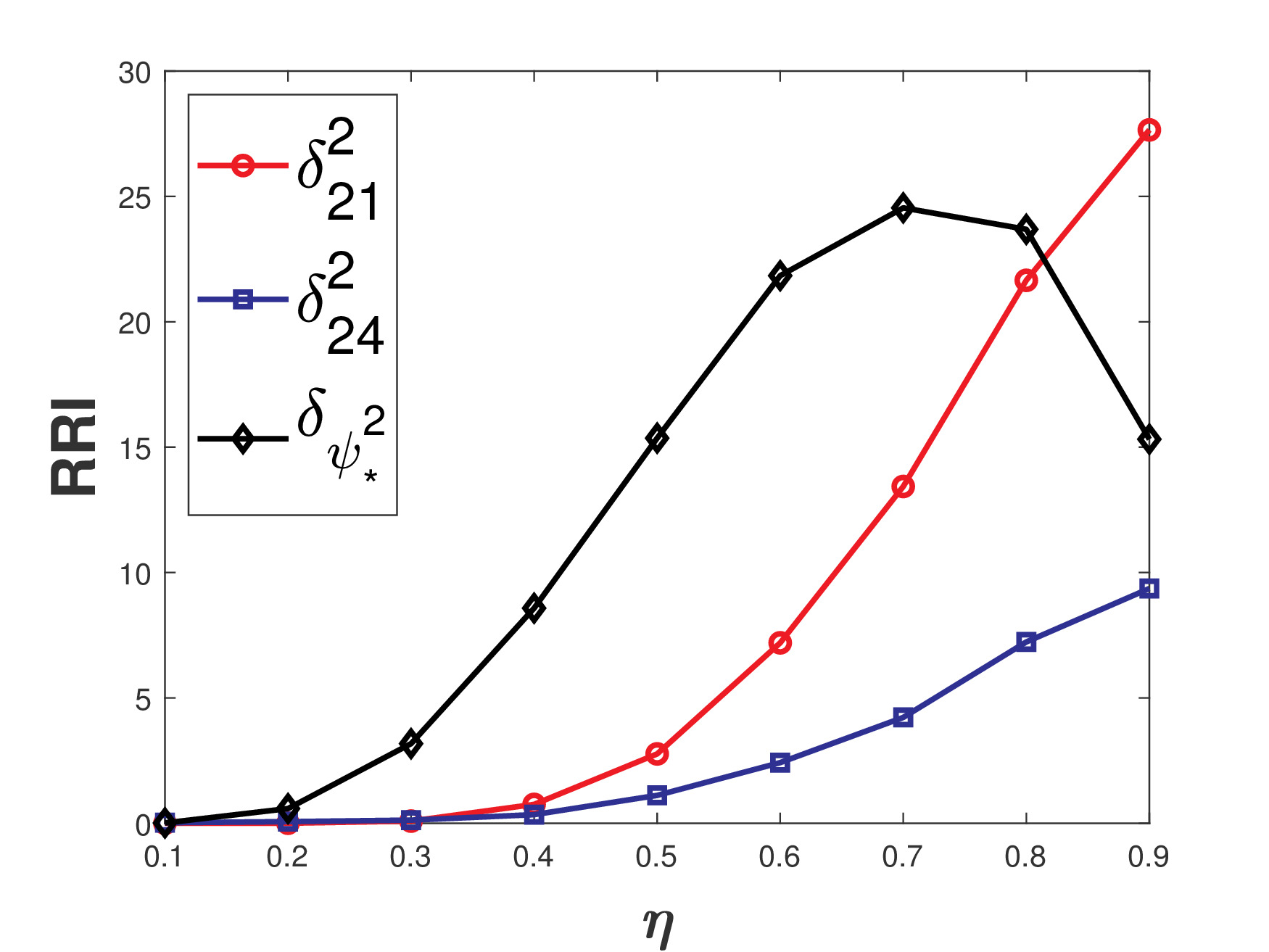}
		\caption{$(p_1,p_2)=(8,10)$, $(\mu_1,\mu_2)=(0,0.2)$}
	\end{subfigure}
	\hfill
	\begin{subfigure}{0.32\textwidth}
		\centering
		\includegraphics[height=3.9cm, width=\textwidth]{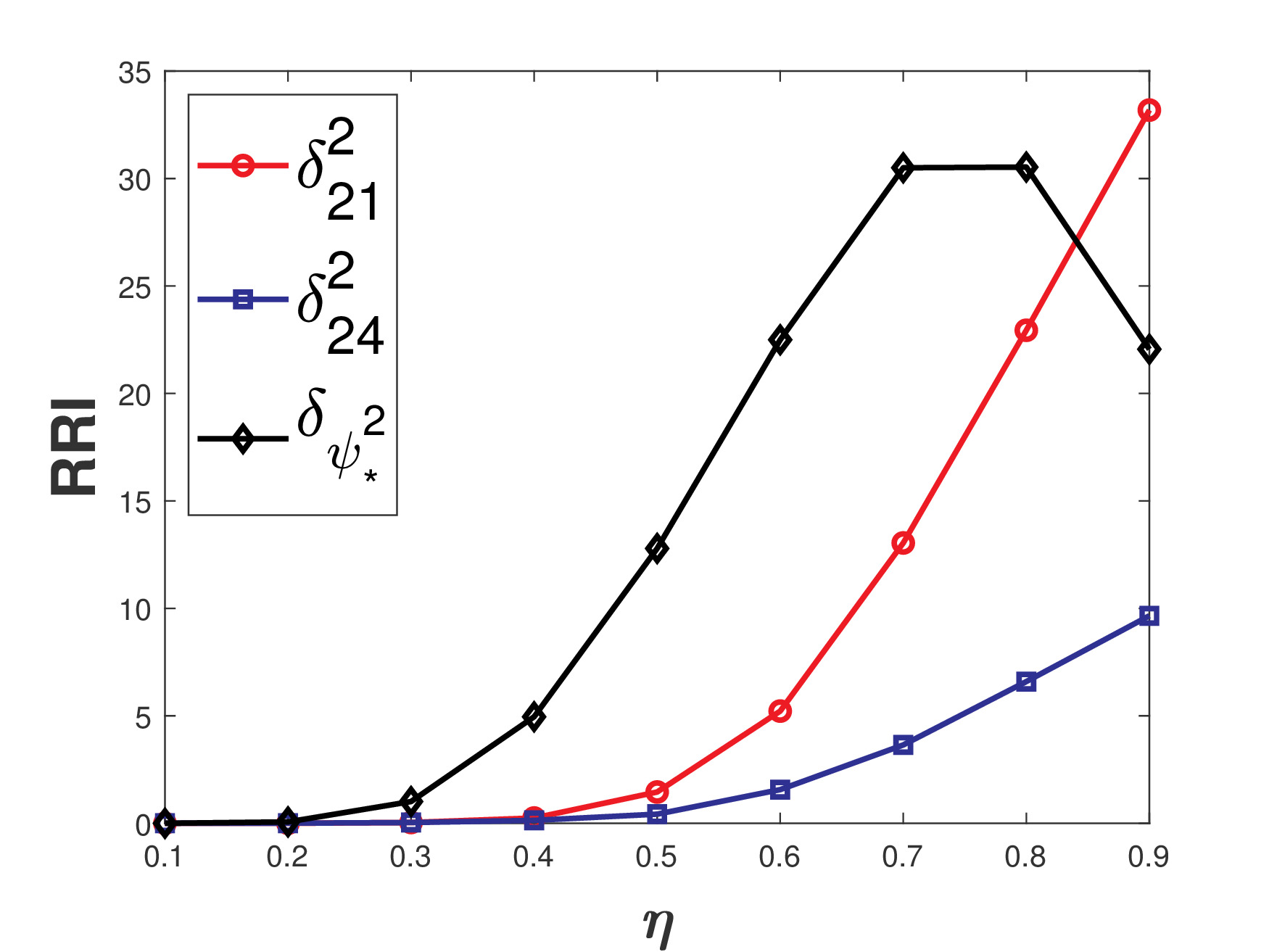}
		\caption{$(p_1,p_2)=(16,12)$, $(\mu_1,\mu_2)=(0.3,0.5)$}
	\end{subfigure}
	\hfill
	\begin{subfigure}{0.32\textwidth}
		\centering
		\includegraphics[height=3.9cm, width=\textwidth]{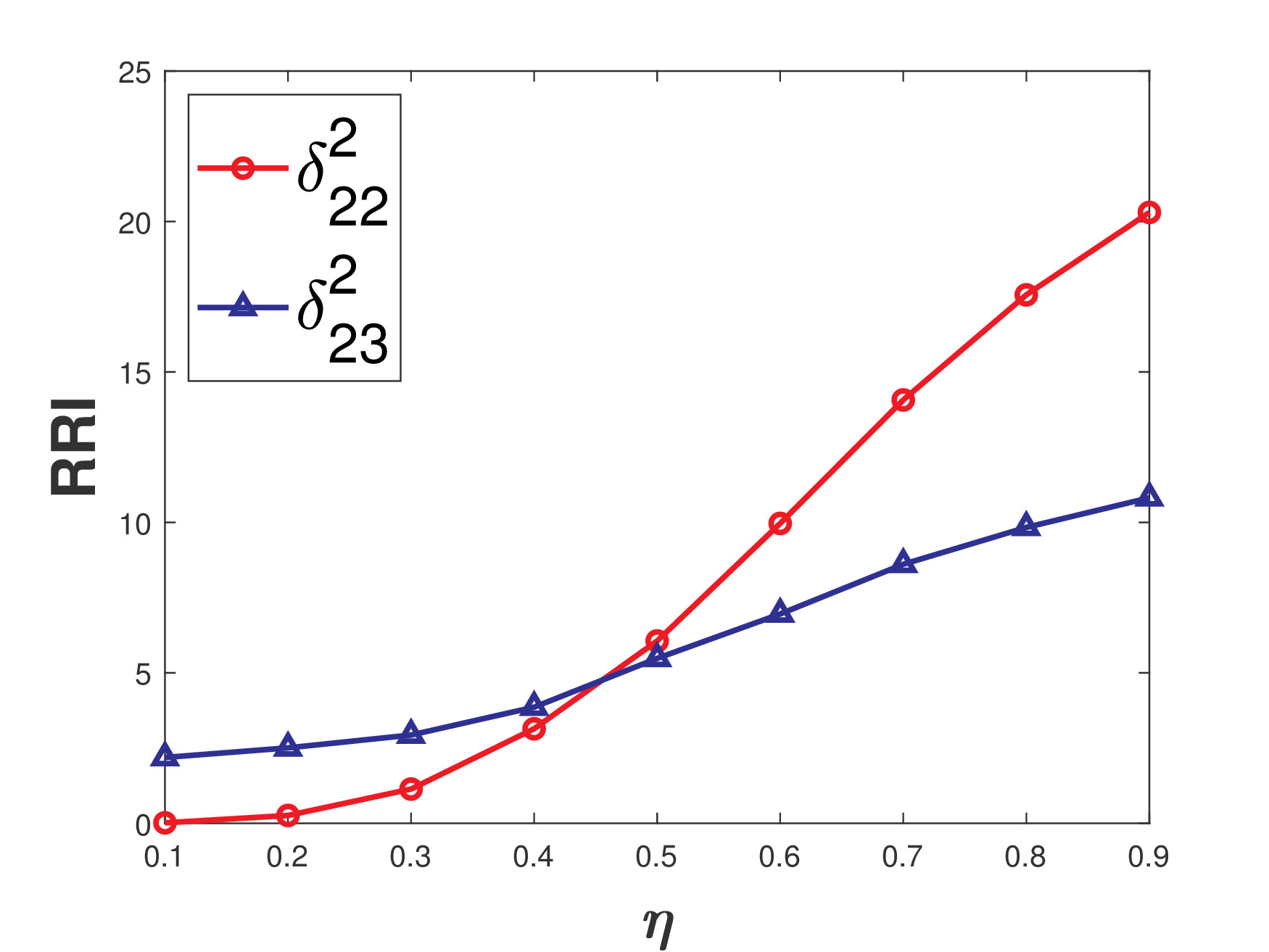}
		\caption{$(p_1,p_2)=(5,5)$, $(\mu_1,\mu_2)=(0,0)$}
	\end{subfigure}
	\hfill
	\begin{subfigure}{0.32\textwidth}
		\centering
		\includegraphics[height=3.9cm, width=\textwidth]{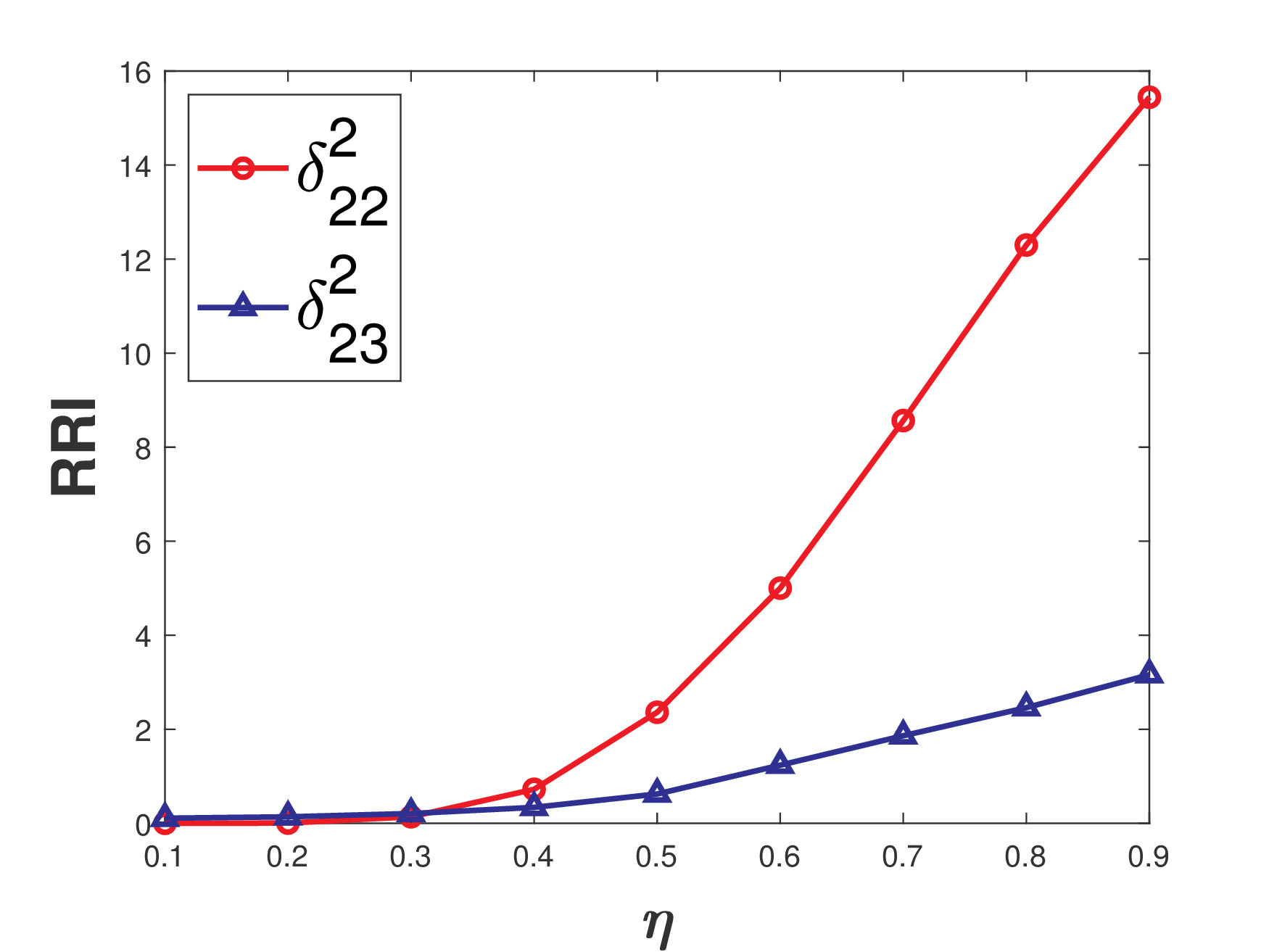}
		\caption{$(p_1,p_2)=(6,8)$, $(\mu_1,\mu_2)=(0,0.3)$}
	\end{subfigure}
	\hfill
	\begin{subfigure}{0.32\textwidth}
		\centering
		\includegraphics[height=3.9cm, width=\textwidth]{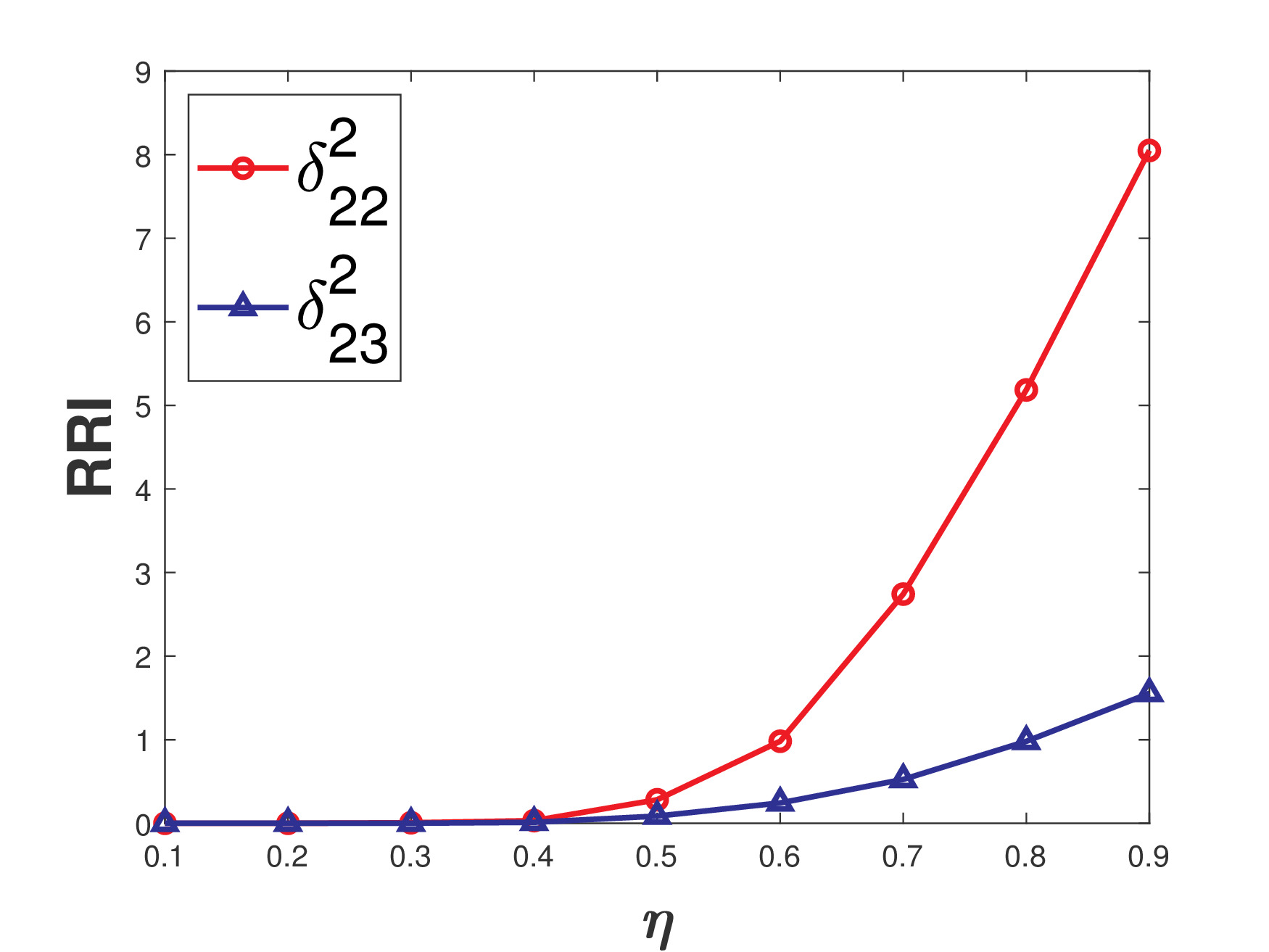}
		\caption{$(p_1,p_2)=(9,12)$, $(\mu_1,\mu_2)=(0.15,0.25)$}
	\end{subfigure}
    \end{subfigure}
\caption{RRI of different estimators with respect to BAEE for $\sigma_2^2$ under $L_2(t)$.}\label{fig6}
\end{figure}

\begin{figure}[htbp]

\begin{subfigure}{\textwidth}
	\centering
	\begin{subfigure}{0.32\textwidth}
		\centering
		\includegraphics[height=3.9cm, width=\textwidth]{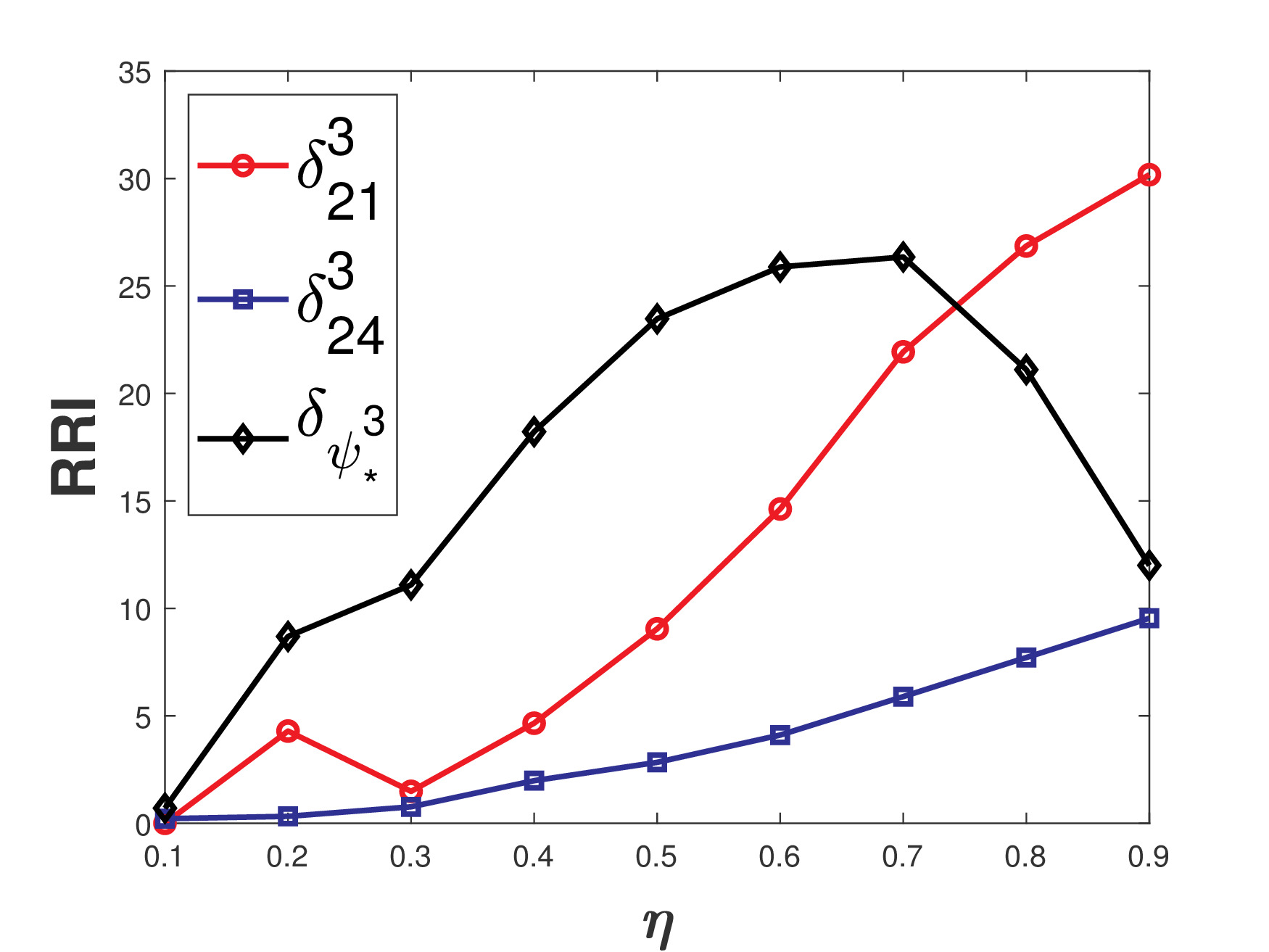}
		\caption{$(p_1,p_2)=(5,7)$, $(\mu_1,\mu_2)=(-0.5,-0.2)$}
	\end{subfigure}
	\hfill
	\begin{subfigure}{0.32\textwidth}
		\centering
		\includegraphics[height=3.9cm, width=\textwidth]{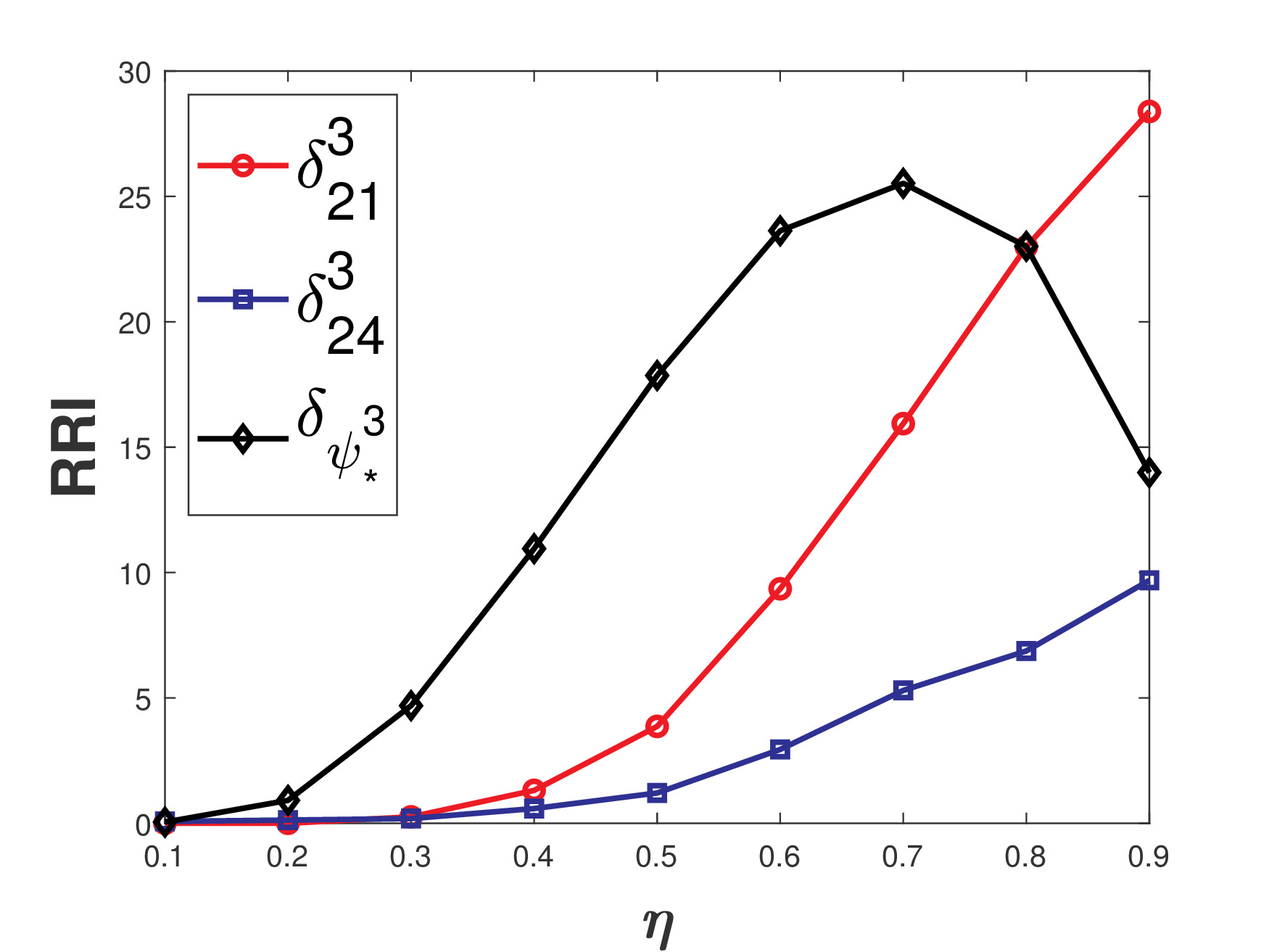}
		\caption{$(p_1,p_2)=(8,10)$, $(\mu_1,\mu_2)=(0,0.2)$}
	\end{subfigure}
	\hfill
	\begin{subfigure}{0.32\textwidth}
		\centering
		\includegraphics[height=3.9cm, width=\textwidth]{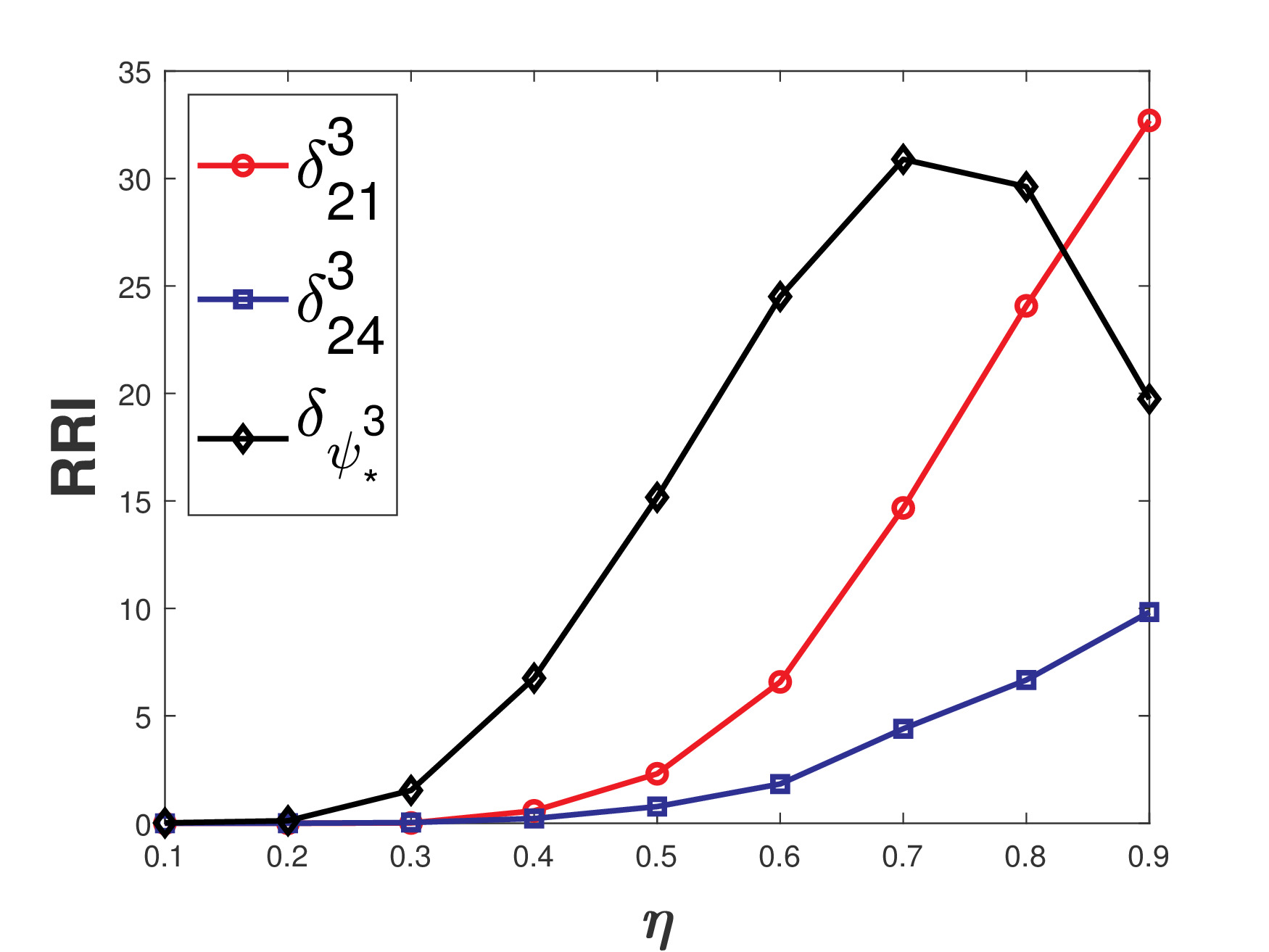}
		\caption{$(p_1,p_2)=(16,12)$, $(\mu_1,\mu_2)=(0.3,0.5)$}
	\end{subfigure}
	\hfill
	\begin{subfigure}{0.32\textwidth}
		\centering
		\includegraphics[height=3.9cm, width=\textwidth]{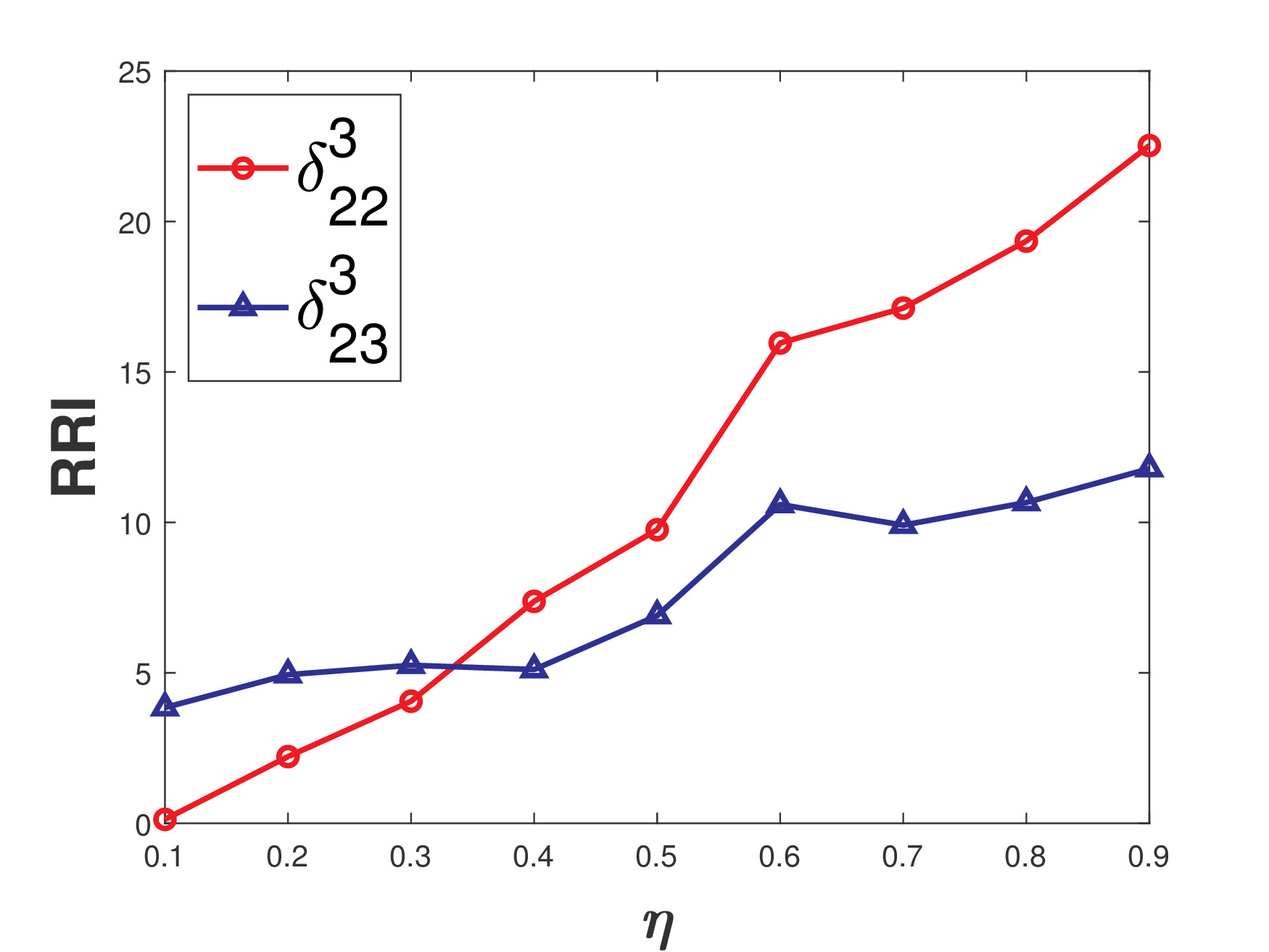}
		\caption{$(p_1,p_2)=(5,5)$, $(\mu_1,\mu_2)=(0,0)$}
	\end{subfigure}
	\hfill
	\begin{subfigure}{0.32\textwidth}
		\centering
		\includegraphics[height=3.9cm, width=\textwidth]{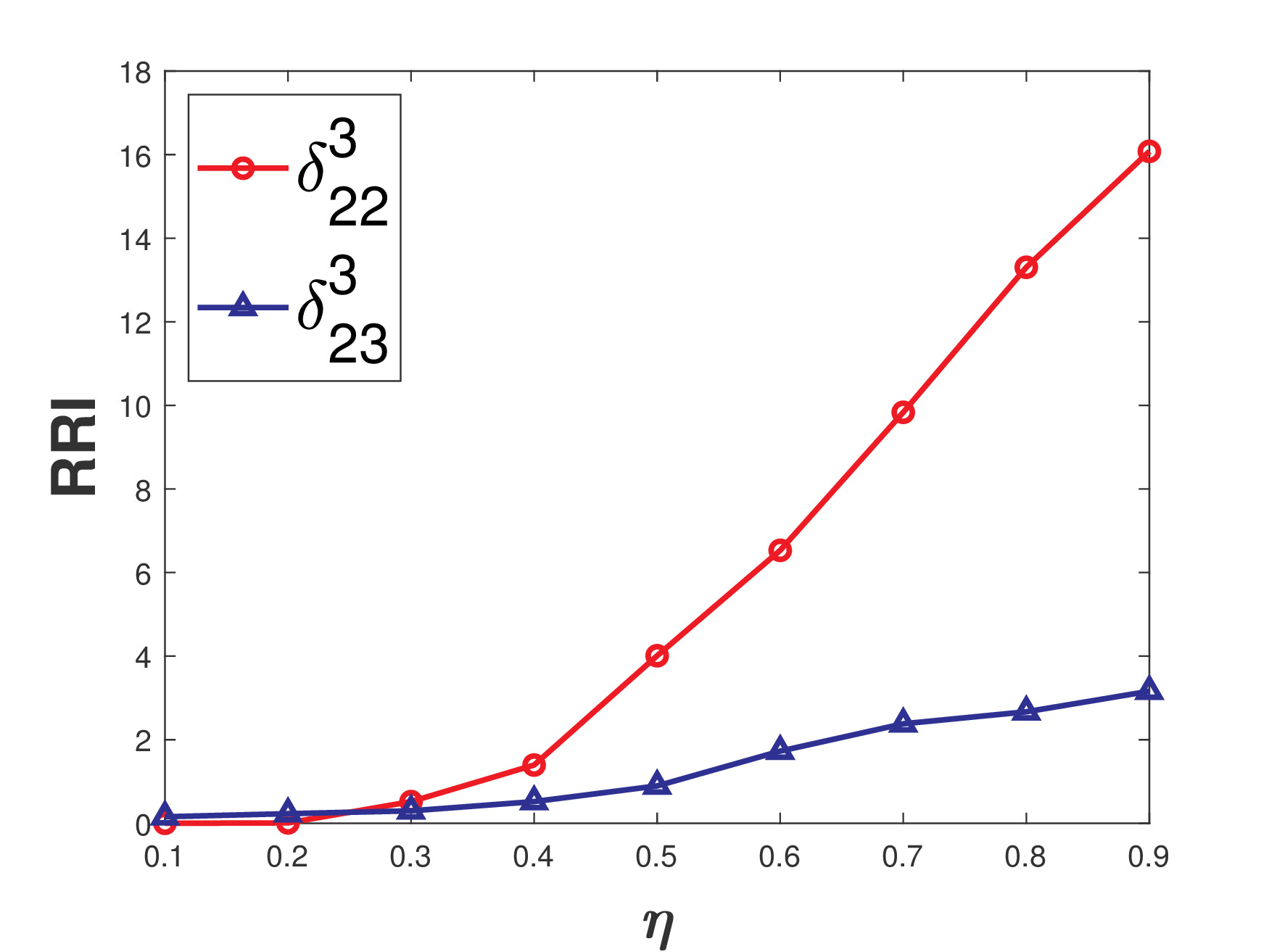}
		\caption{$(p_1,p_2)=(6,8)$, $(\mu_1,\mu_2)=(0,0.3)$}
	\end{subfigure}
	\hfill
	\begin{subfigure}{0.32\textwidth}
		\centering
		\includegraphics[height=3.9cm, width=\textwidth]{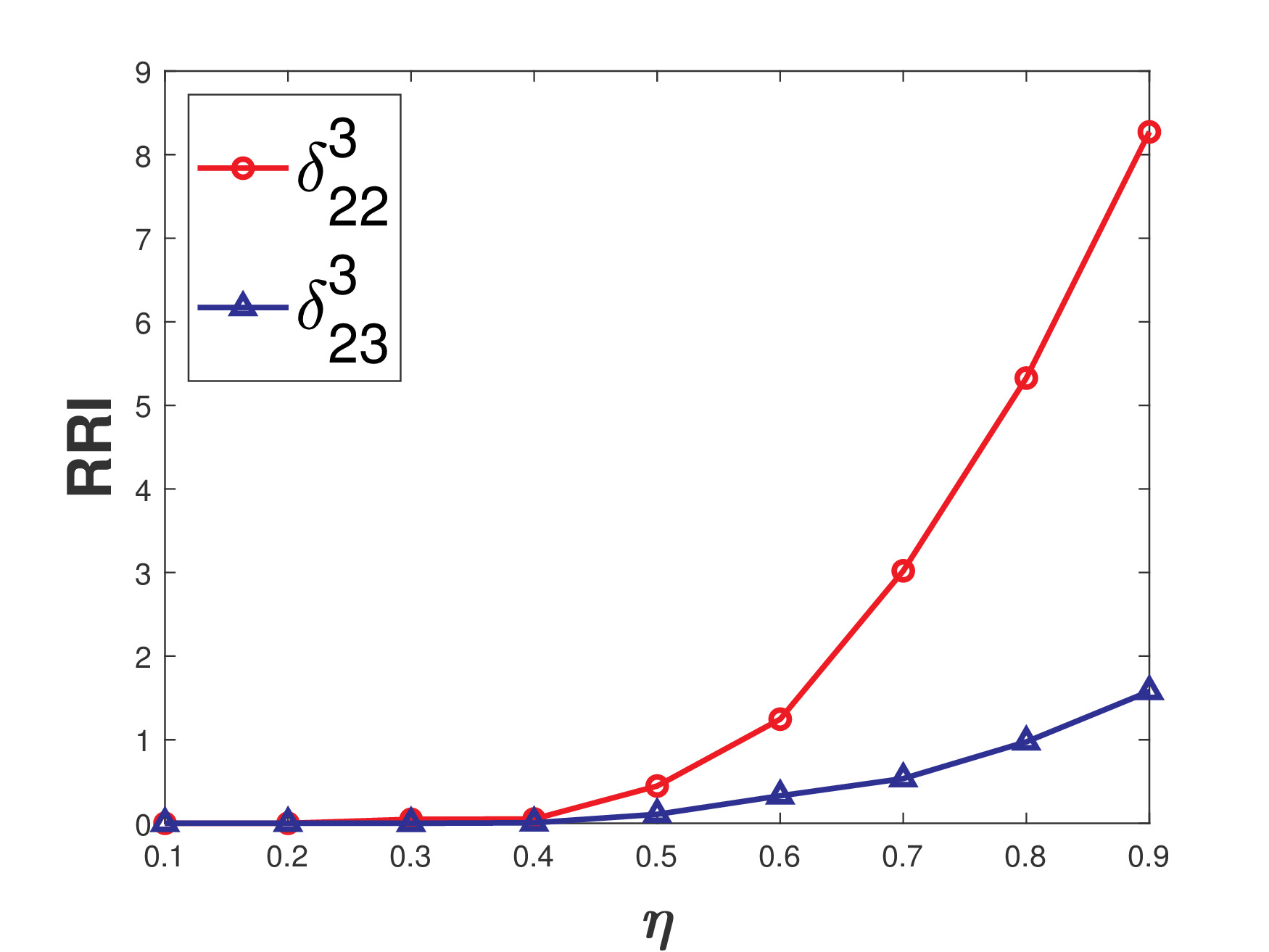}
		\caption{$(p_1,p_2)=(9,12)$, $(\mu_1,\mu_2)=(0.15,0.25)$}
	\end{subfigure}
	\end{subfigure}
	\caption{RRI of different estimators with respect to BAEE for $\sigma_2^2$ under $L_3(t)$.}\label{fig7}
\end{figure}

\begin{figure}[h!]
	
\centering
	\begin{subfigure}{\textwidth}
		\centering
		\begin{subfigure}{0.32\textwidth}
			\centering
			\includegraphics[height=3.9cm, width=\textwidth]{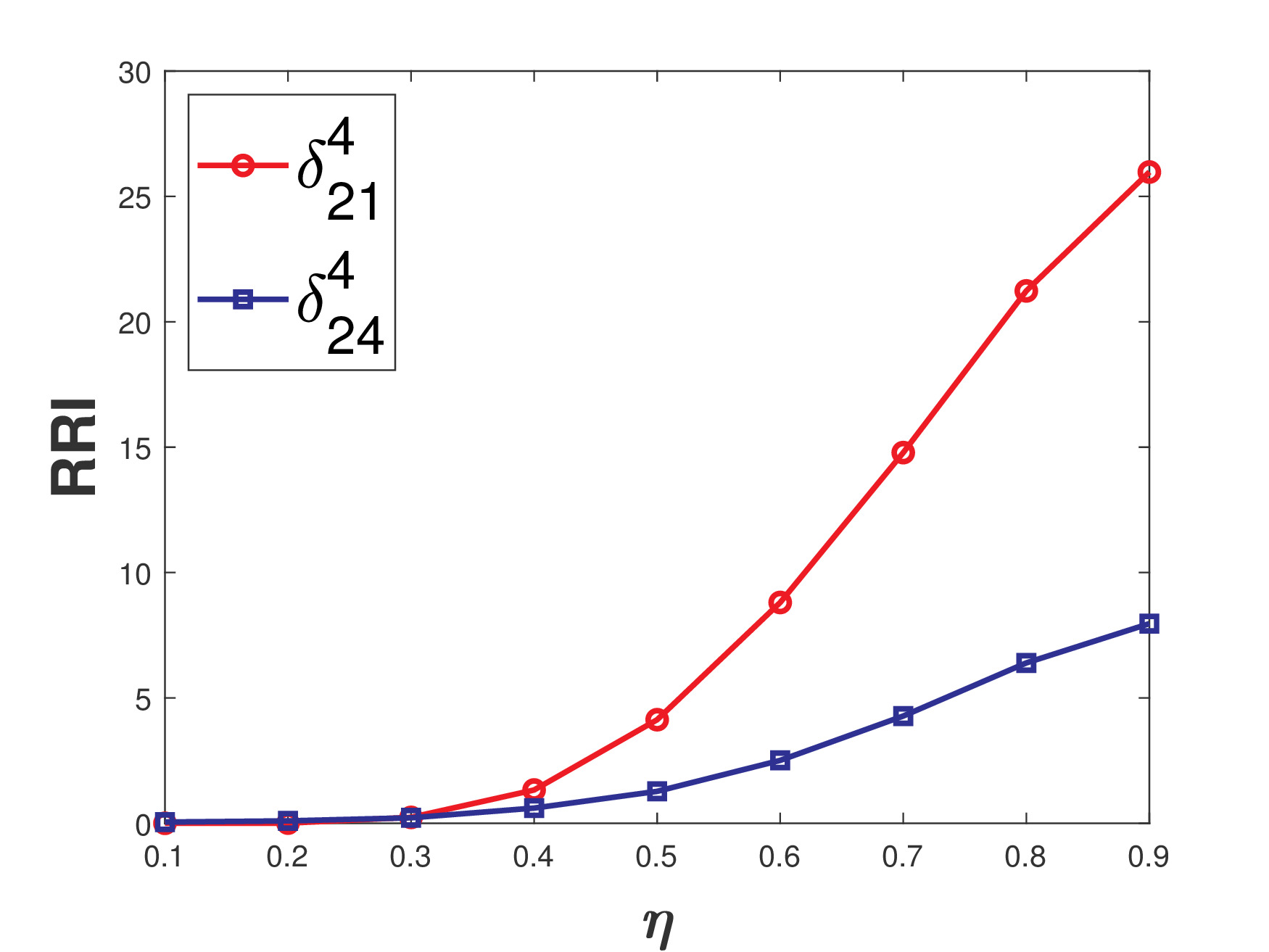}
			\caption{$a=-2$, $(p_1,p_2)=(5,7)$, $(\mu_1,\mu_2)=(-0.5,-0.2)$}
		\end{subfigure}
		\hfill
		\begin{subfigure}{0.32\textwidth}
			\centering
			\includegraphics[height=3.9cm, width=\textwidth]{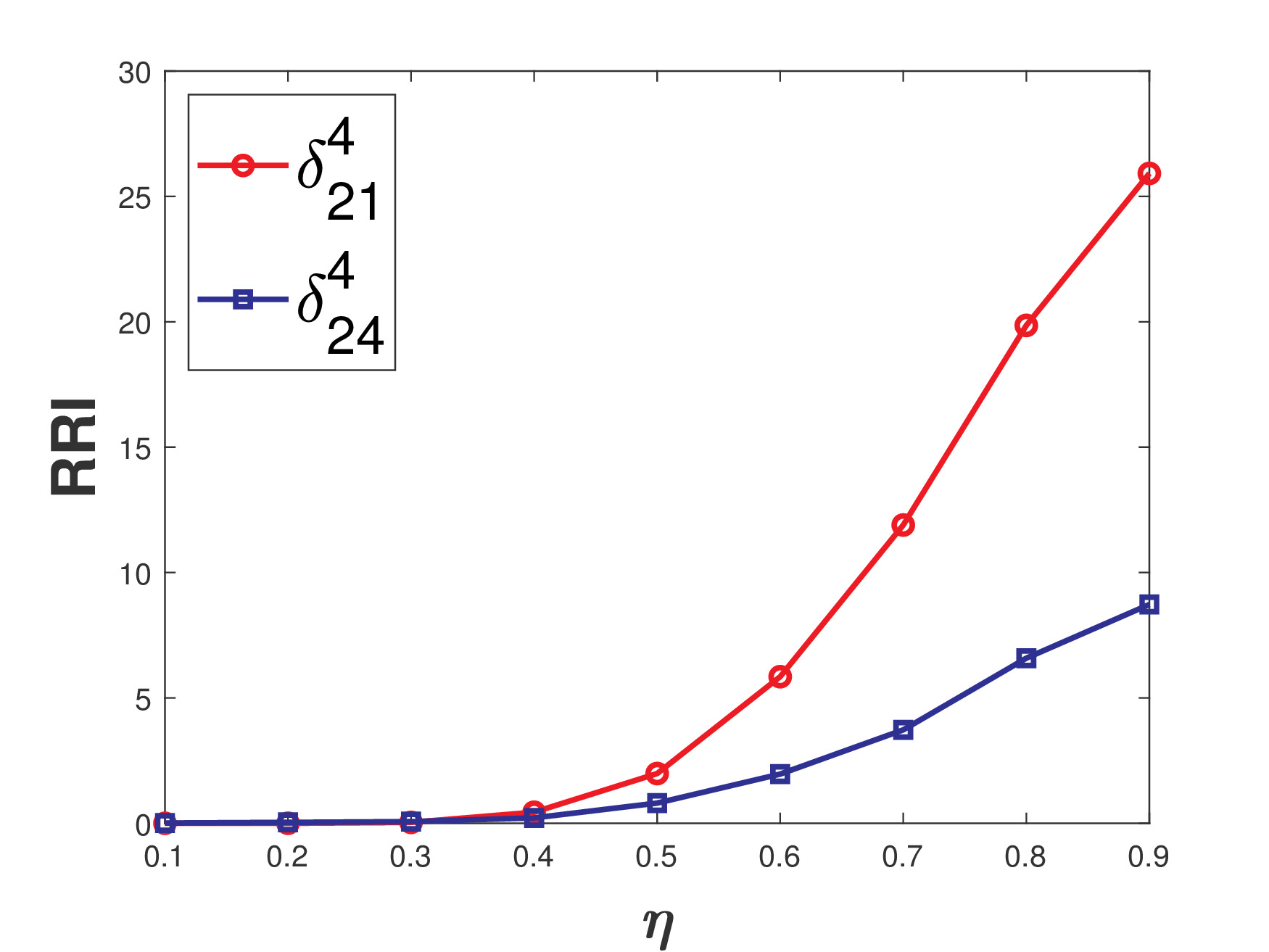}
			\caption{$a=-2$, $(p_1,p_2)=(8,10)$, $(\mu_1,\mu_2)=(0,0.2)$}
		\end{subfigure}
		\hfill
		\begin{subfigure}{0.32\textwidth}
			\centering
			\includegraphics[height=3.9cm, width=\textwidth]{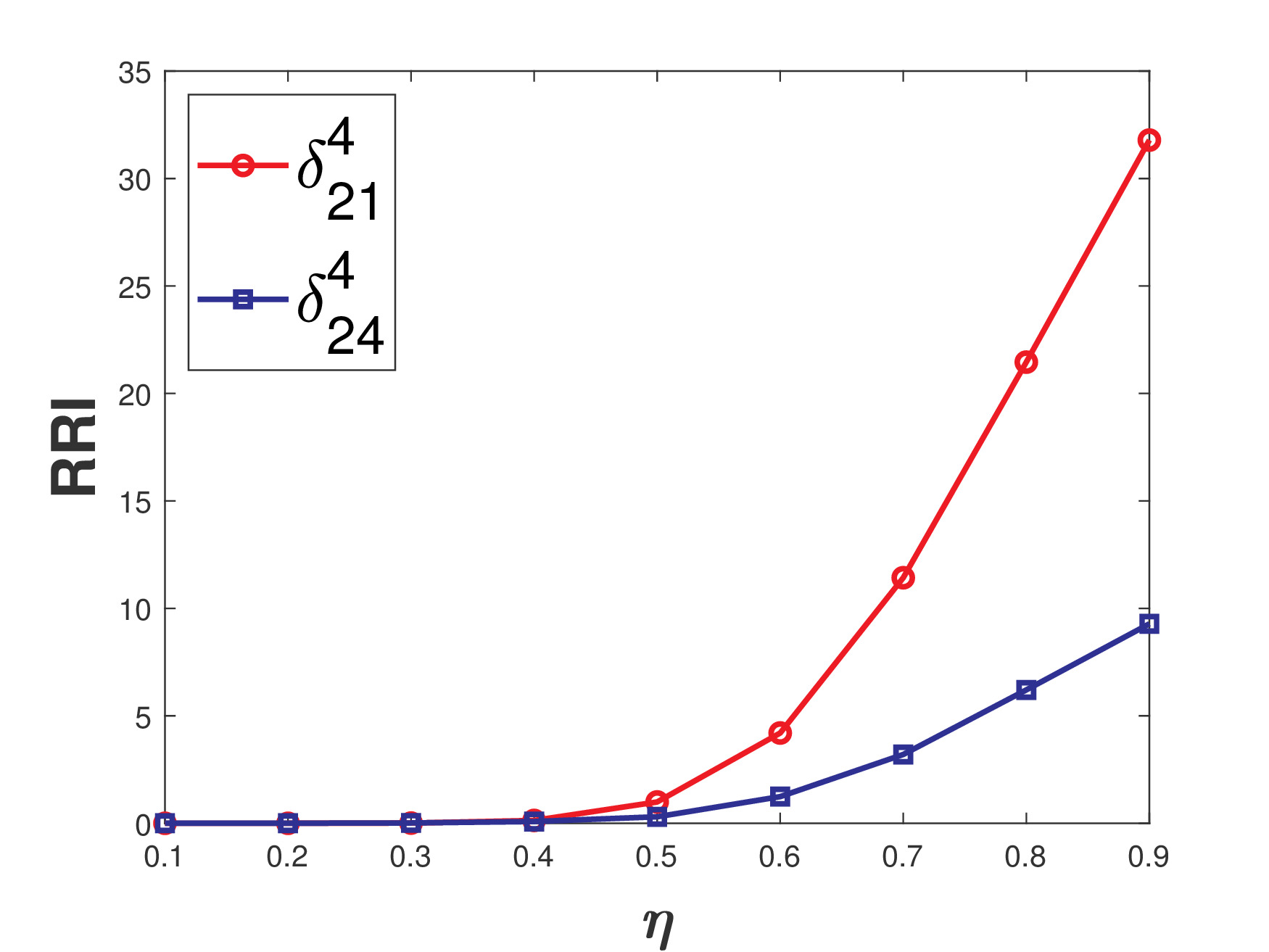}
			\caption{$a=-2$, $(p_1,p_2)=(16,12)$, $(\mu_1,\mu_2)=(0.3,0.5)$}
		\end{subfigure}
		\hfill	\begin{subfigure}{0.32\textwidth}
			\centering
			\includegraphics[height=3.9cm, width=\textwidth]{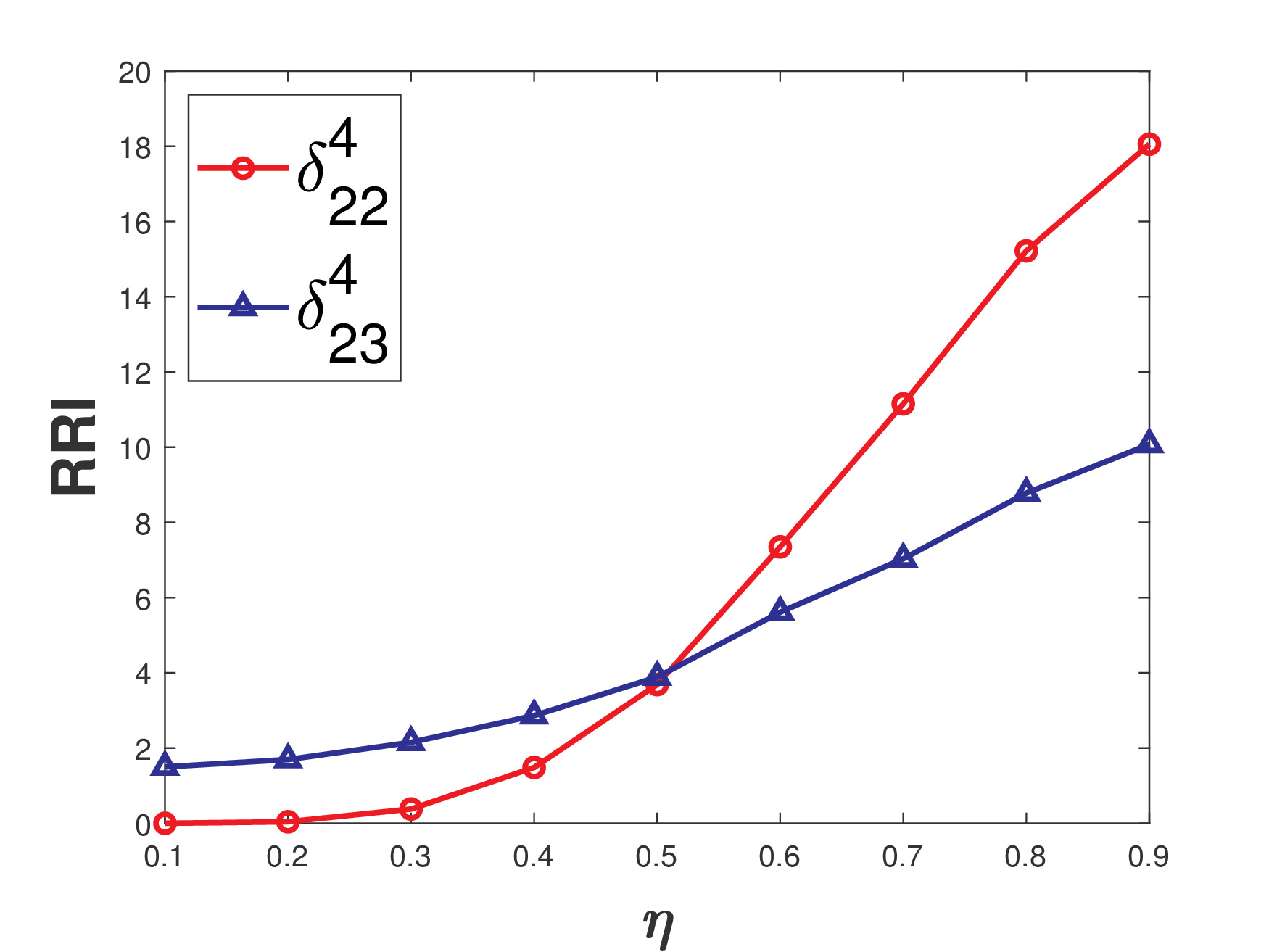}
			\caption{$a=-2$, $(p_1,p_2)=(5,5)$, $(\mu_1,\mu_2)=(0,0)$}
		\end{subfigure}
		\hfill
		\begin{subfigure}{0.32\textwidth}
			\centering
			\includegraphics[height=3.9cm, width=\textwidth]{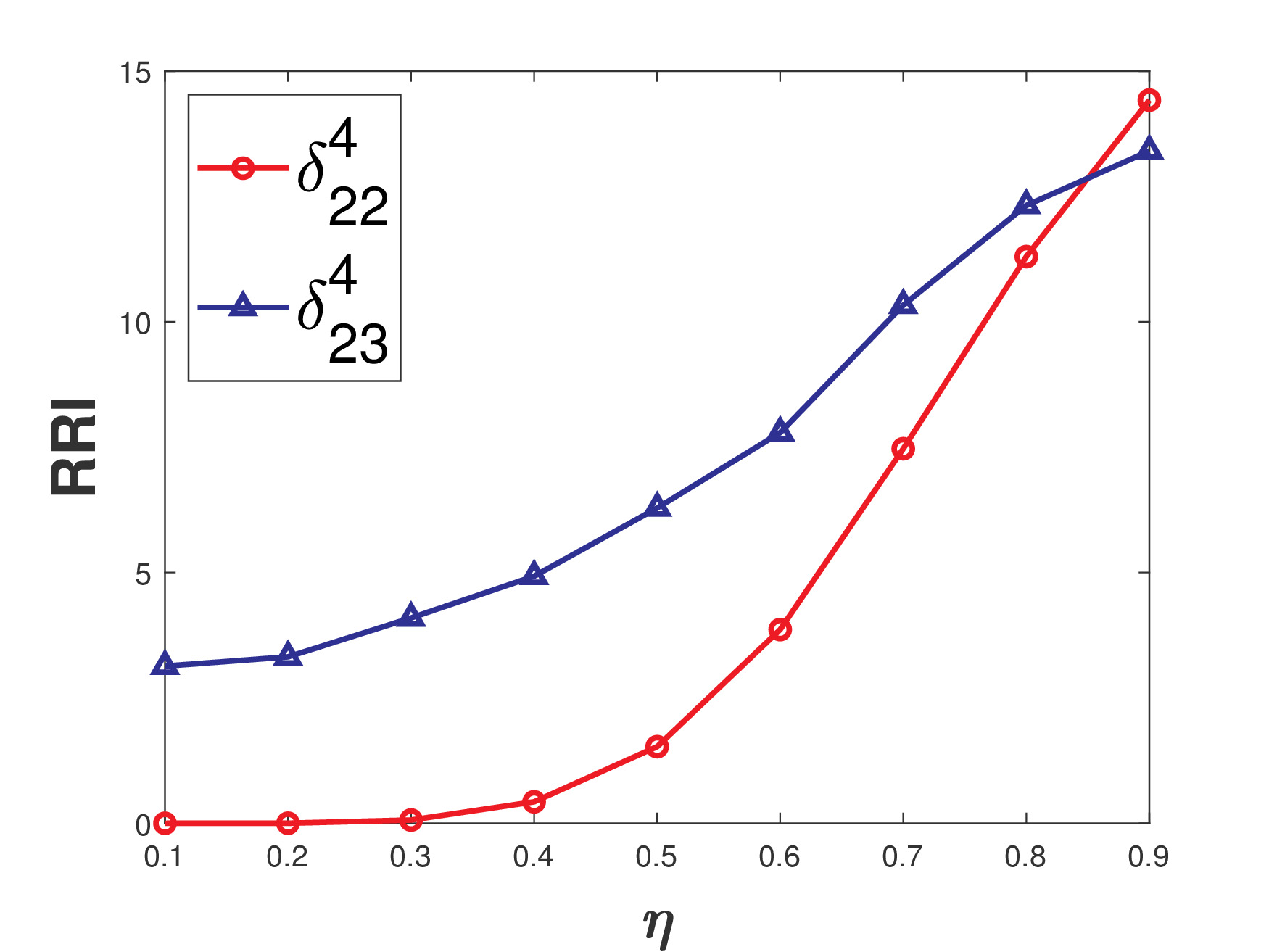}
			\caption{$a=-2$, $(p_1,p_2)=(6,8)$, $(\mu_1,\mu_2)=(0,0.3)$}
		\end{subfigure}
		\hfill
		\begin{subfigure}{0.32\textwidth}
			\centering
			\includegraphics[height=3.9cm, width=\textwidth]{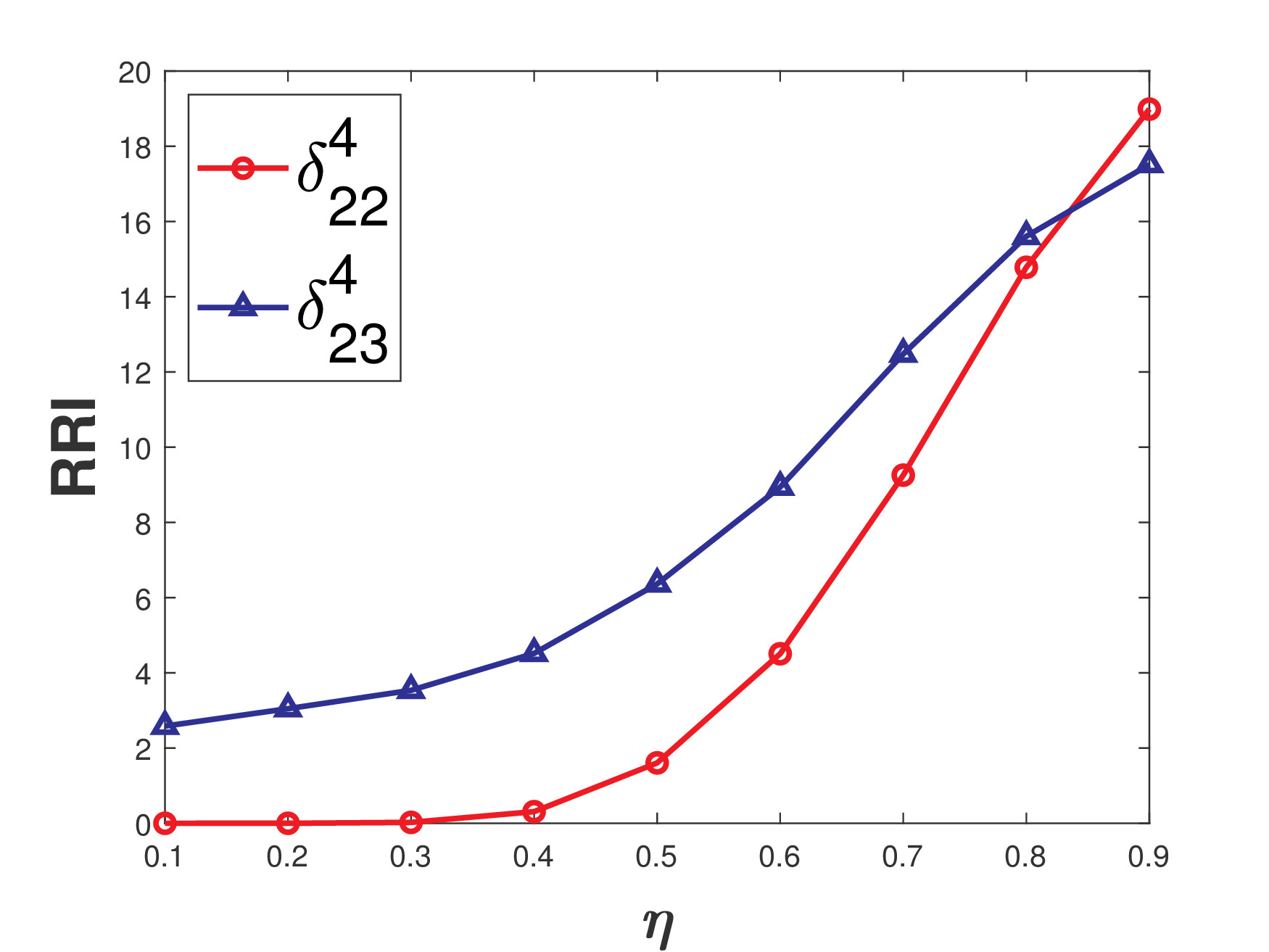}
			\caption{$a=-2$, $(p_1,p_2)=(9,12)$, $(\mu_1,\mu_2)=(0.15,0.25)$}
		\end{subfigure}
	\end{subfigure}
\caption{RRI of different estimators with respect to BAEE for $\sigma_2^2$ under $L_4(t)$.}\label{fig8}
\end{figure}

	In this section, we will compare the risk performance of the improved estimators proposed in the previous sections with respect to various scale-invariant loss functions. For this purpose we have generated $60000$ random samples from two normal populations $N(\mu_1,\sigma_1^2)$ and $N(\mu_2,\sigma_2^2)$ for various values of $(\mu_1,\mu_2)$ and $(\sigma_1,\sigma_2)$. Observed that the risk of estimators depends on the parameters $\sigma_1$ and $\sigma_2$ through $\eta=\sigma_1/\sigma_2$. The performance measure of the improved estimators has been studied using relative risk improvement (RRI) with respect to BAEE. 
The relative risk improvement of the estimators $\delta$ with the respect to $\delta_{0}$ is defined as 
$$\mbox{RRI}(\delta) = \frac{\mbox{Risk}(\delta_0)-\mbox{Risk}(\delta)}{\mbox{Risk}(\delta_{0})} \times 100.$$
In the simulation study, we have considered the case $k=2$. We have plotted the RRI of the improved estimator of $\sigma_1^2$ in Figure \ref{fig1}, \ref{fig2}, \ref{fig3} and \ref{fig4} under the loss functions $L_1$, $L_2$, $L_3$ and $L_4$. We now present the following observations from Figure \ref{fig1}, which corresponds to the quadratic loss function $L_1(t)$. 
\begin{enumerate}
	\item[(i)] The RRI of $\delta^1_{11}$ and $\delta^1_{14}$ are increasing functions of $\eta$ but $\delta_{\phi_{*}^1}$ is not monotone in $\eta$. The improvement region of $\delta^1_{11}$ larger than $\delta^1_{14}$ for all values of $\eta$.  
	 The risk performance $\delta^1_{11}$ and $\delta^1_{14}$ is better when $(\mu_1,\mu_2)$ is closed to $(0,0)$.
	
	\item [(iv)] The RRI of $\delta_{\phi_{*}^1}$ increases when $\eta\le 0.6$ (approximately) and decreases otherwise. However, $\delta_{\phi_{*}^1}$ achieve the highest risk improvement region compared to $\delta^1_{11}$ and $\delta^1_{14}$ for all values of $\eta$.
	
	\item [(v)] The risk performance of $\delta_{\phi_{*}^1}$ is better than $\delta^1_{11}$ and $\delta^1_{14}$ in the region $0.1 \leq \eta \leq 0.73$ (approximately) and under performed when $\eta \geq 0.73$ (approximately). Furthermore $\delta_{11}^1$ dominates $\delta_{\phi_{*}^1}$ as well as $\delta_{14}^1$ when $\eta\geq 0.73$ (approximately).
	
	
	\item [(vi)] The RRI $\delta^1_{12}$ and $\delta^1_{13}$ are increasing function of $\eta$. The improvement region for these estimators becomes smaller when sample sizes are increased and the value of $(\mu_1,\mu_2)$ deviates from $(0,0)$. However the risk performance of $\delta_{12}^1$ is better than $\delta_{13}^1$ for any values of $\eta$. Furthermore, in the Figure \ref{fig4}, under the loss function $L_4(t)$, the estimator $\delta_{13}^4$ is not an increasing function of $\eta$. 
\end{enumerate}
We observe similar behaviour in the simulation results for the entropy loss function $L_2(t)$, the symmetric loss function $L_3(t)$ and the Linex loss function $L_4(t)$. For the Linex loss function, we plotted the graphs for different values of $a = -2, -1, 1, 2$. However, the Figure \ref{fig4} shows only the case for $a=-2$, while the remaining plots are provided in the supplementary material. The RRI of the improved estimators with the respect to BAEE for the $\sigma_2^2$ under $L_1(t)$, $L_2(t)$, $L_3(t)$ and $L_4(t)$ is shown in Figure \ref{fig5}, \ref{fig6}, \ref{fig7} and \ref{fig8} respectively. We now discuss the following observation for the quadratic loss function $L_1(t)$ based on Figure \ref{fig5}.
\begin{enumerate}
	\item[(i)] The relative risk improvements of $\delta^1_{21}$ and $\delta^1_{22}$ are increasing function $\eta$. However, $\delta_{\psi_{*}^1}$ is not  strictly increasing in $\eta$; it increasing when $\eta$ lies between $0.1$ to $0.8$ (approximately) and decreasing for $\eta \ge 0.8$ (approximately).
	The improvement region of $\delta_{21}^1$ is greater than that of $\delta_{24}^1$ for all values of $\eta$. However, $\delta_{\psi_{*}^1}$ shows the highest improvement region compared to the $\delta_{2}^1$ and $\delta_{22}^1$.
	\item[(iii)] The RRI of $\delta_{22}^1$ is an increasing function of $\eta$, whereas $\delta_{23}^1$ is not necessarily monotone in $\eta$ (see Figure \ref{fig5}).
	
	\item[(iv)] The improvement regions of $\delta_{22}^1$ and $\delta_{22}^2$ become smaller as the sample size increases or as the parameter values $(\mu_1,\mu_2)$ deviate further from $(0,0)$ (An opposite behavior can be observed under the loss function $L_4(t)$ in Figure \ref{fig8}).
	
	\item[(v)] When $(\mu_1,\mu_2)$ are sufficiently close to $(0,0)$, the risk performance of $\delta_{\psi_{*}^1}$ is significantly better than that of the other estimators within the domain $0.1\leq \eta \leq 0.8$ (approximately). 
	
\end{enumerate}
We observe similar patterns under the entropy loss function $L_2(t)$, symmetric loss function $L_3(t)$ and linex loss function $L_4(t)$.

In conclusion, overall  performance of the estimators $\delta_{\phi_{*}^1}$, $\delta_{\phi_{*}^2}$, $\delta_{\phi_{*}^3}$ and $\delta_{\phi_{*}^4}$ are better than the other competing estimators for estimating $\sigma_1^k$ and similarly for $\sigma_2^k$. Therefore, we recommend these estimators for use in real-life applications.

\captionsetup[subfigure]{justification=centering, singlelinecheck=on, font=small}
\captionsetup[subfigure]{justification=centering, singlelinecheck=on, font=small}
\captionsetup[subfigure]{justification=centering, singlelinecheck=on, font=small}
\captionsetup[subfigure]{justification=centering, singlelinecheck=on, font=small}
\captionsetup[subfigure]{justification=centering, singlelinecheck=on, font=small}
\captionsetup[subfigure]{justification=centering, singlelinecheck=on, font=small}
\captionsetup[subfigure]{justification=centering, singlelinecheck=on, font=small}
\captionsetup[subfigure]{justification=centering, singlelinecheck=on, font=small}
						
\clearpage
\section{Data analysis} \label{sec5}
This section presents a real-life data analysis to illustrate our findings. In particular, we obtain the estimates of $\sigma_i^4$ and $\sigma_i^5$ for $i=1,2$. We have taken the data set form from https:// data.opencity.in/dataset/ bengaluru-rainfall and https://data.opencity.in/dataset/hyderabad-rainfall-data. This data reports the total annual rainfall (in mm) in Bengaluru and Hyderabad from $1985$ to $2000$, respectively. The datasets are given below.\\

\noindent \textbf{Bengaluru (Data-I):}$~$
$634,$ $1145.1,$ $ 798.6,$ $1221,$  $905.1,$ $613.1,$ $1350.5,$ $826.3,$ $1069,$ $587.2,$  $1068.4,$ $1172.9,$ 1229.8, 1431.8,  1014, 1193.9\\

\noindent \textbf{Hyderabad (Data-II):}$~$
550.280, 682.652,  978.866, 828.710, 867.701, 964.704, 836.222, 638.196, 673.357, 792.315, 1166.311,
953.552, 781.910, 879.267, 535.661, 959.343.
\\

Using the Kolmogorov-Smirnov test at a significance level of $0.05$, we find that both datasets satisfy the normality assumption with p-values of $0.8654$ and $0.9764$ for the first and second datasets, respectively. Here we assume that $\sigma_1 \le \sigma_2$. Based on these data, the summarized data are as follows: $p_1=16$, $p_2=16$, $X_1=1016.2937$, $X_2=818.0654$, $S_1=1038675.0494$, and $S_2=438664.9655$, Where $X_1$ and $X_2$ are the sample mean of the Data-I and and Data-II respectively. The quantities $S_1=\sum_{i=1}^{16}\left(X_{1i}-X_1\right)^2$ and $S_2=\sum_{i=1}^{16}\left(X_{2i}-X_2\right)^2$ are total sum of squares for Data-I and Data-II, where $X_{1i}$ and $X_{2i}$ denote individual observations from Data-I and Data-II respectively. Using these statistics, we have computed the values of the estimators several of $\sigma_1^2$, $\sigma_2^2$, $\sigma_1^4$, and $\sigma_2^4$, and the values of the estimators are tabulated in Tables \ref{table1}, \ref{table2}, \ref{table3}, and \ref{table4} respectively.

\begin{table}[h!]
	\centering
	\caption{Values of the estimators of $\sigma_1^2$.}\label{table1}
\begin{tabular}{@{}cccccc@{}}
			\toprule
			& $\delta_{01}$ & $\delta_{11}$ & $\delta_{12}$ &  $\delta_{13}$ & $\delta_{\phi_{*}}$ \\
			\hhline{======}
			$L_1(t)~~~~~~~~~~~~~$  & $6.10099\times10^4$& $4.6167\times10^4$ & $6.1099\times10^4$ & $6.10999\times10^4$ &  $4.2395\times10^4$ \\
			$L_2(t)~~~~~~~~~~~~~$  &$6.9245\times10^4$  &  $4.9245\times10^4$ & $6.9245\times10^4$  & $6.9245\times10^4$  &  $4.6295\times10^4$\\
			$L_3(t)~~~~~~~~~~~~~$  & $7.4381\times10^{4}$ & $5.0973\times10^{4}$ & $7.4381\times10^{4}$ & $7.4381\times10^{4}$ & $4.7976\times10^{4}$\\ 
			$L_4(t)\ (a=-2)$  & $6.8885\times10^{4}$ & $4.9176\times10^{4}$ & $6.8885\times10^{4}$ & $6.8885\times10^{4}$ & -\\
			$L_4(t)\ (a=-1)$  & $6.4838\times10^{4}$ & $4.7640\times10^{4}$ & $6.4838\times10^{4}$ & $6.4838\times10^{4}$ & -\\
			$L_4(t)\ (a=1)~~$  & $5.7641\times10^{4}$ & $4.4754\times10^{4}$ & $5.7641\times10^{4}$ & $5.7641\times10^{4}$ & -\\
			$L_4(t)\ (a=2)~~$  & $5.4443\times10^{4}$ & $4.3398\times10^{4}$ & $5.4443\times10^{4}$ & $5.4443\times10^{4}$ & -\\
			\bottomrule
	\end{tabular}
\end{table}
\clearpage
\vspace{0.3cm}
\begin{table}[h!]
	\centering
		\caption{Values of the estimators of  $\sigma_2^2$.}\label{table2}
\begin{tabular}{@{}cccccc@{}}
			\toprule
			& $\delta_{02}$ & $\delta_{21}$ & $\delta_{22}$ &  $\delta_{23}$ &  $\delta_{\psi_{*}}$ \\
			\hhline{======}
			$L_1(t)~~~~~~~~~~~~~$  &$2.5804\times10^{4}$&  $4.6167\times10^{4}$ & $2.5804\times10^{4}$ & $2.5804\times10^{4}$ & $5.3220\times10^{4}$ \\
			$L_2(t)~~~~~~~~~~~~~$ &$2.9244\times10^{4}$ & $4.9245\times10^{4}$ & $2.9244\times10^{4}$ & $2.9244\times10^{4}$ & $5.7994\times10^{4}$\\ 
			$L_3(t)~~~~~~~~~~~~~$ & $3.1413\times10^{4}$ & $5.0973\times10^{4}$ & $3.1413\times10^{4}$ & $3.1413\times10^{4}$ & $6.0887\times10^{4}$\\ 
			$L_4(t)\ (a=-2)$ & $2.9092\times10^{4}$ & $4.9176\times10^{4}$ & $2.9092\times10^{4}$ & $2.9092\times10^{4}$ & -\\ 
			$L_4(t)\ (a=-1)$  & $2.7383\times10^{4}$ & $4.7640\times10^{4}$ & $2.7383\times10^{4}$ & $2.7383\times10^{4}$ & -\\
			$L_4(t)\ (a=1)~~$  & $2.4344\times10^{4}$ & $4.4754\times10^{4}$ & $2.4344\times10^{4}$ & $2.4344\times10^{4}$ & -\\
			$L_4(t)\ (a=2)~~$  & $2.2993\times10^{4}$ & $4.3398\times10^{4}$ & $2.2993\times10^{4}$ & $2.2993\times10^{4}$ & -\\
			\bottomrule
	\end{tabular}
\end{table}

\begin{table}[h!]
 \centering
 	\caption{Improved estimator values for $\sigma_1^4$.}\label{table3}
\begin{tabular}{@{}cccccc@{}}
	\toprule
 & $\delta_{01}$ & $\delta_{11}$ & $\delta_{12}$ &  $\delta_{13}$ & $\delta_{\phi_{*}}$ \\
\hhline{======}
		$L_1(t)$  & $2.7039\times10^9$& $1.7831\times10^9$ & $2.7039\times10^9$ & $2.7039\times10^9$ &  $1.5559\times10^9$ \\
		$L_2(t)$  &$4.2308\times10^9$  &  $2.2735\times10^9$ & $4.2308\times10^9$  & $4.2308\times10^9$  &  $2.0043\times10^9$\\
		$L_3(t)$  & $5.6496\times10^{9}$ & $2.6107\times10^{9}$ & $5.6496\times10^{9}$ & $5.6496\times10^{9}$ & $2.3124\times10^{9}$\\ 
		\bottomrule
\end{tabular}
\end{table}

\vspace{0.3cm}
\begin{table}[h!]
\centering
	\caption{Improved estimator values for $\sigma_2^4$.}\label{table4}
	\begin{tabular}{@{}cccccc@{}}
		 		\toprule
		 		& $\delta_{02}$ & $\delta_{21}$ & $\delta_{22}$ &  $\delta_{23}$ & $\delta_{\psi_{*}}$ \\
		 		\hhline{======}
		 		$L_1(t)$  & $0.4823\times10^9$& $1.7831\times10^9$ & $0.4823\times10^9$ & $0.4823\times10^9$ & $2.2617\times10^9$ \\
		 		$L_2(t)$  &$0.7546\times10^9$& $2.2735\times10^9$ & $0.7546\times10^9$ & $0.7546\times10^9$ & $3.0865\times10^9$\\ 		
		 		$L_3(t)$  & $1.0077\times10^9$ & $2.6107\times10^9$ & $1.0077\times10^9$ & $1.0077\times10^9$ & $3.7581\times10^9$\\ 
		 		\bottomrule
		 	\end{tabular}
 
\end{table}


 \section{Conclusions}   
 In this manuscript, we consider the problem of estimating the positive power of the ordered variance of two normal populations when means satisfy certain restrictions. The estimation problem has been studied with respect to a general bowl-shaped scale-invariant loss function. We propose sufficient conditions under which we obtain estimators dominating the BAEE. We have obtained various \cite{stein1964}-type improved estimators that improve upon the BAEE. Further, a class of improved estimators has been presented using the IERD approach of \cite{kubokawa1994double}. We observed that the boundary estimator of this class is the \cite{brewster1974improving}-type estimator. Moreover, we showed that the \cite{brewster1974improving}-type improved estimator is a generalized Bayes estimator. We have obtained the expression of the improved estimator for quadratic, entropy, symmetric loss, and Linex to demonstrate an immediate application.
  Further, a simulation study is conducted to compare the risk performance of the proposed estimators. For $k=2$, we evaluated the performance of various improved estimators of $\sigma_1^2$ and $\sigma_2^2$ under quadratic, entropy, symmetric, and Linex losses. The \cite{brewster1974improving}-type estimators 
 perform better than others when $\eta<0.7$ approximately and $(\mu_1,\mu_2)$ are close to zero. However, for  $\eta>0.7$, Stein-type estimators perform better. Finally a data analysis is given. In the data analysis we have obtained the values of the estimators of $\sigma_i^2$ and $\sigma_i^4$, $i=1,2$. Furthermore, for the Linex loss function, we conducted the analysis for different values of the parameter $a$, specifically $a = -2, -1, 1, 2$.
 
 \section{Appendix}
 \begin{lemma}
 	The function $f(x;r)=\frac{\Gamma(x)}{\Gamma(x+r)}$ is strictly decreasing (increasing) in $x$ for all $x>0$ and fixed $r>0\ (r<0)$.
 \end{lemma}						
 \noindent\textbf{Proof.} Let $g(x)= \log f(x;r) = \log\left(\Gamma(x)\right) - \log\left(\Gamma(x+r)\right)$. The derivative of $\log\left(\Gamma(x)\right)$ is the digamma function $\psi(x)$. So, $g'(x)=\psi(x)-\psi(x+r)$. Now for all $x>0$, and any fixed $r>0$, the digamma function satisfies : $\psi(x+r)>\psi(x) \implies \psi(x)-\psi(x+r)<0$. So $g'(x)=\psi(x)-\psi(x+r)<0$ for all $x>0$, $r>0$. Hence $g(x)$ is strictly decreasing i.e., $f(x;r)=\frac{\Gamma(x)}{\Gamma(x+r)}$ is strictly decreasing.
						
\section*{Funding}
 Lakshmi Kanta Patra thanks the Science and Engineering Research Board, India for providing financial support to carry out this research  with project number MTR/2023/000229.
\section*{Disclosure statement}
No potential conflict of interest was reported by the authors.
\section*{Acknowledgments}
The authors would like to thank the anonymous referees for their comments, which have led to significant improvements in the paper. 
\bibliography{bib_normal}

@article{kubokawa1994double,
	title={Double shrinkage estimation of ratio of scale parameters},
	author={Kubokawa, Tatsuya},
	journal={Annals of the Institute of Statistical Mathematics},
	volume={46},
	pages={95--116},
	year={1994},
	publisher={Springer}
}

@article{brewster1974improving,
  title={Improving on equivariant estimators},
  author={Brewster, J F and Zidek, JV},
  journal={The Annals of Statistics},
  volume={2},
  number={1},
  pages={21--38},
  year={1974},
  publisher={JSTOR}
}

@article{kubokawa1994unified,
  title={A unified approach to improving equivariant estimators},
  author={Kubokawa, Tatsuya},
  journal={The Annals of Statistics},
  volume={22},
  number={1},
  pages={290--299},
  year={1994},
  publisher={JSTOR}
}

@article{stein1964,
  title={Inadmissibility of the usual estimator for the variance of a normal distribution with unknown mean},
  author={Stein, Charles},
  journal={Annals of the Institute of Statistical Mathematics},
  volume={16},
  number={1},
  pages={155--160},
  year={1964},
  publisher={Springer}
}

@article{lehmann2005testing,
  title={Testing Statistical Hypotheses, Springer, New York},
  author={Lehmann, EL and Romano, JP},
  journal={MR2135927},
  year={2005}
}

@article{patra2021componentwise,
  title={Componentwise estimation of ordered scale parameters of two exponential distributions under a general class of loss function},
  author={Patra, Lakshmi Kanta and Kumar, Somesh and Petropoulos, Constantinos},
  journal={Statistics},
  volume={55},
  number={3},
  pages={595--617},
  year={2021},
  publisher={Taylor \& Francis}
}

@article{misra1994estimation,
  title={Estimation of ordered location parameters: the exponential distribution},
  author={Misra, Neeraj and Singh, Harshinder},
  journal={Statistics: A Journal of Theoretical and Applied Statistics,},
  volume={25},
  number={3},
  pages={239--249},
  year={1994},
  publisher={Taylor \& Francis}
}

@article{misra2002smooth,
  title={Smooth estimators for estimating order restricted scale parameters of two gamma distributions},
  author={Misra, Neeraj and Choudhary, PK and Dhariyal, ID and Kundu, D},
  journal={Metrika},
  volume={56},
  pages={143--161},
  year={2002},
  publisher={Springer}
}

@article{tripathy2013estimating,
  title={Estimating common standard deviation of two normal populations with ordered means},
  author={Tripathy, Manas Ranjan and Kumar, Somesh and Pal, Nabendu},
  journal={Statistical Methods \& Applications},
  volume={22},
  number={3},
  pages={305--318},
  year={2013},
  publisher={Springer}
}

@article{vijayasree1995componentwise,
	title={Componentwise estimation of ordered parameters of $k (\geq 2)$ exponential populations},
	author={Vijayasree, G and Misra, Neeraj and Singh, Harshinder},
	journal={Annals of the Institute of Statistical Mathematics},
	volume={47},
	pages={287--307},
	year={1995},
	publisher={Springer}
}

@article{kumar1988simultaneous,
	title={Simultaneous estimation of ordered parameters},
	author={Kumar, Somesh and Sharma, Divakar},
	journal={Communications in Statistics-Theory and Methods},
	volume={17},
	number={12},
	pages={4315--4336},
	year={1988},
	publisher={Taylor \& Francis}
}

@article{petropoulos2017estimation,
	title={Estimation of the order restricted scale parameters for two populations from the Lomax distribution},
	author={Petropoulos, Constantinos},
	journal={Metrika},
	volume={80},
	number={4},
	pages={483--502},
	year={2017},
	publisher={Springer}
}

@book {MR0326887,
	AUTHOR = {Barlow, R. E. and Bartholomew, D. J. and Bremner, J. M. and
	Brunk, H. D.},
	TITLE = {Statistical inference under order restrictions. {T}he theory
	and application of isotonic regression},
	NOTE = {},
	PUBLISHER = {John Wiley \& Sons},
	YEAR = {1972},
	PAGES = {xii+388},
	MRCLASS = {62FXX (62GXX 62JXX)},
	MRNUMBER = {0326887},
	MRREVIEWER = {Lionel Weiss},
}

@book {MR961262,
	AUTHOR = {Robertson, Tim and Wright, F. T. and Dykstra, R. L.},
	TITLE = {Order restricted statistical inference},
	SERIES = {},
	PUBLISHER = {John Wiley \& Sons},
	YEAR = {1988},
	PAGES = {xx+521},
	ISBN = {0-471-91787-7},
	MRCLASS = {62-02 (62E15 62G30 62J99)},
	MRNUMBER = {961262},
	MRREVIEWER = {Roy T. Saint Laurent},
}

@book {MR2265239,
	AUTHOR = {van Eeden, Constance},
	TITLE = {Restricted parameter space estimation problems. Admissibility and minimaxity properties},
	VOLUME = {188},
	SERIES = {Lecture Notes in Statistics},
	PUBLISHER = {Springer, New York},
	YEAR = {2006},
	PAGES = {x+167},
	ISBN = {978-0-387-33747-0; 0-387-33747-4},
	MRCLASS = {62-02 (62C15 62C20 62F10 62F30)},
	MRNUMBER = {2265239},
	MRREVIEWER = {Jon Stene},
	}

@article{oono2005estimation,
		title={Estimation of two order restricted normal means with unknown and possibly unequal variances},
		author={Oono, Youhei and Shinozaki, Nobuo},
		journal={Journal of statistical planning and inference},
		volume={131},
		number={2},
		pages={349--363},
		year={2005},
		publisher={Elsevier}
	}

@article{oono2006class,
		title={On a class of improved estimators of variance and estimation under order restriction},
		author={Oono, Youhei and Shinozaki, Nobuo},
		journal={Journal of statistical planning and inference},
		volume={136},
		number={8},
		pages={2584--2605},
		year={2006},
		publisher={Elsevier}
	}

@article{jana2015estimation,
		title={Estimation of ordered scale parameters of two exponential distributions with a common guarantee time},
		author={Jana, N and Kumar, S},
		journal={Mathematical Methods of Statistics},
		volume={24},
		pages={122--134},
		year={2015},
		publisher={Springer}
	}

@article{garg2023componentwise,
		title={Componentwise equivariant estimation of order restricted location and scale parameters in bivariate models: a unified study},
		author={Garg, Naresh and Misra, Neeraj},
		journal={Brazilian Journal of Probability and Statistics},
		volume={37},
		number={1},
		pages={101--123},
		year={2023},
		publisher={Brazilian Statistical Association}
	}

@article{misra1997estimation,
		title={On estimation of the common mean of $k ( \ge2)$ normal populations with order restricted variances},
		author={Misra, Neeraj and van der Meulen, Edward C},
		journal={Statistics \& probability letters},
		volume={36},
		number={3},
		pages={261--267},
		year={1997},
		publisher={Elsevier}
	}

@article{chang2002comparison,
		title={A comparison of restricted and unrestricted estimators in estimating linear functions of ordered scale parameters of two gamma distributions},
		author={Chang, Yuan-Tsung and Shinozaki, Nobuo},
		journal={Annals of the Institute of Statistical Mathematics},
		volume={54},
		pages={848--860},
		year={2002},
		publisher={Springer}
	}

@article{chang2017estimation,
		title={Estimation of two ordered normal means when a covariance matrix is known},
		author={Chang, Yuan-Tsung and Fukuda, Kazufumi and Shinozaki, Nobuo},
		journal={Statistics},
		volume={51},
		number={5},
		pages={1095--1104},
		year={2017},
		publisher={Taylor \& Francis}
	}

@article{ma2013estimation,
		title={Estimation of order-restricted means of two normal populations under the LINEX loss function},
		author={Ma, Tiefeng and Liu, Shuangzhe},
		journal={Metrika},
		volume={76},
		pages={409--425},
		year={2013},
		publisher={Springer}
	}

@article{jana2022estimation,
		title={Estimation of order restricted standard deviations of normal populations with a common mean},
		author={Jana, Nabakumar and Chakraborty, Ankur},
		journal={Statistics},
		volume={56},
		number={4},
		pages={867--890},
		year={2022},
		publisher={Taylor \& Francis}
	}
\end{document}